\documentclass{article}
\usepackage{amssymb, amsmath}
\usepackage{hyperref}

\def\Ric{\mathop{\rm Ric}}

\def\Hess{\mathop{\rm Hess}}

\def\dist{\mathop{\rm dist}}

\def\Riem{\mathop{\rm Rm}}

\def\Diam{\mathop{\rm Diam}}

\def\Vol{\mathop{\rm Vol}}
\def\supp{\mathop{\rm supp}}
\def\div{\mathop{\rm div}}

\newcommand{\qed}{\hfill$\Box$}
\newtheorem{theorem}{Theorem}[section]
\newtheorem{proposition}[theorem]{Proposition}
\newtheorem{lemma}[theorem]{Lemma}
\newtheorem{technicallemma}[theorem]{Technical Lemma}
\newtheorem{corollary}[theorem]{Corollary}

\newtheorem{defi}[theorem]{Definition}

\newtheorem{prob}[theorem]{Problem}
\newtheorem{hypothesis}[theorem]{Hypothesis}

\date{\small\it\today}
\title{Moduli Spaces of critical Riemanian Metrics with $L^{n\over 2} $ norm
curvature bounds}
\author{Xiuxiong Chen and Brian Weber}

\begin{document}

\maketitle
\begin{abstract}
We consider the moduli space of the extremal K\"ahler metrics on
compact manifolds.  We show that under the conditions of two-sided
total volume bounds, $L^{\frac{n}{2}}$-norm bounds on $\Riem$, and
Sobolev constant bounds, this Moduli space can be compactified by
including (reduced) orbifolds with finitely many singularities. Most
of our results go through for certain other classes of critical
Riemannian metrics.
\end{abstract}
\tableofcontents

\section{Introduction}\label{Introduction}

A K\"ahler metric is called extremal if the complex gradient of its
scalar curvature is a holomorphic vector field.  This includes the
more famous K\"ahler Einstein metrics and constant scalar curvature
K\"ahler (cscK) metrics as special cases, though one would like to
understand the structure of extremal metrics as well. In this note,
we propose to study the weak compactness of extremal K\"ahler
metrics in a bounded family of K\"ahler classes together with bounds
on the $L^{n\over 2}$ norm of Riemannian curvature and on the
Sobolev constants. The extremal K\"ahler metric equation is
naturally a 6th order equation on K\"ahler potential, and its
compactness properties are difficult to study directly. We
essentially decompose the extremal condition into three
inter-related second order equations as below:
\begin{eqnarray}
&&\triangle\Riem \;=\; \nabla^2\Ric  \,+\, \Riem*\Riem \label{RiemEq}\\
&&\triangle\Ric \;=\; \Ric * \Riem \,+\, \nabla X \label{RicEq}\\
&&\triangle X \;=\; \Riem * X. \label{XEq}
\end{eqnarray}
The $``*"$ stands for tensor contraction between two multi-index
tensors (more elaboration on this below) and $X$ is a vector field
related to the critical Riemannian metric\footnote{For the extremal
K\"ahler metrics, $X$ is the complex gradient vector field of the
scalar curvature function. }.  A large class of critical metrics
satisfy these three coupled equations, for instance CSC Bach-flat
metrics, harmonic curvature metrics, and Einstein metrics, all of
which have been studied before. Below we show that another class of
metrics, the extremal K\"ahler metrics, also satisfy these
equations.

More specifically we study the weak compactness of the space
$\mathcal{M} = \mathcal{M}(n,C_S,\Lambda,\nu,\delta)$ of critical
metrics  (where $X$ is non-trivial) that satisfy\footnote{For
complex surfaces, the only assumption is the Sobolev constant. The
others are either {\it a priori} or can be derived from {\it a
priori} constraints. Moreover, there is a large open set of K\"ahler
classes where also the Sobolev constant is {\it a priori} bounded
for the extremal representatives, c.f. Section 2.5. }
\begin{itemize}
\item[{\it i})] energies are bounded: $\int_M |\Riem|^{\frac{n}{2}}  \le
\Lambda$
\item[{\it ii})] volumes are bounded from below: $\Vol M \ge \nu$, and
\item[{\it iii})] diameters are bounded from above: $\dist_M(x,y) \le
\delta$, all $x,y\in M$.
\item[{\it iv})] the Sobolev constant $C_g$ on $(M,g)$ has a uniform  bound,
$C_g \le C_S$.
\end{itemize}
The Sobolev inequality referred to here controls the embedding
$W^{1,2}\hookrightarrow{L}^{\frac{2n}{n-2}}$, and usually takes the
form
\begin{eqnarray*}
\left(\int\phi^{2\gamma}\right)^{1/\gamma}\le{C_g}\int|\nabla\phi|^2
\,+\,\frac{A}{\left(\Vol{M}\right)^{2/n}}\int\phi^2,
\label{SobolevStatement}
\end{eqnarray*}
where $\gamma=\frac{n}{n-2}$ and $\phi\in{C}^1$. In fact one often
takes $max\,(C_g,A)$ to be the Sobolev constant. We require the
simplified form of the inequality,
\begin{eqnarray*}
\left(\int\phi^{2\gamma}\right)^{1/\gamma}&\le&C_g\int|\nabla\phi|^2.
\end{eqnarray*}
If one assumes $\Vol(\supp\phi)$ is smaller than $A/2\sqrt{\nu}$,
then we can use it in this form. In Section
\ref{SubsectionUniformSobolevBound} we show how sometimes $A$ and
$C_S$ are automatically controlled.

In this paper, we study the weak compactness in all dimensions of
our ``critical metrics'', which satisfy conditions ({\it i})-({\it
iv}) above. There is a substantial body of prior compactness results
which we build on. The case of CSC Bach-flat, harmonic curvature,
and CSC K\"ahler metrics was considered in \cite{TV1}, \cite{TV2}.
Recent work of Anderson's \cite{And2}, \cite{And3} elaborates on
this theme. These works in turn can be traced back to work of M.
Anderson \cite{And1}, G. Tian \cite{Tia1}, and Bando-Kasue-Nakajima
\cite{BKN} on the moduli space of Einstein metrics on four
dimensional manifolds with $L^2$ norm curvature bound. These in turn
were natural extensions of earlier work of J. Cheeger \cite{Che} and
later M. Gromov \cite{Gro}, which explored geometric and topological
control on manifolds with various pointwise bounds on curvature.
Readers are encouraged to read \cite{CT} for more references.

To analyze the inter-play of the three coupled equations, one must
obtain some {\it a priori} bounds on the $X = \nabla^{(1,0)} R$ ($R$
indicates scalar curvature throughout). Without using the assumption
of Sobolev constant bound, we derive an $L^2$ norm bound on $X$ and
$\nabla X$ in all dimensions (cf. Lemma \ref{L2X}). This is
important for both geometrical and analytical reasons. Analytically,
this $W^{2,2}$ bound on $R$, together with a bound on
$\|\Riem\|_{n\over 2}$, serves as our starting point for a weak
compactness argument on the moduli space of extremal metrics.
Geometrically, the $L^\infty$ bound on scalar curvature (likewise
the $L^2$ bound on $X$) is a consequence of the scalar curvature map
(from complex structure) being a moment map (if interpreted
correctly). It is more difficult to to understand what $\int_M
|\nabla \bar \nabla R|^2$ represents geometrically however. A
natural question is whether all $W^{k,2}$ norms of the scalar
curvature function are {\it a priori} bounded.


Perhaps the main technical theorem we prove is the usual
$\epsilon$-regularity
\begin{theorem} (cf. Theorem \ref{LocalCurvatureBounds})
Assume $g$ is a critical metric on a Riemannian manifold. When
$a>\frac{n}{2}$ and $q\in\{0,1,\dots\}$, there exists $\epsilon_0 =
\epsilon_0(C_S,a,q,n)$ and $C=C(C_S,a,q,n)$ so that
$$
\int_{B(o,r)}|\Riem|^{\frac{n}{2}} \;\le\; \epsilon_0
$$
implies
\begin{eqnarray}
&&\left(\int_{B(o,r/2)}|\nabla^q\,X|^a\right)^\frac1a
\;\le\;Cr^{-q-3+\frac{n}{a}}\left(\int_{B(o,r)}|R|^{\frac{n}{2}}
\right)^\frac2n\\
&&\left(\int_{B(o,r/2)}|\nabla^q\Ric|^a\right)^\frac1a
\;\le\;Cr^{-q-2+\frac{n}{a}}\left(\int_{B(o,r)}|\Ric|^{\frac{n}{2}}
\right)^\frac2n\\
&&\left(\int_{B(o,r/2)}|\nabla^q\Riem|^a\right)^\frac1a
\;\le\;Cr^{-q-2+\frac{n}{a}}\left(\int_{B(o,r)}|\Riem|^{\frac{n}{2}}
\right)^\frac2n.
\end{eqnarray}
\end{theorem}

This is obtained by interactive use of the three equations.   From a
purely technical point of view, the case $n > 4$ is more complicated
than the case of $n=4$ (in the smooth case at least). For $n>4$, we
derive all three estimates simultaneously using an induction
argument (see appendix). The proof is lengthy and technical and we
hope it can be shortened in the future.

The main theorems we prove are:
\begin{theorem}(cf. Theorem \ref{RiemBounds})
Assume $g$ is a critical metric on a Riemannian manifold. Then there
exists an $\epsilon_0=\epsilon_0(C_S,n,p)$ and $C=C(C_S,n,p)$ so
that $\int_{B_r} |\Riem|^{\frac{n}{2}} \le \epsilon_0$ implies
\begin{eqnarray*}
\sup_{B(o,r/2)} |\nabla^p\Riem| \;\le\;Cr^{-p-2}\left(\int_{B(o,r)}
|\Riem|^{\frac{n}{2}}\right)^\frac2n.
\end{eqnarray*}
\end{theorem}
And, specializing to the case of extremal K\"ahler manifolds,
\begin{theorem}[Orbifold compactness] \label{FinalCompactnessTheorem0}
(cf. Theorem \ref{FinalCompactnessTheorem}) Assume
$\{(M_\alpha,J_\alpha,\omega_\alpha)\}$ is a family of compact
extremal K\"ahler manifolds that satisfy conditions ({\it i}) -
({\it iv}). Then a subsequence converges in the Gromov-Hausdorff
topology to a (reduced) compact extremal K\"ahler orbifold. Further,
there is a bound $C_1=C_1(\Lambda,C_S,n)$ on the number of
singularities, and a bound $C_2=C_2(C_S,n)$ on the order of any
orbifold group.
\end{theorem}
If the family does not consist of extremal metrics but their metrics
satisfy the elliptic system (\ref{RiemEq}), (\ref{RicEq}),
(\ref{XEq}) and conditions ({\it i})-({\it iv}),
this theorem still holds, except that the singularities are only of
orbifold type $C^0$, and are not necessarily reduced (meaning a
tangent cone could be a one-point union of standard cones over
various $\mathbb{S}^3/\Gamma$). There is a variety of classes of
metrics that satisfy (\ref{RiemEq}), (\ref{RicEq}), and (\ref{XEq}),
for instance the CSC Bach-flat and harmonic curvature metrics
considered in \cite{TV1}, where in fact $X=0$.

A nontrivial step in proving orbifold compactness is to prove a
uniform upper bound on the local volume ratio.  If there is a {\sl
pointwise} lower bound on Ricci curvature, then this upper bound is
automatic from the Bishop-Gromov comparison theorem. We do not
assume such curvature lower bounds, so we prove that volume growth
is uniformly bounded by generalizing a result of Tian-Viaclovsky's
\cite{TV1}, \cite{TV2} to cover our class of critical metrics in all
dimensions. In \cite{TV1} Tian-Viaclovsky proved that complete
manifolds with bounded energy, bounded Sobolev constant, and
quadratic curvature decay $|\Riem|=o(r^{-2})$ have finitely many ALE
ends and therefore a global upper bound on volume growth. This
represented a major advance; previous results had required a nearly
unusable strengthening of the curvature decay condition. In
\cite{TV2} they use this to prove uniform volume ratio bounds on
compact manifolds with certain critical metrics, without pointwise
bounds on Ricci curvature.



Recall that a specified structure, say a differential manifold
structure or a vector bundle structure, is said to exist on an
orbifold if it exists at all manifold points and, after lifting, can
be completed on any local orbifold cover. In the 4-dimensional case,
in the absence of additional rigidity, the analytic methods
presently known are only strong enough to show that the orbifold
metric is continuous (see \cite{And2}).

Showing that the completion of the orbifold metric (on a smooth
orbifold cover) is $C^\infty$ requires a way to remove apparent
point-singularities. In higher dimensions, powerful analytic
techniques, developed originally to remove singularities in
Yang-Mills instantons, suffice to remove the singularities in our
case as well (e.g. Lemma \ref{L2Lemma}, Proposition \ref{LpLemma}).
The critical case is real dimension $4$, where these analytic
techniques fail. Here one needs the geometry itself to provide
additional rigidity. We find this rigidity, in the case of extremal
K\"ahler metrics, in a partially improved Kato inequality (Lemma
\ref{ImprovedKatoLemma}), which we take advantage of using
Uhlenbeck's broken Hodge gauge technique (\cite{Uhl}, \cite{Tia1},
\cite{TV1}).

In \cite{TV1} an improved Kato inequality was shown to hold for
4-dimensional CSC Riemannian metrics with $\delta{W}^{+}=0$ in the
(sharp) form
\begin{eqnarray*}
|\nabla|E||^2\;\le\;\frac23|\nabla{E}|^2,
\end{eqnarray*}
where $E$ is the trace-free Ricci tensor. This is actually a
consequence of the theory of Kato constants developed in \cite{Bra}
and \cite{CGH}. This is sufficient for applications to K\"ahler
geometry, where for instance constant scalar curvature implies that
$W^{+}$ is constant. We are able to use a direct argument to
partially recover an improved inequality. Specifically, we get
\begin{eqnarray}
2|\nabla|\nabla{E}||^2\;\le\;\frac14|\nabla\nabla{E}|^2\,+\,|\overline
\nabla\nabla{E}|^2.
\end{eqnarray}
This does not quite give sufficient control on the Hessian of $E$;
see Proposition \ref{ImprovedKatoLemma} and its use in Proposition
\ref{ImprovedEllipticUsingKato}. As a result, the removable
singularity theorem becomes correspondingly more complicated, and we
need to utilize Uhlenbeck's technique in slightly different manner.
Our Kato inequality represents a mild extension of the existing
theory, the main difference being that we are forced to consider a
$U(n)$, not $SO(n)$ decomposition of tensors. As usual, the improved
Kato inequality yields an improved elliptic inequality, which (via
Uhlenbeck's method) allows one to improve the behavior of $|\Riem|$
at singularities or at infinity.

{\bf Remark.} In an interesting recent work \cite{CS}, a
corresponding precompactness result for K\"ahler-Ricci solitons was
derived with the additional assumption of pointwise Ricci curvature
bounds. These bounds on Ricci curvature in \cite{CS} can be removed
as in our case. The details will be found in a forthcoming paper
\cite{Web2}.

{\bf Organization.} In section \ref{AnalyticLemmas} we consider the
steps necessary for attaining moduli space compactness under our
assumptions, and establish the analytic lemmas needed to overcome
these difficulties. In section \ref{CurvatureLemmas} we state the
necessary estimates and outline the Moser iteration argument needed
for weak compactness. In section \ref{WeakCompactness} we give the
weak compactness argument; we also give the proof of the volume
growth upper bound, and state a gap theorem for ALE extremal
metrics. We also present our adaptation of the argument for
attaining uniform volume growth bounds. Some details will be omitted
from various arguments, as they are nearly identical to those found
elsewhere.

{\bf Acknowledgements}. The first named author has been studying the
moduli space of extremal K\"ahler metrics of complex surfaces since
1999, although this research program is on and off because the issue
of geometric collapsing quickly asserts itself. Special thanks go to
S. Donaldson for many illuminating discussions with the first named
author on this topic. Our work here is also partly motivated by the
work of Tian-Viaclovsky \cite{TV1} \cite{TV2}; and we are grateful
for  insightful discussions with both of them.  Part of this work
was done while the authors were visiting Peking University and the
University of Science and Technology of China in the summer of 2006,
and we wish to thank both Universities for their generous
hospitality.

\section{A quick introduction to K\"ahler geometry}
\subsection{Setup of notations}

Let $M$ be an $n$-dimensional compact K\"ahler manifold. A K\"ahler
metric can be given by its K\"ahler form $\omega$ on $M$. In local
coordinates $z^1, \cdots, z^n $, this $\omega$ has the form
\[
\omega \;=\; \sqrt{-1} \displaystyle \sum_{i,j=1}^n\;g_{i \bar{j}}
d\,z^i\wedge d\,\bar{z}^{j}  \;>\; 0,
\]
where $\{g_{i\bar{j}}\}$ is a positive definite Hermitian matrix
function. The K\"ahler condition requires that $\omega$ is a closed
positive (1,1)-form, or in other words, that
\[
{{\partial g_{i \bar{k}}} \over {\partial z^{j}}} \;=\; {{\partial
g_{j \bar{k}}} \over {\partial z^{i}}}\qquad {\rm and}\qquad
{{\partial g_{k \bar{i}}} \over {\partial\bar{z}^{j}}} \;=\;
{{\partial g_{k \bar{j}}} \over
{\partial\bar{z}^{i}}}\qquad\forall\;i,j,k=1,2,\cdots, n.
\]
The Hermitian metric corresponding to $\omega$ is given by
\[
\sqrt{-1} \;\displaystyle \sum_1^n \; {g}_{\alpha \bar{\beta}} \;
d\,z^{\alpha}\;\otimes d\, \bar{z}^{\beta}.
\]
For simplicity we will often denote by $\omega$ the corresponding
K\"ahler metric. The K\"ahler class of $(M,\omega)$ is the
cohomology class $[\omega]$ in $H^2(M,{\cal R})$. The curvature
tensor is
\[
{\Riem}_{i\bar{j}k\bar{l}} = - {{\partial^2g_{i\bar{j}}} \over
{\partial z^{k} \partial\bar{z}^{l}}} + \displaystyle \sum_
{p,q=1}^n g^{p\bar{q}} {{\partial g_{i \bar{q}}} \over {\partial
z^{k}}} {{\partial g_{p\bar{j}}} \over {\partial \bar{z}^{l}}},
\qquad\forall\;i,j,k,l=1,2,\cdots n.
\]
The Ricci curvature of $\omega$ is locally given by
\[
{\Ric}_{i\bar{j}}\;=\;-{{\partial}^2\log\det(g_{k\bar{l}})\over
{\partial{z}^i\partial\bar{z}^j}},
\]
so its Ricci curvature form is
\[
\Ric\;=\;\sqrt{-1}\displaystyle\sum_{i,j=1}^n\;{\Ric}_{i\bar{j}}d\,z^i
\wedge d\,\bar{z}^{j} \;=\;-\sqrt{-1}
\partial\overline{\partial}\log\,\det(g_{k\bar{l}}).
\]
It is a real, closed (1,1)-form.

\subsection{Historic background and motivation}
In 1982, E. Calabi \cite{Cal1}  proposed to study the critical
metrics of the so called ``Calabi energy'' in  each K\"ahler class:
\[
Ca(\omega) = \displaystyle \int_M\;(R- \underline{R})^2\,\omega^n.
\]
The critical metrics for this functional (the so-called {\sl
extremal} K\"ahler metrics) satisfy the following equation
\[
R_{,\bar\alpha\bar\beta} = 0, \qquad \forall\; \alpha, \; \beta = 1,
2, \cdots, n.
\]
In other words, the extremal K\"ahler metrics are just those where
the complex gradient field of the scalar curvature functions is a
holomorphic vector field. This class includes the K\"ahler-Einstein
metrics, and more generally the constant scalar curvature (cscK)
metrics. The famous conjecture of Calabi states that if the first
Chern class ($C_1$) has a definite sign, then there is a K\"ahler
Einstein metric in the canonical K\"ahler class. The celebrated work
of T. Aubin \cite{Aub} ($C_1<0$), S. T. Yau \cite{Yau} ($C_1 < 0$
and $C_1 =0$) and G. Tian \cite{Tia1} ($C_1 > 0$ for complex
surfaces) settles the Calabi conjecture in these respective cases.
The remaining case ($C_1>0$ and dimension $>2$) is much more
complicated (\cite{Tia2}). In the 1980s, when he introduced the
notion of extremal K\"ahler metrics, E. Calabi initially expected
that there would exists an extremal K\"ahler metric in each K\"ahler
class. This conjecture of Calabi is known to be false as stated
since there are certain algebraic obstructions to the existence of
extremal K\"ahler metrics (\cite{Lev}). We know our list of
obstructions is incomplete however, as Tian \cite{Tia2} constructed
a example where the known obstructions vanish but there is no cscK
metric.

There is relatively little progress on the general existence problem
using PDE methods, although there is very active research in
utilizing the special symmetric structure of underlying K\"ahler
manifold as well as in deploying subtle implicit function methods
(cf. \cite{Cal1} \cite{LS1} \cite{LS2} \cite{ACGT} \cite{AP}
\cite{APS} \cite{Fine} and references therein) to construct (or
prove the existence of) extremal K\"ahler metrics. The present work
is a movement in this direction using geometric methods.


\subsection{Derivation of some useful formulas}

First we show how to derive the elliptic system (\ref{RiemEq}),
(\ref{RicEq}), and (\ref{XEq}).  We note that the first equation
holds for any Riemannian manifold, though the derivation in the
K\"ahler case is simpler.  We compute in unitary frames
\begin{eqnarray*}
{\Riem}_{i\bar{j}k\bar{l},m\bar{m}}
&=&{\Riem}_{i\bar{j}m\bar{l},k\bar{m}}\\
&=&{\Riem}_{i\bar{j}m\bar{l},\bar{m}k}
\,+\,{\Riem}_{k\bar{m}i\bar{s}}{\Riem}_{s\bar{j}m\bar{l}}
\,-\,{\Riem}_{k\bar{m}s\bar{j}}{\Riem}_{i\bar{s}m\bar{l}}\\
&&\qquad\qquad+\,{\Riem}_{k\bar{m}m\bar{s}}{\Riem}_{s\bar{j}s\bar{l}}
\,-\,{\Riem}_{k\bar{m}s\bar{l}}{\Riem}_{i\bar{j}m\bar{s}}\\
&=&{\Ric}_{i\bar{j},\bar{l}k}
\,+\,{\Riem}_{k\bar{m}i\bar{s}}{\Riem}_{s\bar{j}m\bar{l}}
\,-\,{\Riem}_{k\bar{m}s\bar{j}}{\Riem}_{i\bar{s}m\bar{l}}\\
&&\qquad\quad+\,{\Ric}_{k\bar{s}}{\Riem}_{s\bar{j}s\bar{l}}
\,-\,{\Riem}_{k\bar{m}s\bar{l}}{\Riem}_{i\bar{j}m\bar{s}}.
\end{eqnarray*}
When the exact form of the expression is not important we will
denote a linear combination of traces of tensor products of $S$ and
$T$ simply by $S*T$. Using this of abbreviation, we write
\begin{eqnarray*}
\triangle\Riem &=& \Riem*\Riem \,+\,\nabla\overline\nabla\Ric.
\end{eqnarray*}
Next we work with the Ricci tensor, and note that a simplification
of $\triangle\Ric$ is possible in the K\"ahler case because we are
allowed additional permutations of indices.
\begin{eqnarray*}
{\Ric}_{i\bar{j},m\bar{m}} &=&{\Ric}_{m\bar{j},i\bar{m}}\\
&=&{\Ric}_{m\bar{j},\bar{m}i}
\,+\,{\Ric}_{i\bar{s}}{\Ric}_{s\bar{j}}
\,-\,{\Riem}_{i\bar{m}s\bar{j}}{\Ric}_{m\bar{s}}\\
&=&R_{,\bar{j}i} \,+\,{\Ric}_{i\bar{s}}{\Ric}_{s\bar{j}}
\,-\,{\Riem}_{i\bar{m}s\bar{j}}{\Ric}_{m\bar{s}}.
\end{eqnarray*}
The computation for $\Ric_{i\bar{j},\bar{m}m}$ is similar. Using the
notation $X=\overline\nabla{R}$, we get
\begin{eqnarray*}
\triangle\Ric &=&\Riem*\Ric \,+\,\nabla{X}.
\end{eqnarray*}
In the extremal case we can actually get an elliptic equation for
$X$. Recalling that $\overline\nabla{X}=0$ for extremals, a
commutator formula gives
\begin{eqnarray*}
X_{,m\bar{m}} &=&R_{,\bar{i}m\bar{m}}
\;=\;-{\Ric}_{s\bar{i}}R_{,\bar{s}}\\
\triangle{X} &=&\Ric*X.
\end{eqnarray*}
Essentially the same computation gives that
\begin{eqnarray*}
\nabla^2X&=&\Riem*X.
\end{eqnarray*}

\subsection{{\it A priori} bounds on the extremal vector field} \label
{SubsectionBoundsOnExtremalField} In this section we establish
preliminary local estimates for $|X|$ and $|\nabla X|$. It is well
known that, given a K\"ahler manifold and a K\"ahler class, then the
$L^\infty$ norm of its scalar curvature function is {\it a priori}
bounded. Moreover, the extremal vector field $X$ is determined up to
conjugation.  However, one does not expect that the length of $|X|$
with respect to varying extremal metrics has any kind of bound.  We
are pleasantly surprised that, without use of the Sobolev
inequality, we can directly bound $L^2(X)$. By Fatou's lemma, this
result will hold on any manifold-with-singularities that arises as
the limit of such manifolds. Extremal K\"ahler metrics have
automatic upper and lower bounds on scalar curvature which depends
on the complex structure and K\"ahler class. Using this fact, we can
prove
\begin{proposition}\label{L2X}
Assume $M$ is a compact manifold and that $X = \overline{\nabla}R
\triangleq R_{,\bar{i}}\,d\bar{z}^i$ is a holomorphic covector
field. Then
$$
\int_M |X|^2  \;\le\;  C\sup|R| \int_M |\Ric|^2.
$$
and
$$
\int_M |\nabla X|^2  \;\le\;  C\sup|R|^2 \int_M |\Ric|^2.
$$
for a constant $C=C(n)$.
\end{proposition}
\underline{\sl Pf}
\newline
We deal with $L^2(|\nabla X|)$ first. We use formula (\ref{XEq}) in
a more specific form,
$$
R_{,\bar{i}j\bar{j}} = R_{,\bar{i}\bar{j}j} +
{\Riem}_{\bar{j}j\bar{i}k} R_{,\bar{k}} =
-{\Ric}_{k\bar{i}}R_{,\bar{k}},
$$
and integration by parts.  Note that $\nabla X =
R_{,\bar{i}j}+R_{,\bar{i}\bar{j}} = R_{,\bar{i}{j}}$.
\begin{eqnarray*}
\int |\nabla X|^2 &=& \int R_{,i\bar{j}}R_{,\bar{i}j} \;=\; - \int
R_{,i\bar{j}j}R_{,\bar{i}}
\;=\; - \int {\Ric}_{i\bar{k}} R_{,k} R_{,\bar{i}} \\
&=& \int {\Ric}_{i\bar{k},\bar{i}} R_{,k} R  \,+\,  \int
{\Ric}_{i\bar {k}} R_{,k\bar{i}} R \;=\; \int R_{,\bar{k}} R_{,k} R
\,+\, \int {\Ric}_{i\bar{k}} R_{,k \bar{i}} R
\\
&=& \int |X|^2 R \,+\, \int \left<\Ric,\nabla X\right> R \;\le\;
\int |X|^2 R \,+\, \frac12 \int |\Ric|^2 R^2 \,+\, \frac12  \int
|\nabla X|^2
\\
\int |\nabla X|^2 &\le& 2\int|X|^2 R \,+\, \int |\Ric|^2 R^2
\end{eqnarray*}
We use
$$
\frac12\triangle R^2 \,=\, |\nabla R|^2 + R\triangle R \,\ge\, |
\nabla R|^2 - \frac{R}{\sqrt{n}}|\Hess R|^2 \,=\, |X|^2 - \frac{R}
{\sqrt{n}}|\nabla X|.
$$
Then
\begin{eqnarray*}
\int|X|^2 &\le& \frac{1}{\sqrt{n}}\int R|\nabla X| \;\le\;
\left(\frac1n\int R^2 \right)^\frac12\left(\int|\nabla
X|^2\right)^\frac12
\\
&\le& \left( \frac1n\int R^2 \right)^\frac12  \left( 2 \int |X|^2 R
\, +\, \int |\Ric|^2 R^2 \right)^\frac12
\\
&\le& \frac{\sup|R|}{n} \int R^2  \,+\, \frac{1}{2\sup|R|} \int
|X|^2  R \,+\, \frac{1}{4\sup|R|} \int |\Ric|^2 R^2
\\
\int |X|^2 &\le& \frac{2 \sup|R|}{n} \int R^2 \,+\,
\frac{\sup|R|}{2}  \int |\Ric|^2
\end{eqnarray*}
Thus
$$
\int |X|^2  \;\le\; C\sup|R|\int |\Ric|^2.
$$
\hfill$\Box$

\subsection{Uniform Sobolev constant bound} \label{SubsectionUniformSobolevBound}
The large scale aim of this
research program is to contribute to the understanding the
Yau-Tian-Donaldson conjecture and the Calabi conjecture. The most
immediate natural application is the special case of complex
surfaces with K\"ahler class in the so-called generalized Tian cone.
Let us first define
\begin{defi}
The K\"ahler class $\omega$ of a complex surface is in Tian's cone
if
\[
c_1(M)^2\,-\,{2\over3}{{(c_1(M)\cdot[\omega])^2}\over{[\omega]^2}}\;>
\;0.
\]
\end{defi}
\indent\indent A striking observation (\cite{Tia3}, \cite{TV2}) of
Tian's is that a positive cscK metric in the Tian cone automatically
has a uniform Sobelev constant bound. One can modify this to include
the case of extremal K\"ahler metrics: We say a surface's K\"ahler
class lies in the generalized Tian cone if
\begin{eqnarray}
c_1(M)^2\,-\,{2\over
3}\left({{(c_1(M)\cdot[\omega])^2}\over[\omega]^2}
\,+\,\frac{1}{64\pi^2}\|{\cal{F}}\|^2\right) \;>\; 0 \label{TiansKC}
\end{eqnarray}
Here $\|\cal F\|$ is the norm of the Calabi-Futaki invariant
\cite{Fut} in a Mabuchi-Futaki invariant metric \cite{FM}; see
\cite{Chn2} for the definition of this norm. More importantly,
extremal metrics in this modified Tian cone sometimes enjoy similar
properties. In other words, some extremal K\"ahler metrics in a
bounded region of the modified Tian cone have bounds ({\it i})-({\it
iv}) {\it a priori}.

To make sense of this assertion, use
\begin{eqnarray*}
&&\frac{(C_1\cdot[\omega])^2}{[\omega]\cdot[\omega]}\;=\;\frac{1}{32
\pi^2}\frac{1}{\Vol}\left(\int{R}\right)^2
\end{eqnarray*}
and
\begin{eqnarray*}
&&C_1^2\;=\;\frac{1}{96\pi^2}\int\left(R^2-12|E|^2\right)
\,+\,\frac{1}{48\pi^2}\int{R}^2,
\end{eqnarray*}
where $E$ indicates the trace-free Ricci tensor. If the
representative metric happens to be extremal, it turns out that
\begin{eqnarray*}
\|\mathcal{F}\|^2&=&2\left(\int{R}^2-\frac{1}{\Vol}\left(\int{R}
\right)^2\right).
\end{eqnarray*}
A glance at the Chern-Gauss-Bonnet formula for $\chi$ indicates that
$\int\left(R^2-12|E|^2\right)$ is a conformal invariant on
4-manifolds, so when (\ref{TiansKC}) holds, we get a bound on the
square of the Yamabe minimizer. It is well-known that the Sobolev
constant is bounded in the conformal class of a positive Yamabe
minimizer (ref!!), where the constant $A$ in
(\ref{SobolevStatement}) is controlled by the Yamabe constant and
$L^\infty(R)$. So assuming (\ref{TiansKC}) and a positive Yamabe
constant there is abound on the Sobolev constant. Such a bound
holds, for example, on del Pezzo surfaces.


\subsection{Future work}
Due to LeBrun-Simanca \cite{LS1}, it is known that the set of
K\"ahler classes (and bounded complex structures) which admit
extremal K\"ahler metrics is open in the K\"ahler cone. This
suggests that it is possible to pursue the existence of the extremal
K\"ahler metrics using the method of continuity.  In a subsequent
work, we want to study
\begin{prob}\label{Problem1}
Let $\{[\omega_n]\}$ be a sequence of K\"ahler classes which
converges to a limiting K\"ahler class $[\omega_\infty]$. Suppose
that the limiting K\"ahler class is K stable, and suppose that
$\{g_n\}$ is a sequence of extremal K\"ahler metrics in
$\{[\omega_n]\}$ respectively. If the $g_i$ all satisfy conditions
(i)-(iv), do we have a smooth limit as $i\rightarrow\infty$? In
other words, will orbifold singularities fail to develop?
\end{prob}

A special case of this problem, perhaps more natural, is
\begin{prob}\label{Problem2}
In complex dimension 2, can we solve problem \ref{Problem1}? What
about in the interior of the generalized Tian's K\"ahler
cone\footnote{See Section 2.5 for definition.}? What happens at the
border of this modified K\"ahler cone?  What if we don't assume the
limiting class is stable?
\end{prob}

\begin{prob}\label{Problem3}
If we remove the assumption of uniform bound on Sobolev constant,
does some version of Theorem \ref{FinalCompactnessTheorem0} still
hold? What if we restrict to surfaces only?
\end{prob}

In a series of remarkable works \cite{Don1} \cite{Don2} \cite{Don3},
S. K. Donaldson initiated the study of the existence of extremal
K\"ahler metrics on toric surfaces; see also \cite{Zhu} for further
work in toric varieties. This program might be viewed as parallel to
the one described out here. Addressing the problem in full
generality would mean tackling one of two essential difficulties: on
the one hand the lack of 2-dimensional symmetry in general, and on
the other the lack of Sobolev constant control in general. The work
of Cheeger-Tian [CT] on 4 dimensional Einstein manifolds may shed
some light on this problem.

\begin{prob}\label{Problem4}
What can we say about Theorem \ref{FinalCompactnessTheorem0} if we
assume bounds on the $L^2$ (instead of $L^{n\over 2}$) norm of
Riemannian curvature?
\end{prob}
In extremal K\"ahler geometry this is an especially natural
question, as the $L^2$ norm has {\it a priori} bounds, from which we
don't know how to obtain $L^{n\over 2}$ bounds. There are many
important works in this direction by M. Anderson, J. Cheeger, T.
Colding, G. Tian, and others.  Readers are encouraged to browse
\cite{CCT} or \cite{CT} for further details and references.

\section{Analytic Lemmas} \label{AnalyticLemmas}

The results of this section hold for complete manifolds with certain
kinds of singular points, what Anderson calls ``curvature
singularities.''  Specifically,
\begin{defi} Assume $M$ is a length space with a set
$S=\bigcup_{j=1}^{N}\{p_j\}$ such that $M - S$ is a smooth
Riemannian manifold.  If $S$ is the smallest such set, we call it
the singular set of $M$.  If for each $p_j$ there is an
$\epsilon_j>0$ and numbers $0 < \underline{v}_j \le \overline{v}_j$
with the property that $\underline{v}_j r^{n} \le \Vol B(p_j,r) \le
\overline{v}_j r^n$ for $0\le r < \epsilon_j$, then we call $M$ a
manifold-with-singularities, and call the $p_j$ curvature
singularities.
\end{defi}
Our goal in this section is to establish the tools we shall need
later to establish the pointwise bounds for the Riemannian curvature
tensor on manifolds-with-singularities.  This provides the first
step in both the weak compactness and the removable singularity
theorems.

Moser iteration with the elliptic inequality $\triangle{u}\ge-fu-g$
(roughly the form of (\ref{RiemEq}), (\ref{RicEq}), and (\ref{XEq}))
requires the a priori conditions that $u\in L^2$ and $f,g\in L^p$
for some $p>n/2$. We will have only $p=n/2$ a priori. Essentially by
exploiting the {\sl non}linear structure of the system
(\ref{RiemEq}), (\ref{RicEq}), (\ref{XEq}), with methods pioneered
in \cite{BKN}, \cite{Tia1}, \cite{And1}, we can bootstrap $f$ and
$g$ into the needed $L^p$ spaces. The presence of singularities
complicates this, the main difficulty being that integration by
parts leaves an uncontrollable residue at singularities. The first
task is partially recovering integration by parts, which is possible
for functions that are differentiable away from the singular set and
in $L_{loc}^{n/(n-1)}$ at the singularities.

{\bf Remark.} The fact that the Sobolev inequality continues to hold
for $W^{1,2}$-functions across the singular points, assuming the
local upper bound on volume growth, is by now a well known.
\begin{lemma}[Sobolev inequality for $W^{1,2}$ functions]\label{Sobolev}
Assume the Sobolev inequality $\left(\int_U
v^{2\gamma}\right)^\frac1\gamma \,\le\, C_S \int_U |\nabla v|^2 $
holds for all domains $U$ with closure $\overline{U}$ compact and
disjoint from the singular set, with $\Vol U \le \frac12 \Vol M$ if
$\Vol M$ is finite, and with $v\in C_c^1(U)$. Then the Sobolev
inequality holds for functions $v\in W_0^{1,2}(U)$ even if
$\overline{U}$ contains singular points.
\end{lemma}
\underline{\sl Pf} See, for instance, the proof of Theorem 5.1 in
\cite{BKN}\qed

\begin{lemma}[Integration by parts] \label{IntegrationByParts}
Assume $X$ is any vector field with compact support which is smooth
outside the singular set. If $|X|\in L^{\frac{n}{n-1}}$ (or just
$|X|=O(r^{-(n-1)})$ near singularities) and either
$\int(\div(X))_{-}$ or $\int(\div(X))_{+}$ is finite, we retain the
divergence formula: $\int_M d\left(i_X dVol\right) = 0$.
\end{lemma}
\underline{\sl Pf}

Without loss of generality, we assume $(\div(X))_{-}$ is integrable,
and we assume there is only one singularity, at $o$.  For small
values of $r$, let $\phi_r \ge 0$ be a test function with $\phi_r
\equiv 1$ outside $B(o,2r)$, $\phi_r \equiv 0$ inside $B(o,r)$, and
$|\nabla \phi_r| \le 2/r$. Possibly $\int_M\div(X)=+\infty$, but in
any case the Dominated Convergence Theorem and Fatou's lemma give
\begin{eqnarray*}
\int_M\div(X)&=&\int_M(\div(X))_{+} - \int_M(\div(X))_{-} \\
&\le& \lim_{r\rightarrow0^{+}}\int_M\phi_r(\div(X))_{+}
\,-\,\lim_{r\rightarrow0^{+}}\int_{M}\phi_r(\div(X))_{-}\\
&=& \lim_{r\rightarrow 0^{+}} \int_M \phi_r\,\div(X).
\end{eqnarray*}
But
\begin{eqnarray*}
|\int_M \phi_r\,\div(X)| &=& |\int_M \left<X,\nabla \phi_r\right>| \\
&\le& \left(\int_M |\nabla\phi_r|^n\right)^\frac1n
\left(\int_{\supp(\nabla\phi_r)}
|X|^{\frac{n}{n-1}}\right)^{\frac{n-1}{n}}
\end{eqnarray*}
Since $\left(\int_M |\nabla\phi_r|^n\right)^\frac1n \;\le\; \frac2r
\left(\Vol B(o,r)\right)^\frac1n$ is finite and
$\int_{\supp(\nabla\phi_r)} |X|^{\frac{n}{n-1}}$ is bounded as
$r\rightarrow 0$, we get that $\lim_{r\rightarrow0} \int_M
\phi_r\,\div(X)$ is bounded. Therefore
\begin{eqnarray*}
\int_M\div(X) &<& \infty,
\end{eqnarray*}
which proves that indeed $\int(\div(X))_{+}<\infty$. The DCT now
gives
\begin{eqnarray*}
\int_M\div(X)&=&\lim_{r\rightarrow 0^{+}} \int_M \phi_r\,\div(X)
\;=\; 0.
\end{eqnarray*}
\qed



We eventually wish to prove that the curvature singularities are
``removable,'' in the sense that the Riemann curvature tensor has
uniform pointwise bounds in the neighborhood of any singular point.
Thus the singularity may be topologically nontrivial, but its metric
structure will be controlled, and in all cases one can prove such a
singularity will be a Riemannian orbifold point of regularity at
least $C^0$.

The first step in the removable singularity theorem is establishing
that $|\Riem|\in{L}_{loc}^p$ for some $p>n/2$. A result of
\cite{BKN} is that if $|\Riem|=O(r^{-2+\alpha})$ for any $\alpha>0$,
one can construct coordinates with $C^{1,1}$ bounds on metric
components. In fact, if one can obtain just ${C}^{1,\alpha}$
coordinates, one has access to harmonic coordinates (\cite{DK}) and
a bootstrapping argument can commence, which we give in section
\ref{WeakCompactness}.

For dimensions 6 and up, we obtain $|\Riem|\in{L}_{loc}^p$ using
analytic methods first developed in \cite{Sib}. Sibner's original
purpose was to prove removable singularity theorems for Yang-Mills
instantons, a problem closely related to ours. This method was used
again by Cao-Sesum in \cite{CS} to remove singularities on
K\"ahler-Ricci solitons. Sibner's theorem is really only useful in
dimension 5 and higher; in the Yang-Mills case other methods were
used in dimensions 2, 3, and 4. We use other methods in dimension 4
as well; see section \ref{SubsectionUhlenbeck}. The proof below does
have some limited applicability in dimensions $3$ and $4$.

\begin{lemma}[$u^k \in L^2$ implies $\nabla u^k \in L^2$]\label{L2Lemma}
(\cite{Sib}) Assume 2-sided volume growth bounds, Sobolev constant
bounds, and $\triangle u \ge -fu$ where $f\in{L}^{n/2}(B-\{o\})$ and
$u\ge0$. If $k>\frac12\frac{n}{n-2}$, then $u^k \in L^2(B-\{o\})$
implies $\nabla u^k \in L^2(B-\{o\})$.
\end{lemma}
\underline{\sl Pf}\\
\indent The idea is to dampen the growth near the singularity while
retaining weakly an elliptic inequality. Assume $\frac12<q_0\le{q}$,
to be chosen later. We set up a test function as follows. Let
$$
F(t) =
\begin{cases}
t^q & {\rm if} \; 0 \le t \le l \\
\frac{1}{q_0}\left(ql^{q-q_0} t^{q_0} + (q_0-q)l^q\right) & {\rm if}
\; l \le t,
\\
\end{cases}
$$
and set $G(t) = F(t) F'(t)$. We shall need the following easily
verified facts:
\begin{eqnarray}
&&F(t) \le \frac{q}{q_0}l^{q-q_0}t^{q_0} \label{PropertyI} \\
&&q_0F(t) \le tF'(t) \label{PropertyII} \\
&&tF'(t) \le qF(t) \label{PropertyIII} \\
&&G'(t) \ge \frac{2q_0-1}{q_0}(F'(t))^2 \label{PropertyIV}.
\end{eqnarray}
For a test function $\zeta$,
\begin{eqnarray}
\int \left<\nabla\zeta,\nabla u\right>
&\le&\int\zeta{u}f\label{vEqn}
\end{eqnarray}
Choose $\zeta = \eta^2 G(u)$ for our test function. We have to
assume $\eta \equiv 0$ across any singularities in order to make
integration by parts work. The trick will be to make $u$ disappear.
We have
\begin{eqnarray*}
\left<\nabla u,\nabla\zeta\right>
\;\ge\;2\eta{F}(u)F'(u)\left<\nabla{u},\nabla\eta\right>
\,+\,\frac{2q_0-1}{q_0}\eta^2(F'(u))^2|\nabla{u}|^2,
\end{eqnarray*}
so combining with (\ref{vEqn}) and simplifying gives our main
inequality:
\begin{eqnarray}
&&\int2\eta{F}(u)\left<\nabla\eta,\nabla{F}(u)\right>
\,+\,\frac{2q_0-1}{q_0}\int\eta^2|\nabla F(u)|^2\\
&&\quad\le\;q\int\eta^2(F(u))^2f\label{MainIneq}
\end{eqnarray}
We deal with the terms one by one. The first term on the left of
(\ref{MainIneq}) is easily dispatched with Schwartz:
\begin{eqnarray}
\int2\eta{F}(u)\left<\nabla\eta,\nabla{F}(u)\right>
&\ge&-\mu\int\eta^2|\nabla{F}(u)|^2
\,-\,\frac1\mu\int|\nabla\eta|^2|F(u)|.
\end{eqnarray}
For the term on the right of (\ref{MainIneq}), H\"older and Sobolev
give
\begin{eqnarray*}
&&\int\eta^2(F(u))^2f
\;\le\;\left(\int{f}^{\frac{n}{2}}\right)^\frac2n\left(\int\eta^{\frac
{2n}{n-2}}(F(u))^{\frac{2n}{n-2}}\right)^{\frac{n-2}{n}}\\
&&\quad\le\;2C_S^2\left(\int{f}^{\frac{n}{2}}\right)^\frac2n\left
(\int|\nabla\eta|^2F(u)^2\right)
\,+\,2C_S^2\left(\int{f}^{\frac{n}{2}}\right)^\frac2n\left(\int\eta^2|
\nabla{F}(u)|^2\right)
\end{eqnarray*}
Putting everything back into (\ref{MainIneq}) and simplifying now
gives
\begin{eqnarray}
&&\left(\frac{2q_0-1}{q_0}-\mu-2qC_S^2\left(\int{f}^{\frac{n}{2}}
\right)^\frac2n\right)\int\eta^2|\nabla
F(u)|^2\nonumber\\
&&\quad\le\;\left(\frac1\mu+2qC_S^2\left(\int{f}^{\frac{n}{2}}\right)^
\frac2n\right)\int|\nabla\eta|^2|F(u)|^2
\end{eqnarray}
After $q$ and $q_0$ are chosen, we choose $\mu$ to be small and
choose the cutoff $\eta$ so that
$\left(\int_{\supp(\eta)}{f}^{\frac{n}{2}}\right)^\frac2n$ is also
small. Under these conditions we get that
\begin{eqnarray}
\int\eta^2|\nabla{F}(u)|^2\;\le\;C\int|\nabla\eta|^2|F(u)|^2,
\label{MainIneqSimplified}
\end{eqnarray}
where $C = C(q_0)$, provided that
$\|f\|_{L^{\frac{n}{2}}(\supp(\eta))}\le{C}(C_S,q,q_0)$.

We want to ferret out the contribution at the singularity, so
replace $\eta$ with $\eta \eta_\epsilon$, where now $\eta\equiv1$
across the singularity, and $\eta_\epsilon \ge 0$ is another cutoff
function with $\eta_\epsilon \equiv 1$ outside $B(o,2\epsilon)$,
$\eta_\epsilon \equiv 0$ inside $B(o,\epsilon)$, and
$|\nabla\eta_\epsilon|\le 2/\epsilon$. Using $F(v) \le
\frac{q}{q_0}l^{q-q_0}u^{q_0}$ and applying H\"older again:
\begin{eqnarray*}
\int(\eta\eta_\epsilon)^2|\nabla F(u)|^2
&\le&C\left(\frac{q}{q_0}l^{q-q_0}\right)^2\left(\int|\nabla\eta_
\epsilon|^n\right)^\frac2n\left(\int_{\supp(\nabla\eta_\epsilon)}u^
{\frac{2nq_0}{n-2}}\right)^{\frac{n-2}{n}}
\\
&&+\,C\int|\nabla\eta|^2|F(u)|^2.
\end{eqnarray*}
Now choose $q_0>\frac12$ so $q_0=k(n-2)/n$ (here we use the
hypothesis that $k>\frac12\frac{n}{n-2}$). Then $\frac{2nq_0}{n-2} =
2k$ and so $u^{\frac{2nq_0}{n-2}}$ is locally integrable. As
$\epsilon\rightarrow0$ we get
\begin{eqnarray*}
\left(\int|\nabla\eta_\epsilon|^n\right)^\frac2n
\left(\int_{\supp(\nabla\eta_\epsilon)}u^{\frac{2nq_0}{n-2}}\right)^
{\frac{n-2}{n}} \rightarrow0,
\end{eqnarray*}
So $\int\eta^2|\nabla F(u)|^2\le C\int|\nabla\eta|^2|F(u)|^2$. Now
letting also $l\rightarrow\infty$ we finally get
\begin{eqnarray}
\int\eta^2|\nabla{u^q}|^2
\;\le\;C\int|\nabla\eta|^2|u^q|^2.\label{SibnersNablaTransfer}
\end{eqnarray}
Choosing $q=k$ (so automatically $q>q_0$), we have our result:
\begin{eqnarray*}
\nabla{u}^k \in L^2.
\end{eqnarray*}
\hfill$\Box$

\begin{proposition}[$L^p$-regularity]\label{LpLemma}
Assume $\triangle{u}\ge-fu-g$, $u\ge0$ in $B-\{o\}$, with
$f,g\in{L}^{n/2}(B-\{o\})$, and assume 2-sided volume growth bounds
at the singular point and a finite Sobolev constant. If $u\in
L^q(B-\{o\})$ for some $q>\frac{n}{n-2}$, then $u\in L^p(B-\{o\})$
for all $\infty>p\ge{q}$. Explicitly, with $a>q>\frac{n}{n-2}$,
there exists $\epsilon_0=\epsilon_0(q,a,C_S)$, $C=C(q,a,C_S,n)$ so
that $\int_{B(o,r)}f^{\frac{n}{2}}\le\epsilon_0$ implies
\begin{eqnarray}
\left(\int_{B(o,r/2)}u^a\right)^\frac1a
\;\le\;Cr^{\frac{n}{a}-\frac{n}{q}}\left(\int_{B(o,r)}u^q\right)^\frac1q
\,+\,Cr^{\frac{n}{a}}\left(\int_{B(o,r)}g^{\frac{n}{2}}\right)^\frac2n.
\end{eqnarray}
\end{proposition}
\underline{\sl Pf}\\
We must pay special attention to any use of integration by parts;
otherwise the argument is standard. Assume $p>1$. Replace $u$ by
$u+\mathcal{C}\|g\|_{L^{\frac{n}{2}}(B(o,r))}$ and $f$ by
$f+\frac{1}{\mathcal{C}}\frac{g}{\|g\|_{L^{\frac{n}{2}}(B(o,r))}}$,
where $\mathcal{C}$ is some number to be chosen later; it will be
roughly $\left(a^2/(a-1)\right)^{n/2}$.  Then $\triangle u \ge -fu$.
We get
$$
\left(\int \left(\eta^2
u^p\right)^{\frac{n}{n-2}}\right)^{\frac{n-2} {n}} \;\le\; 2 C_S^2
\int|\nabla \eta|^2 u^p \,+\, 2 p^2 C_S^2 \int \eta^2  u^{p-2}
|\nabla u|^2.
$$
The last term reads $C\int\eta^2|\nabla{u}^{\frac{p}{2}}|^2,$ which
we can estimate using (\ref{SibnersNablaTransfer}).  This estimate
requires that $\left(\int_{B(o,r)}|f|^{\frac{n}{2}}\right)^\frac2n$
be small compared to $p$ and $C_S$, which, incidentally requires
choosing $\mathcal{C}$. We get
\begin{eqnarray}
\left(\int\left(\eta^2u^{p}\right)^{\frac{n}{n-2}}\right)^{\frac{n-2}
{n}} &\le&C\int|\nabla\eta|^2u^p\label{LocalIneq}
\end{eqnarray}
where $C = C(p,C_S)$. Iterating this inequality will give $u\in L^p$
for all $q\le{p}<\infty$.

We carry this out explicitly. With $0<k<1$ and an appropriate choice
of test functions $\phi$, (\ref{LocalIneq}) implies
\begin{eqnarray*}
&&\left(\int_{B(o,kr)}u^{p\gamma}\right)^\frac1\gamma\;\le\;Cr^{-2}
\int_{B(o,r)}u^p,
\end{eqnarray*}
with $C=C(p,k,C_S)$, and iterating, we get
\begin{eqnarray*}
&&\left(\int_{B(o,k^{i+1}r)}u^{p\gamma^{i+1}}\right)^{\frac{1}{\gamma^
{i+1}}}\;\le\;Cr^{\frac{n}{\gamma^{i+1}}-n}\int_{B(o,r)}u^p,
\end{eqnarray*}
with $C=C(p,k,i,C_S)$.  Now choose $i$ so
$p\gamma^i\,\le\,a\,<\,p\gamma^{i+1}$. Then
\begin{eqnarray*}
&&\int_{B(o,k^{i+1}r)}u^a \;\le\;
r^{n-\frac{na}{p}}\left(\int_{B(o,k^{i+1}r)}u^{p\gamma^i}\right)^
{\frac{p\gamma^{i+1}-a}{p\gamma^{i+1}-p\gamma^i}}\left(\int_{B(o,k^{i
+1}r)}u^{p\gamma^{i+1}}\right)^{\frac{a-p\gamma^i-a}{p\gamma^{i+1}-p
\gamma^i}}\\
&&\quad\quad\le\;Cr^{n-\frac{na}{p}}\left(\int_{B(o,kr)}u^p\right)^
{\gamma^i\frac{p\gamma^{i+1}-a}{p\gamma^{i+1}-p\gamma^i}}\left(\int_{B
(o,r)}u^p\right)^{\gamma^{i+1}\frac{a-p\gamma^i}{p\gamma^{i+1}-p
\gamma^i}}\\
&&\quad\quad\le\;Cr^{n-\frac{na}{p}}\left(\int_{B(o,r)}u^p\right)^
{\frac{a}{p}},
\end{eqnarray*}
where $C=C(p,k,a,C_S)$. Now lastly put
$u+\mathcal{C}\|g\|_{L^{\frac{n}{2}}(B(o,r))}$ back in for $u$:
\begin{eqnarray*}
&&\left(\int_{B(o,kr)}u^a\right)^\frac1a
\;\le\;\left(\int_{B(o,r)}(u+\mathcal{C}\|g\|)^a\right)^\frac1a\\
&&\quad\quad\le\;Cr^{\frac{n}{a}-\frac{n}{p}}\left(\int_{B(o,r)}(u+
\mathcal{C}\|g\|)^p\right)^\frac1p\\
&&\quad\quad\le\;Cr^{\frac{n}{a}-\frac{n}{p}}\left(\int_{B(o,r)}u^p
\right)^{\frac1p}
\,+\,Cr^{\frac{n}{a}}\left(\int_{B(o,r)}g^{\frac{n}{2}}\right)^
\frac2n.\,
\end{eqnarray*}
for $C=C(a,n,k,C_S)$. \qed

\section{Regularity of sectional curvature}

\subsection{Statement of the curvature estimates} \label {CurvatureLemmas}

In this section we state our main curvature integral estimates, and
actually establish them in the low order case. The method of proof
is standard, but establishing the estimates in the possible presence
of singularities is more complicated. At smooth points, Propositions
\ref{GoodXEstimate}, \ref{GoodRicEstimate}, \ref{GoodRiemEstimate},
and \ref{GoodNablaRicAndXEstimate} give the result for small values
of $q$. The subject of Sections \ref{SubsectionRemCurvSingsDim6} and
\ref{SubsectionUhlenbeck} is to prove the $q=0$ case at singular
points. The rest of the long, unenlightening proof by induction is
consigned to the appendix.
\begin{theorem} \label{LocalCurvatureBounds}
Assume $g$ is an extremal K\"ahler metric on a Riemannian
manifold-with-singularities.  When $a>\frac{n}{2}$, and
$q\in\{0,1,\dots\}$, there exists $\epsilon_0 =
\epsilon_0(C_S,a,q,n)$ and $C=C(C_S,a,q,n)$ so that
$$
\int_{B(o,r)}|\Riem|^{\frac{n}{2}} \;\le\; \epsilon_0
$$
implies
\begin{eqnarray}
&&\left(\int_{B(o,r/2)}|\nabla^q\,X|^a\right)^\frac1a
\;\le\;Cr^{-3-q+\frac{n}{a}}\left(\int_{B(o,r)}|R|^{\frac{n}{2}}
\right)^\frac2n\label{XReg}\\
&&\left(\int_{B(o,r/2)}|\nabla^q\Ric|^a\right)^\frac1a
\;\le\;Cr^{-2-q+\frac{n}{a}}\left(\int_{B(o,r)}|\Ric|^{\frac{n}{2}}
\right)^\frac2n\label{RicReg}\\
&&\left(\int_{B(o,r/2)}|\nabla^q\Riem|^a\right)^\frac1a
\;\le\;Cr^{-2-q+\frac{n}{a}}\left(\int_{B(o,r)}|\Riem|^{\frac{n}{2}}
\right)^\frac2n.\label{RiemReg}
\end{eqnarray}
In the presence of singularities estimate (\ref{RiemReg}) holds if
$q=0$ and $n\ge4$, and estimates (\ref{RicReg}) and (\ref{XReg})
hold if $q=0,1$ and $n\ge6$.  In all other cases the estimates hold
if $B(o,r)$ consists of manifold points.
\end{theorem}

We begin the induction argument for Proposition
\ref{LocalCurvatureBounds} at smooth points, for $|\Riem|$, $|\Ric|$
and $|\nabla\Ric|$, and $|X|$ and $|\nabla{X}|$. It is worth noting
that the arguments here work in real dimension $n\ge3$.

\begin{proposition} \label{GoodXEstimate}
If $p>\frac{n}{2}$ and $B(o,r)$ consists of smooth points, there
exists $\epsilon_0=\epsilon_0(p,C_S,n)$ and $C=C(p,C_S,n)$ so that
$\int_{B(o,r)}|\Ric|^{\frac{n}{2}} \le \epsilon_0$ implies
\begin{eqnarray*}
\left(\int_{B(o,r/2)}|X|^p\right)^\frac1p \;\le\;
Cr^{\frac{n}{p}-3}\left(\int_{B(o,r)}|R|^{\frac{n}{2}}\right)^\frac2n
\end{eqnarray*}
\end{proposition}
\underline{\sl Pf}\\
\indent This is basically a local version of Proposition \ref{L2X}.
We obtain the estimates in a series of steps. First,
\begin{eqnarray}
\int\phi^2|X|^2
&=&-2\int\phi{R}\left<\nabla\phi,X\right>\,-\,\int\phi^2R\triangle{R}
\nonumber\\
\int\phi^2|X|^2 &\le&4\int|\nabla\phi|^2R^2
\,+\,2\left(\int{R}^2\right)^\frac12\left(\int\phi^4|\triangle{R}|^2
\right)^\frac12\label{BoundingScalFirst}
\end{eqnarray}
Then we estimate the last term, using
$|\triangle{R}|^2=R_{,m\bar{m}}R_{,k\bar{k}}$. We get
\begin{eqnarray*}
\int\phi^4|\triangle{R}|^2
&=&-4\int\phi^3\triangle{R}\left<\nabla\phi,X\right>\,-\,\int\phi^4R_
{,m\bar{m}k}R_{\bar{k}}\\
&=&-4\int\phi^3\triangle{R}\left<\nabla\phi,X\right>\,+\,\int\phi^4
{\Ric}_{k\bar{s}}R_{,s}R_{\bar{k}}\\
\int\phi^4|\triangle{R}|^2
&\le&16\int\phi^2|\nabla\phi|^2|X|^2\,+\,2\int\phi^4\Ric(X,X).
\end{eqnarray*}
It is also possible to estimate the $\Ric(X,X)$ term:
\begin{eqnarray*}
&&\int\phi^4\Ric(X,X)\\
&&\quad=\;-4\int\phi^3R\,{\Ric}_{k\bar{s}}\phi_{,s}R_{,\bar{k}}
\,-\,\int\phi^4R\,{\Ric}_{k\bar{s},s}R_{,\bar{k}}
\,-\,\int\phi^4R\,{\Ric}_{k\bar{s}}R_{,\bar{k}s}\\
&&\quad\le\;2\int\phi^2|\nabla\phi|^2|X|^2
\,+\,2\int\phi^4|R|^2|\Ric|^2 \,-\,\int\phi^4{R}|X|^2
\,+\,\frac14\int\phi^4|\nabla{X}|^2
\end{eqnarray*}
Finally we have to estimate the $|\nabla{X}|^2$ term.
\begin{eqnarray*}
\int\phi^4|\nabla{X}|^2
&=&-4\int\phi^3R_{,i\bar{j}}\phi_{,j}R_{,\bar{i}}
\,-\,\int\phi^4R_{,i\bar{j}j}R_{,\bar{i}}\\
&=&-4\int\phi^3R_{,i\bar{j}}\phi_{,j}R_{,\bar{i}}
\,+\,\int\phi^4{\Ric}_{i\bar{s}}R_{,s}R_{,\bar{i}}\\
\int\phi^4|\nabla{X}|^2 &\le&16\int\phi^2|\nabla\phi|^2|X|^2
\,+\,2\int\phi^4\Ric(X,X).
\end{eqnarray*}
Now we successively put these estimates back. First we get
\begin{eqnarray*}
\int\phi^4\Ric(X,X) &\le&12\int\phi^2|\nabla\phi|^2|X|^2
\,+\,4\int\phi^4|R|^2|\Ric|^2 \,-\,2\int\phi^4{R}|X|^2.
\end{eqnarray*}
Note that this also provides
\begin{eqnarray*}
\int\phi^4|\nabla{X}|^2 &\le&40\int\phi^2|\nabla\phi|^2|X|^2
\,+\,8\int\phi^4|R|^2|\Ric|^2 \,-\,4\int\phi^4{R}|X|^2
\end{eqnarray*}
and
\begin{eqnarray}
&&\int\phi^4|\triangle{R}|^2 \nonumber\\
&&\quad\le\;40\int\phi^2|\nabla\phi|^2|X|^2
\,+\,8\int\phi^4|R|^2|\Ric|^2
\,-\,4\int\phi^4{R}|X|^2.\label{TriangleScalInitial}
\end{eqnarray}
Using the Sobolev inequality, we can do something with the final
term:
\begin{eqnarray*}
\int\phi^4{R}|X|^2
&\le&\left(\int{R}^{\frac{n}{2}}\right)^\frac2n\left(\int\phi^{4
\gamma}|X|^{2\gamma}\right)^\frac1\gamma\\
&\le&4C_S\left(\int{R}^{\frac{n}{2}}\right)^\frac2n\int\phi^2|\nabla
\phi|^2|X|^2
\,+\,2C_S\left(\int{R}^{\frac{n}{2}}\right)^\frac2n\int\phi^4|\nabla {X}|^2\\
&\le&4C_S\left(\int{R}^{\frac{n}{2}}\right)^\frac2n\int\phi^2|\nabla
\phi|^2|X|^2
\,+\,16C_S\left(\int{R}^{\frac{n}{2}}\right)^\frac2n\int\phi^4|R|^2| \Ric|^2\\
&&\,-\,8C_S\left(\int{R}^{\frac{n}{2}}\right)^\frac2n\int\phi^4{R}|X| ^2\\
\int\phi^4{R}|X|^2
&\le&C\left(\int{R}^{\frac{n}{2}}\right)^\frac2n\int\phi^2|\nabla\phi|
^2|X|^2
\,+\,C\left(\int{R}^{\frac{n}{2}}\right)^\frac2n\int\phi^4|R|^2|\Ric|^2
\end{eqnarray*}
Remarkably the constant $C$ is bounded independently of the Sobolev
constant. Thus
\begin{eqnarray}
\int\phi^4|\triangle{R}|^2 &\le&C\int\phi^2|\nabla\phi|^2|X|^2
\,+\,C\int\phi^4|R|^2|\Ric|^2. \label{TriangleScalIneq}
\end{eqnarray}
Returning to (\ref{BoundingScalFirst}), we get
\begin{eqnarray*}
\int\phi^2|X|^2 &\le&4\int|\nabla\phi|^2R^2
\,+\,C\left(\int{R}^2\right)^\frac12\left(\int\phi^2|\nabla\phi|^2|X|
^2\right)^\frac12\\
&&+\,C\left(\int{R}^2\right)^\frac12\left(\int\phi^4|R|^2|\Ric|^2
\right)^\frac12.
\end{eqnarray*}
Using $|\nabla\phi|\le\frac2r$ gives us
\begin{eqnarray}
\int\phi^2|X|^2 &\le&Cr^{-2}\int{R}^2
\,+\,C\left(\int{R}^2\right)^\frac12\left(\int\phi^4|R|^2|\Ric|^2
\right)^\frac12.\label{L2XPrettyGoodEstimate}
\end{eqnarray}

We must deal with the final term. In the case $n\ge8$, we can easily
deal with the final term:
\begin{eqnarray*}
\int\phi^4|R|^2|\Ric|^2
&\le&\left(\int|\Ric|^{\frac{n}{2}}\right)^{\frac4n}\left(\int|R|^
{\frac{2n}{n-4}}\right)^{\frac{n-4}{n}}\\
&\le&\left(\Vol\supp\phi\right)^{\frac{n-8}{n}}\left(\int|\Ric|^{\frac
{n}{2}}\right)^{\frac4n}\left(\int|R|^{\frac{n}{2}}\right)^\frac4n
\end{eqnarray*}
The case $n=6$ is more difficult. H\"older gives
\begin{eqnarray*}
\int\phi^4|\Ric|^2|R|^2
&\le&\left(\int\phi^4|\Ric|^3\right)^\frac23\left(\int\phi^4|R|^6
\right)^\frac13,
\end{eqnarray*}
and we use the Sobolev inequality to get
\begin{eqnarray*}
\left(\int\phi^4|R|^6\right)^\frac23
&\le&4C_S\int\phi^2|\nabla\phi|^2|R|^4
\,+\,4C_S\int\phi^4|R|^2|\nabla{R}|^2.
\end{eqnarray*}
Now integration by parts on the last term yields
\begin{eqnarray*}
\left(\int\phi^4|R|^6\right)^\frac23
&\le&16C_S\int\phi^2|\nabla\phi|^2|R|^4
\,-\,4C_S\int\phi^4R^3\triangle{R}.
\end{eqnarray*}
Using H\'older's inequality and (\ref{TriangleScalIneq}) gives
\begin{eqnarray*}
\left(\int\phi^4|R|^6\right)^\frac23
&\le&16C_S\int|\nabla\phi|^2|R|^4
\,+\,4C_S\left(\int\phi^4R^6\right)^\frac12\left(\int\phi^4|\triangle
{R}|^2\right)^\frac12\\
&\le&C\int|\nabla\phi|^2|R|^4
\,+\,C\left(\int\phi^4R^6\right)^\frac12\left(\int\phi^2|\nabla\phi|
^2|X|^2\,+\,\int\phi^4|R|^2|\Ric|^2\right)^\frac12\\
\left(\int\phi^4|R|^6\right)^\frac23 &\le&C\int|\nabla\phi|^4|R|^2
\,+\,C\left(\int\phi^2|\nabla\phi|^2|X|^2\right)^2\,+\,C\left(\int
\phi^4|R|^2|\Ric|^2\right)^2
\end{eqnarray*}
Putting this back in, we get
\begin{eqnarray*}
\int\phi^4|\Ric|^2|R|^2
&\le&C\left(\int\phi^4|\Ric|^3\right)^\frac23\int|\nabla\phi|^4|R|^2
\,+\,C\left(\int\phi^4|\Ric|^3\right)^\frac23\int\phi^2|\nabla\phi|^2|
X|^2.
\end{eqnarray*}

Finally we work with the case $n=4$. We use simply
\begin{eqnarray*}
\int\phi^4|\Ric|^2|R|^2
&\le&\left(\int\phi^4|\Ric|^4\right)^\frac12\left(\int\phi^4R^4\right)
^\frac12,
\end{eqnarray*}
and use the Sobolev inequality to get
\begin{eqnarray*}
\left(\int\phi^4|\Ric|^4\right)^\frac12
&\le&2C_S\int|\nabla\phi|^2|\Ric|^2
\,+\,2C_S\int\phi^2|\nabla\Ric|^2.
\end{eqnarray*}
Using integration by parts on the last term lets us obtain
\begin{eqnarray*}
\left(\int\phi^4|\Ric|^4\right)^\frac12
&\le&C\int|\nabla\phi|^2|\Ric|^2 \,+\,C\int\phi^2|\Riem||\Ric|^2
\,+\,C\int\phi^2|\Ric||\nabla{X}|.
\end{eqnarray*}
Using our expressions for $\int\phi^4|\nabla{X}|^2$ and
$\int{R}|X|^2$ allows us to obtain
\begin{eqnarray*}
\left(\int\phi^4|\Ric|^4\right)^\frac12
&\le&C\int|\nabla\phi|^2|\Ric|^2
\,+\,C\left(\int|\Ric|^2\right)^\frac12\left(\int\phi^2|\nabla\phi|^2|
X|^2\right)^\frac12.
\end{eqnarray*}
The Sobolev inequality applied to $\int\phi^4|R|^4$ gives us
\begin{eqnarray*}
\left(\int\phi^4R^4\right)^\frac12 &\le&C\int|\nabla\phi|^2R^2
\,+\,C\int\phi^2|X|^2
\end{eqnarray*}

Putting the estimates for $\int\phi^4R^2|\Ric|^2$ in the cases
$n=4$, $n=6$, and $n\ge8$ into (\ref{L2XPrettyGoodEstimate}) lets us
conclude, regardless of dimension, that
\begin{eqnarray*}
\int\phi^2|X|^2 &\le&Cr^{-2}\int{R}^2
\,+\,r^{n-6}\left(\int{R}^\frac{n}{2}\right)^\frac4n
\end{eqnarray*}

The conclusion now follows from Proposition \ref{LpLemma}. \qed

\begin{proposition} \label{GoodRicEstimate}
If $p>\frac{n}{2}$ and $B(o,r)$ consists of smooth points, there
exists $\epsilon_0=\epsilon_0(p,C_S,n)$ and $C=C(p,C_S,n)$ so that
$\int_{B(o,r)}|\Riem|^{\frac{n}{2}} \le \epsilon_0$ implies
\begin{eqnarray*}
\left(\int_{B(o,r/2)}|\Ric|^p\right)^\frac1p
&\le&Cr^{\frac{n}{p}-2}\left(\int_{B(o,r/2)}|\Ric|^{\frac{n}{2}}
\right)^\frac2n
\end{eqnarray*}
\end{proposition}
\underline{\sl Pf}\\
\indent We use integration by parts to get
\begin{eqnarray*}
\int\phi^l|\nabla{X}|^k
&=&\int\phi^l|\nabla{X}|^{k-2}\left<\nabla{X},\nabla{X}\right>\\
&=&-l\int\phi^{l-1}|\nabla{X}|^{k-2}\left<\nabla\phi\otimes{X},\nabla
{X}\right>\\
&&-(k-2)\int\phi^l|\nabla{X}|^{k-3}\left<\nabla|\nabla{X}|\otimes{X},
\nabla{X}\right>\\
&&-\phi^l\int|\nabla{X}|^{k-2}\left<X,\triangle{X}\right>
\end{eqnarray*}
Using H\"older's inequality on the first term and using
$\nabla^2{X}=\Riem*X$, we get
\begin{eqnarray*}
\int\phi^l|\nabla{X}|^k
&\le&C(k,l)\int\phi^{l-2}|\nabla\phi|^2|\nabla{X}|^{k-2}|X|^2
\,+\,C(k,l)\phi^l\int|\nabla{X}|^{k-2}|X|^2|\Riem|
\end{eqnarray*}
Assuming $k<n$ we can use H\"older's inequality again to get
\begin{eqnarray}
\int\phi^l|\nabla{X}|^k
&\le&\left(\int\phi^{\frac{ln}{n-k}}|X|^{\frac{kn}{n-k}}\right)^{\frac
{n-k}{n}}.\label{APreliminaryForNablaX}
\end{eqnarray}
This holds in particular when $k=\frac{n}{2}$. Now the inequality
$\triangle|\Ric|\ge-{C}|\Riem||\Ric|-C|\nabla{X}|$ yields the
conclusion, via Proposition \ref{LpLemma}. \qed

\begin{proposition}\label{GoodRiemEstimate}
If $p>\frac{n}{2}$ and $B(o,r)$ consists of smooth points, there
exists $\epsilon_0=\epsilon_0(p,C_S,n)$ and $C=C(p,C_S,n)$ so that
$\int_{B(o,r)}|\Riem|^{\frac{n}{2}} \le \epsilon_0$ implies
\begin{eqnarray*}
\left(\int_{B(o,r/2)}|\Riem|^p\right)^\frac1p
&\le&Cr^{\frac{n}{p}-2}\left(\int_{B(o,r/2)}|\Riem|^{\frac{n}{2}}
\right)^\frac2n
\end{eqnarray*}
\end{proposition}
\underline{\sl Pf}\\
\indent Following the calculation leading up to
(\ref{LPNormaKnowledge}), we get that
\begin{eqnarray*}
&&\left(\int\phi^{k\gamma}|\Riem|^{k\gamma}\right)^\frac1\gamma\\
&&\quad\le{C}\int\phi^{l-2}|\nabla\phi|^2|\Riem|^k
\,+\,C\int\phi^l|\Riem|^{k+1}
\,+\,\int\phi^l|\Riem|^{k-1}|\nabla{X}|
\end{eqnarray*}
holds when $\supp\phi$ consists of smooth points.  The second term
on the right easily combines into the left side when
$\int_{\supp\phi}|\Riem|^{\frac{n}{2}}$ is small, and then using
H\"older's inequality on the rightmost term, we get
\begin{eqnarray*}
\left(\int\phi^{k\gamma}|\Riem|^{k\gamma}\right)^\frac1\gamma
&\le&{C}\int\phi^{l-2}|\nabla\phi|^2|\Riem|^k
\,+\,\left(\int\phi^{\Box}|\nabla{X}|^{\frac{nk}{2k+n-2}}\right)^
{\frac{2k+n-2}{n}}.
\end{eqnarray*}
Noticing that $\frac{nk}{2k+n-2}<n$ and using
(\ref{APreliminaryForNablaX}) gives
\begin{eqnarray*}
\left(\int\phi^{k\gamma}|\Riem|^{k\gamma}\right)^\frac1\gamma
&\le&{C}r^{-2}\int\phi^{l-2}|\Riem|^k.
\end{eqnarray*}
Iterating this inequality yields the conclusion. \qed

\begin{proposition}\label{GoodNablaRicAndXEstimate}
If $p>\frac{n}{2}$ and $B(o,r)$ consists of smooth point, there
exists $\epsilon_0=\epsilon_0(p,C_S,n)$ and $C=C(p,C_S,n)$ so that
$\int_{B(o,r)}|\Riem|^{\frac{n}{2}} \le \epsilon_0$ implies
\begin{eqnarray*}
\left(\int_{B(o,r/2)}|\nabla\Ric|^p\right)^\frac1p
&\le&Cr^{\frac{n}{k}-3}\left(\int_{B(o,r/2)}|\Ric|^{\frac{n}{2}}
\right)^\frac2n\\
\left(\int_{B(o,r/2)}|\nabla{X}|^p\right)^\frac1p
&\le&Cr^{\frac{n}{k}-4}\left(\int_{B(o,r/2)}|R|^{\frac{n}{2}}\right)^
\frac2n\\
\end{eqnarray*}
\end{proposition}
\underline{\sl Pf}\\
\indent Applying the Sobolev inequality, integration by parts, and
the elliptic inequality for $|\Riem|$, we get
\begin{eqnarray*}
&&C\left(\int\phi^{l\gamma}|\nabla\Ric|^{k\gamma}\right)^\frac1\gamma
\;\le\int\phi^{l-2}|\nabla\phi|^2|\nabla\Ric|^k\\
&&\quad\quad+\,\int\phi^l|\nabla\Ric|^{k-2}|\Ric|^2|\Riem|^2
\,+\,\int\phi^l|\nabla\Ric|^{k-2}|\nabla{X}|^2.
\end{eqnarray*}
H\"older's inequality, combined with Proposition
\ref{GoodRiemEstimate} gives
\begin{eqnarray*}
\left(\int\phi^{l\gamma}|\nabla\Ric|^{k\gamma}\right)^\frac1\gamma
&\le&{C}r^{-2}\int\phi^{l-2}|\nabla\Ric|^k
\,+\,\left(\int\phi^{\Box}|\nabla{X}|^{\frac{nk}{k+n-2}}\right)^{\frac
{n}{k+n-2}}
\end{eqnarray*}
Since $\frac{nk}{k+n-2}<n$, using (\ref{APreliminaryForNablaX}), we
get
\begin{eqnarray*}
\left(\int\phi^{l\gamma}|\nabla\Ric|^{k\gamma}\right)^\frac1\gamma
&\le&{C}r^{-2}\int\phi^{l-2}|\nabla\Ric|^k,
\end{eqnarray*}
which we can iterate to get the stated result for $|\nabla\Ric|$.
Now the equation
$\triangle|\nabla{X}|\ge-C|\nabla\Ric||X|-C|\Riem||\nabla{X}|$,
along with Proposition \ref{LpLemma} (which always works at smooth
points), yields the result for $|\nabla{X}|$. \qed

\subsection{Pointwise curvature regularity}

Here we assume that Proposition \ref{LocalCurvatureBounds} has been
entirely proved at smooth points. The beginning of the proof was
undertaken in the previous section. The rest of the proof,
consisting of an induction argument in dimension, is in the
appendix.

\begin{theorem}\label{RiemBounds}
Assume $B(o,r)$ consists of manifold points. There exists an
$\epsilon_0=\epsilon_0(C_S,n,p)$ and $C=C(C_S,n,p)$ so that
$\int_{B_r} |\Riem|^{\frac{n}{2}} \le \epsilon_0$ implies
\begin{eqnarray*}
\sup_{B(o,r/2)}|\nabla^p\Riem|
\;\le\;Cr^{-p-2}\left(\int_{B(o,r)}|\Riem|^{\frac{n}{2}}\right)^\frac2n.
\end{eqnarray*}
\end{theorem}
\underline{\sl Pf}

First we prove a commutator formula.  If $T$ is any tensor, we have
\begin{eqnarray}
\triangle\nabla\,T&=&T_{,imm} \;=\;T_{,mim}
\,+\,\left({\Riem}_{im**}*T\right)_{,m}\label{GeneralCommutator}\\
&=&T_{,mmi} \,+\,\Riem*\nabla\,T \,+\,{\Riem}_{im**,m}*T
\,+\,{\Riem}_{im**}*T_{,m}\nonumber\\
&=&\nabla\triangle\,T \,+\,\Riem*\nabla\,T
\,+\,\nabla\Ric*T.\nonumber
\end{eqnarray}
Here the stars in the subscript positions of $\Riem$ are meant to
indicate a contraction with various indices of $T$. Replacing $T$
with $\nabla^{p-1}T$, an induction argument gives
\begin{eqnarray*}
&&[\triangle,\nabla^p] \;=\;
\sum_{i=0}^{p-1}\nabla^i\Riem*\nabla^{p-i}
\,+\,\sum_{i=1}^q\nabla\Ric*\nabla^{p-i}.
\end{eqnarray*}
Therefore
\begin{eqnarray*}
&&\triangle\nabla^p\Riem
\;=\;\sum_{i=0}^p\nabla^i\Riem*\nabla^{p-i}\Riem
\,+\,\nabla^{p+2}\Ric,
\end{eqnarray*}
so
\begin{eqnarray*}
\triangle|\nabla^p\Riem| &\ge& -C|\Riem||\nabla^p\Riem|\\
&&\quad\quad-\,C\sum_{i=1}^{p-1}|\nabla^i\Riem||\nabla^{p-i}\Riem|
\,-\,C|\nabla^{p+2}\Ric|.
\end{eqnarray*}
With $u=|\nabla^p\Riem|$, $f=C|\Riem|$, and
$g=C\sum_{i=1}^{p-1}|\nabla^i\Riem||\nabla^{p-i}\Riem|
\,+\,C|\nabla^{p+2}\Ric|$, we get the elliptic inequality
$$
\triangle u \;\ge\; -f\,u \,-\,g
$$
which holds everywhere that $u\ne 0$. Proposition
(\ref{LocalCurvatureBounds}) gives that $f,g\in L^{s'}(B(o,r/2))$
for some $s' > n/2$, so theorem 8.15 of \cite{GT} gives
\begin{eqnarray*}
&&\sup_{B(o,r/4)}|\nabla^p\Riem| \;\le\;
Cr^{-2}\left(\int_{B(o,r/2)}|\nabla^p\Riem|^{\frac{n}{2}}\right)^
\frac2n,
\end{eqnarray*}
and so
\begin{eqnarray*}
&&\sup_{B(o,r/4)}|\nabla^p\Riem| \;\le\;
Cr^{-p-2}\left(\int_{B(o,r)}|\Riem|^{\frac{n}{2}}\right)^\frac2n.
\end{eqnarray*}
Applying this for balls $B(o',r/2)$ with $o'\in\partial B(o,3r/8)$,
we get the final form of the result. \qed

\subsection{Removing curvature singularities, $n\ge6$}\label
{SubsectionRemCurvSingsDim6}

Here we undertake the proof Proposition \ref{LocalCurvatureBounds}
in the cases where $B(o,r)$ has curvature singularities and the
dimension satisfies $n\ge6$. We will make use of our original
elliptic system
\begin{eqnarray}
\triangle\Riem &=&\Riem*\Riem \,+\,\nabla^2\Ric\label
{RiemReferenceEllipticIneq}\\
\triangle\Ric &=&\Riem*\Ric \,+\,\nabla{X}\\
\triangle{X} &=&\Ric*X.
\end{eqnarray}
In addition we will use the formulas
\begin{eqnarray}
\nabla^2X&=&\Riem*X\\
\triangle\nabla{X}&=&\nabla\Ric*X\,+\,\Riem*\nabla{X}\\
\triangle\nabla\Ric&=&\Riem*\nabla\Ric\,+\,\Ric*\nabla\Riem\,+\,\Riem*X.
\end{eqnarray}
Our ultimate goal is to show that $|\Riem|\in{L^k}$ for all $k$
despite the singularities. One must only show that
$|\nabla^2\Ric|\in{L}^{\frac{n}{2}}$, and then Proposition
\ref{LpLemma} theory gives $|\Riem|\in{L}^k$.  Showing that
$|\nabla^2\Ric|\in{L}^{\frac{n}{2}}$ isn't too bad at smooth points,
but with singularities we must use a more round-about route. We
already have $|X|\in{L}^k$ (Proposition \ref{LpLemma}). We can show
$|\nabla{X}|\in{L}^{n}$, so Proposition \ref{LpLemma} gives that
$|\Ric|\in{L}^k$.

Now at this stage we try to get estimates for $\Riem$. The model
case is the real-valued system in divergence form
\begin{eqnarray*}
\triangle{u} &\ge& -fu \,-\,\nabla^ig_i,
\end{eqnarray*}
where one gets that $u\in{L}^k$ provided $f\in{L}^{\frac{n}{2}}$ and
$g_i\in{L}^n$. Abusing both notation and the very notion of
divergence, we consider equation (\ref{RiemReferenceEllipticIneq})
to have a nonhomogeneous term in divergence form, namely
$g_i=\nabla\Ric$. If $|\nabla\Ric|\in{L}^n$ we then expect
$|\Riem|\in{L}^k$. This intuition certainly pans out in the smooth
case, but unfortunately the tool in the singular case, Proposition
\ref{LpLemma}, is not built to handle the divergence term. It is
essentially the divergence structure that we exploit in our
argument, however, so it is likely that some improvements can be
made to Proposition \ref{LpLemma} also.

Technically intricate arguments allow us to play estimates for
$|\nabla\Ric|$ and $|\Riem|$ off of each other; we show that
$|\Riem|\in{L}^p$ implies an improved estimate for $|\nabla\Ric|$,
and this improved estimate in turn lets us bootstrap $|\Riem|$ into
higher $L^p$ spaces.


In we use the following the shorthand notation: if $p$ is a number
we use $p^{-}$ to indicate a variable that may have any value less
than $p$, and $p^{+}$ to indicate a variable that may have any value
greater than $p$.
\begin{lemma}\label{XLemmaLimitedMaybe}
Assume $M$ is a manifold-with-singularities. There exist
$\epsilon_0=\epsilon_0(n,k,C_S)$ and $C=C(n,k,C_S)$ so that
$\int_{B(o,r)}|\Ric|^{\frac{n}{2}}\le\epsilon_0$ implies
\begin{eqnarray*}
\left(\int_{B(o,r/2)}|X|^k\right)^{\frac1k}
&\le&Cr^{\frac{n}{k}-3}\left(\int_{B(o,r)}|R|^{\frac{n}{2}}\right)^
\frac2n.
\end{eqnarray*}
\end{lemma}
\underline{\sl Pf}\\
\indent The proof of this lemma in the smooth case, given by
Proposition \ref{GoodXEstimate}, will carry through provided we can
justify the use of integration by parts in case $n\ge6$.

Assume that $|X|\in{L}_{loc}^{p^{-}}$. Assuming that $\supp\phi$
consists of smooth points, we get
\begin{eqnarray*}
\int\phi^p|X|^p
&=&-p\int\phi^{p-1}R\left<\nabla\phi,X\right>|X|^{p-2}\\
&&-\,(p-2)\int\phi^pR\left<\nabla|X|,X\right>|X|^{p-3}
\,-\,\int\phi^pR\triangle{R}|X|^{p-2}\\
&\le&p\int\phi^{p-1}|R||\nabla\phi||X|^{p-1}\\
&&+\,(p-2)\int\phi^p|R||X|^{p-2}|\nabla{X}|
\,+\,\int\phi^p|R||\triangle{R}||X|^{p-2}\\
&\le&C|R|\left(\int|\nabla\phi|^n\right)^\frac1n\left(\int|X|^{(p-1)
\frac{n}{n-1}}\right)^{\frac{n-1}{n}}\\
&&+\,C|R|\left(\int|\nabla{X}|^2\right)^\frac12\left(\int|X|^{2p-4}
\right)^\frac12
\end{eqnarray*}
Now replace $\phi$ by $\phi\cdot\phi_R$ where $\phi\equiv1$ across
the singularity $o$, and $\phi_R$ is a cutoff function with
$\phi_R\equiv0$ inside $B(o,R)$, $\phi_R\equiv1$ outside $B(o,2R)$,
and $|\nabla\phi_R|\le\frac2R$. Then we can take a limit as
$R\rightarrow0$ assuming first that $(p-1)\frac{n}{n-1}<p$ and
second that $2p-4<p$; it suffices to require $p<4$.

Therefore, assuming $p<4$ so that integration by parts works, we can
get
\begin{eqnarray*}
\int\phi^p|X|^p &\le&C\int|\nabla\phi|^p|R|^p
\,+\,C\int\phi^p|R|^{\frac{p}{2}}|\nabla{X}|^{\frac{p}{2}}.
\end{eqnarray*}
Now we can take a limit as $p\rightarrow4$. Using the Dominated
Convergence Theorem on the right side and Fatou's lemma on the left,
we get that this inequality holds for $p=4$ as well. Now we can
repeat the proof of Proposition \ref{GoodXEstimate}. \qed

\begin{lemma}\label{NablaXLemmaLimited}
Assume $M$ is a manifold-with-singularities. If $2\le{k}\le{n}$
there exist $\epsilon_0=\epsilon_0(n,k,C_S)$ and $C=C(n,k,C_S)$ so
that $\int_{B(o,r)}|\Riem|^{\frac{n}{2}}\le\epsilon_0$ implies
\begin{eqnarray*}
\left(\int_{B(o,r/2)}|\nabla{X}|^k\right)^\frac1k
&\le&Cr^{\frac{n}{k}-4}\left(\int_{B(o,r)}|R|^{\frac{n}{2}}\right)^
\frac2n,
\end{eqnarray*}
irrespective of the presence of singularities.
\end{lemma}
\underline{\sl Pf}
\begin{eqnarray*}
\int\phi^l|\nabla{X}|^k
&=&-l\int\phi^{l-1}|\nabla{X}|^{k-2}\left<\nabla{X},\nabla\phi\otimes
{X}\right>\\
&&-(k-2)\int\phi^l|\nabla{X}|^{k-3}\left<\nabla{X},\nabla|\nabla{X}|
\otimes{X}\right>\\
&&-\int\phi^l|\nabla{X}|^{k-2}\left<\triangle{X},X\right>
\end{eqnarray*}
We use $\nabla^2{X}\;=\;\Riem*X$ and $\triangle{X}\;=\;\Ric*X$ to
get
\begin{eqnarray*}
\int\phi^l|\nabla{X}|^k
&\le&\frac{l^2}{2}\int\phi^{l-2}|\nabla\phi|^2|\nabla{X}|^{k-2}|X|^2
\,+\,\frac12\int\phi^l|\nabla{X}|^k\\
&&+(k-2)\int\phi^l|\nabla{X}|^{k-2}|X|^2|\Riem|\\
&&+\int\phi^l|\nabla{X}|^{k-2}|X|^2|\Ric|,
\end{eqnarray*}
so with $C=C(k,l)$ we have
\begin{eqnarray}
\int\phi^l|\nabla{X}|^k
&\le&C\int\phi^{l-2}|\nabla\phi|^2|\nabla{X}|^{k-2}|X|^2\nonumber\\
&&+\,C\int\phi^l|\nabla{X}|^{k-2}|X|^2|\Riem|. \label{PreNablaXSing}
\end{eqnarray}
If $\supp\phi$ has a singularity, we now show that when $k<2p\le{n}$
this still holds. We get
\begin{eqnarray}
\int\phi^l|\nabla{X}|^k
&\le&C\left(\int\phi^l|\nabla{X}|^k\right)^{\frac{k-2}{k}}\left(\int
\phi^{l-k}|\nabla\phi|^k|X|^k\right)^\frac2k\label{NablaXSing}\\
&&+\,C\left(\int\phi^l|\nabla{X}|^k\right)^{\frac{k-2}{k}}\left(\int
\phi^{\frac{2pl}{2p-k}}|X|^{\frac{2pk}{2p-k}}\right)^{\frac{2p-k}{pk}}
\left(\int|\Riem|^p\right)^\frac1p\nonumber
\end{eqnarray}
Replace $\phi$ in the inequality by $\phi\cdot\phi_R$. Assuming $o$
is a singularity we choose the cutoff function $\phi_R$ with the
following properties: $\phi_R\equiv1$ outside $B(o,2R)$,
$\phi\equiv0$ inside $B(o,R)$ and $|\nabla\phi_R|\le\frac2R$. If we
take a limit as $R\rightarrow0$ (so the cutoff function closes in
around the singularity), we get that
$\phi\cdot\phi_R\rightarrow\phi$ and we can use the dominated
convergence theorem on everything except the integral with the
$\nabla(\phi\phi_R)$, which we analyze separately. With $k<n$ we get
\begin{eqnarray*}
\int(\phi\cdot\phi_R)^{l-k}|\nabla(\phi\cdot\phi_R)|^k|X|^k
&\le&\left(\int|\nabla\phi_R|^n\right)^{\frac{k}{n}}\left(\int_{\supp
\nabla\phi_1}\phi^{\frac{ln}{n-k}}|X|^{\frac{kn}{n-k}}\right)^{\frac
{n-k}{n}}\\
&&+\,\int\phi_R^k|\nabla\phi|^k|X|^{k}.
\end{eqnarray*}
Since $\int|\nabla\phi_R|^n$ is bounded and
$|X|\in{L}^{\frac{nk}{n-k}}$ by Lemma \ref{XLemmaLimitedMaybe}, the
first term on the right side goes to zero as $R\rightarrow0$.
Therefore (\ref{NablaXSing}) holds despite the possible presence of
singularities. Using $p=\frac{n}{2}$ in (\ref{NablaXSing}) now gives
\begin{eqnarray*}
\int\phi^l|\nabla{X}|^k &\le&C\int\phi^{l-k}|\nabla\phi|^k|X|^k
\,+\,C\left(\int\phi^{\frac{nl}{n-k}}|X|^{\frac{nk}{n-k}}\right)^
{\frac{n-k}{n}},
\end{eqnarray*}
where $C=C(k,l)$, and Proposition \ref{LpLemma} yields finally (with
$k<n$)
\begin{eqnarray*}
\int_{B(o,r/2)}|\nabla{X}|^k &\le&Cr^{-k}\int_{B(o,r)}|X|^k.
\end{eqnarray*}
The value of $C$ does not degenerate as $k\nearrow{n}$, so we can
take a limit, using Fatou's lemma on the left side, and get the
result for $k=n$ as well. \qed

\begin{theorem}
Assume $M$ is a manifold-with-singularities. For
$\frac{n}{n-2}\le{k}\le{a}<\infty$, there exist
$\epsilon_0=\epsilon_0(C_S,n,a,k)$ and $C=C(C_S,n,a,k)$ so that
$\int_{B(o,r)}|\Riem|^{\frac{n}{2}}\le\epsilon_0$ implies
\begin{eqnarray*}
\left(\int_{B(o,r/2)}|\Ric|^a\right)^\frac1a
&\le&Cr^{\frac{n}{a}-\frac{n}{k}}\left(\int_{B(o,r)}|\Ric|^k\right)^
\frac1k,
\end{eqnarray*}
irrespective of the presence of singularities.
\end{theorem}
\indent\underline{\sl Pf}\\
With
\begin{eqnarray}
\triangle\Ric &=&\Riem*\Ric \,+\,\nabla{X}\\
\triangle|\Ric| &\ge&-|\Riem||\Ric| \,-\,|\nabla{X}|,
\end{eqnarray}
and since $|\nabla{X}|\in{L}^{\frac{n}{2}}$, we can use Proposition
\ref{LpLemma} to get that $|\Ric|\in{L}^k$ for all $k<\infty$.  We
get the local estimates
\begin{eqnarray}
\left(\int_{B(o,r/2)}|\Ric|^a\right)^\frac1a
&\le&Cr^{\frac{n}{a}-\frac{n}{k}}\left(\int_{B(o,r)}|\Ric|^k\right)^
\frac1k.
\end{eqnarray}
\qed

Note that the hypotheses of the following technical lemma hold
because we have independently proven Theorem \ref{NablaPRiemTheorem}
and Theorem \ref{RiemBounds} in the smooth case. The proof is just a
more involved version of the proof of \ref{XLemmaLimitedMaybe}.
\begin{technicallemma} \label{Nabla2RicImprovements}
Assuming $|\nabla^p\Ric|=o(r^{-2-p})$ near singularities, then
$|\nabla^2\Ric|\in{L}^{\frac{n}{3}^{-}}$ and
$|\nabla\Ric|\in{L}^{\frac23n^{-}}$.
\end{technicallemma}
\underline{\sl Pf}\\
\indent First, we know that
$|\nabla^4\Ric|\in{L}^{\frac{n}{6}^{-}}$. Thus assuming
$|\nabla^2\Ric|\in{L}^{p^{-}}$ we have
\begin{eqnarray*}
\int\phi^2|\nabla^3\Ric|^k
&=&-2\int\phi|\nabla^3\Ric|^{k-2}\left<\nabla\phi\otimes\nabla^2\Ric,
\nabla^3\Ric\right>\\
&&-(k-2)\int\phi^2|\nabla^3\Ric|^{k-3}\left<\nabla|\nabla^3\Ric|
\otimes\nabla^2\Ric,\nabla^3\Ric\right>\\
&&-\int\phi^2|\nabla^3\Ric|^{k-3}\left<\nabla^2\Ric,\triangle\nabla^2
\Ric\right>\\
&\le&c\int\phi|\nabla\phi||\nabla^3\Riem|^{k-1}|\nabla^2\Ric|\\
&&+\,c\left(\int|\nabla^3\Ric|^{k^{-}}\right)^{\frac{k-2}{k}\frac{n}
{n-6}^{+}}\left(\int|\nabla^2\Ric|^{\frac{kn}{2n-6k}^{+}}\right)^
{\frac{2n-6k}{k(n-6)}^{-}}\\
&&+\,c\int|\nabla^4\Ric|^{\frac{n}{6}^{-}}
\end{eqnarray*}
which holds across singularities provided $\frac{nk}{n+k}\le{p}$.
Thus $\frac{nk}{2n-6k}<p$ gives $|\nabla^3\Ric|\in{L}^k$.  Since we
can always choose $p>\frac{n}{5}$, we have that
\begin{eqnarray*}
k<\frac{2np}{n+6p} \quad \implies \quad |\nabla^3\Ric|\in{L}^k.
\end{eqnarray*}

We do the same thing for $|\nabla^2\Ric|$. Assume
$|\nabla\Ric|\in{L}^{q^{-}}$. We get
\begin{eqnarray*}
\int\phi^2|\nabla^2\Ric|^k
&\le&2\int\phi|\nabla\phi||\nabla^2\Ric|^{k-1}|\nabla\Ric|\\
&&+\,c\left(\int|\nabla^2\Ric|^{k^{-}}\right)^{\frac{k-2}{k}\frac{m}
{m-1}^{+}}\left(\int|\nabla\Ric|^{\frac{km}{2m-k}^{+}}\right)^{\frac
{2m-k}{k(m-1)}^{-}}\\
&&+\,(k-1)\int\phi^m|\nabla^3\Ric|^{m^{-}}
\end{eqnarray*}
This holds across singularities if $\frac{nk}{n+k}<q$. If
$|\nabla^3\Ric|\in{L}^{m^{-}}$ and $\frac{km}{2m-k}<q$, then
$|\nabla^2\Ric|\in{L}^k$.  Assuming $m<n$, we get
$k<\frac{2mq}{m+q}$ implies $|\nabla^2\Ric|\in{L}^k$.

We do the same thing for $|\nabla\Ric|$. Assume
$|\Ric|\in{L}^{q^{-}}$. We get
\begin{eqnarray*}
\int\phi^2|\nabla\Ric|^k
&\le&2\int\phi|\nabla\phi||\nabla\Ric|^{k-1}|\Ric|\\
&&+\,c\left(\int|\nabla\Ric|^{k^{-}}\right)^{\frac{k-2}{k}\frac{m}
{m-1}^{+}}\left(\int|\Ric|^{\frac{km}{2m-k}^{+}}\right)^{\frac{2m-k}{k
(m-1)}^{-}}\\
&&+\,(k-1)\int\phi^m|\nabla^2\Ric|^{m^{-}}
\end{eqnarray*}
This holds across singularities if $\frac{nk}{n+k}<q$. If
$|\nabla^2\Ric|\in{L}^{m^{-}}$ and $\frac{km}{2m-k}<q$, then
$|\nabla\Ric|\in{L}^k$.  Assuming $m<n$, we get $k<\frac{2mq}{m+q}$
implies $|\nabla\Ric|\in{L}^k$.

The result of these three inequalities is that
\begin{eqnarray}
&&|\nabla^3\Ric|\in{L}^{m^{-}},\;{\rm
and}\quad|\nabla\Ric|\in{L}^{q^{-}} \implies
|\nabla^2\Ric|\in{L}^{\frac{2mq}{m+q}^{-}}\label{ImprovementNabla2Ric}\\
&&|\nabla^2\Ric|\in{L}^{r^{-}} \implies
|\nabla\Ric|\in{L}^{2r^{-}}\label{ImprovementNablaRic}\\
&&|\nabla^2\Ric|\in{L}^{r^{-}} \implies
|\nabla^3\Ric|\in{L}^{\frac{2nr}{n+6r}}\label{ImprovementNabla3Ric}
\end{eqnarray}
Fixing $m$ and iterating, we continue to get increases in $q$ up
until $q=3m$. Thus we get $|\nabla\Ric|\in{L}^{3m^{-}}$ and
$|\nabla^2\Ric|\in{L}^{\frac{3m}{2}^{-}}$.  Then letting $m$ vary
and iterating, we get improvements up until $m=\frac{2n}{9}$.

Therefore $|\nabla\Ric|\in{L}^{\frac23n^{-}}$,
$|\nabla^2\Ric|\in{L}^{\frac13n^{-}}$,
$|\nabla^3\Ric|\in{L}^{\frac29n^{-}}$.

\qed

The next lemma sets up the possibility of using integration by parts
across singularities, but does not give any particular bound for
$L^p(|\Riem|)$.
\begin{technicallemma} \label{Nabla2RicLimited}
Assuming the above lemma, we have $|\Riem|\in{L}^{p}$ and
$|\nabla\Ric|\in{L}^p$ for all $p$.
\end{technicallemma}
\underline{\sl Pf}\\
\noindent We can use the improvement integral bounds of
$|\nabla^2\Ric|$ to our advantage. Sobolev's inequality gives
\begin{eqnarray*}
C\left(\int\phi^{2k}|\Riem|^{k\gamma}\right)^\frac1\gamma
&\le&\int|\nabla\phi|^2|\Riem|^k
\,+\,\int\phi^2|\Riem|^{k-2}|\nabla|\Riem||^2.
\end{eqnarray*}
Choosing $k=\frac{n-2}{2}$ so integration by parts works across
singularities (by Lemma \ref{IntegrationByParts}), so we get
\begin{eqnarray*}
C\left(\int\phi^{2k}|\Riem|^{k\gamma}\right)^\frac1\gamma
&\le&\int|\nabla\phi|^2|\Riem|^k \,+\,\int\phi^2|\Riem|^{k+1}
\,+\,\int\phi^2|\Riem|^{k-1}|\nabla^2\Ric|.
\end{eqnarray*}
Using that $\int|\Riem|^{\frac{n}{2}}$ is small we get
\begin{eqnarray*}
C\left(\int\phi^{2k}|\Riem|^{\frac{n}{2}}\right)^\frac1\gamma
&\le&\int|\nabla\phi|^2|\Riem|^{\frac{n-2}{2}}
\,+\,C\left(\int|\nabla^2\Ric|^{\frac{n}{4}}\right)^\frac4n\left(\int|
\Riem|^{\frac{n}{2}}\right)^{\frac{n-4}{n}}
\end{eqnarray*}
Since $|\nabla^2\Ric|\in{L}^{\frac{n}{4}}$  we get
\begin{eqnarray}
C\left(\int\phi^{2k}|\Riem|^{\frac{n}{2}}\right)^\frac1\gamma
&\le&\int|\nabla\phi|^2|\Riem|^{\frac{n-2}{2}}
\,+\,\left(\int|\nabla^2\Ric|^{\frac{n}{4}}\right)^{\frac{2(n-2)}{n}}.
\label{ImprovedRiem}
\end{eqnarray}
Since $\Riem\in{L}^{n/2}$ this holds across singularities. Using an
argument similar to Theorem 5.8 of \cite{BKN} we can get that
$\int_{B_r}|\Riem|^{\frac{n}{2}}$ decays like
$\left(\int_{B_r}|\nabla^2\Ric|^{\frac{n}{4}}\right)^2$ (the
argument needed here is given in detail in Lemma
\ref{UhlenbeckPrelimEst}). Since
$|\nabla^2\Ric|\in{L}^{\frac{n}{3}^{-}}$, we get that
\begin{eqnarray*}
\int_{B_r}|\Riem|^{\frac{n}{2}} &=&O(r^{\frac{n}{2}^{-}}).
\end{eqnarray*}
Using this in conjunction with Theorem \ref{RiemBounds} gives that
$|\nabla^s\Riem|=O(r^{(-s-1)^{-}})$ near singularities. Note that
all this fails in the case $n<6$, for in that case the use of the
Sobolev inequality that began the discussion would be unavailable to
us. This means that $|\Riem|\in{L}^{n^{-}}$,
$|\nabla\Riem|\in{L}^{\frac{n}{2}^{-}}$,
$|\nabla^2\Riem|\in{L}^{\frac{n}{3}^{-}}$,
$|\nabla^3\Riem|\in{L}^{\frac{n}{4}^{-}}$ etc.

Now we return to (\ref{ImprovementNabla2Ric}),
(\ref{ImprovementNablaRic}), and (\ref{ImprovementNabla3Ric}) from
above.  We now have $|\nabla^3\Ric|\in{L}^{\frac{n}{4}}$ so we can
expect some improvements. We initially have that
$|\nabla\Ric|\in{L}^{\frac23n^{-}}$, as we got from the last
theorem. We can iterate up until $q=\frac34n$, so that
$|\nabla^2\Ric|\in{L}^{\frac{3n}{8}^{-}}.$ We get therefore that
\begin{eqnarray*}
&&\left(\int_{B_r}|\nabla^2\Ric|^{\frac{n}{4}}\right)^2=O(r^{\frac{2n}
{3}^{-}})
\end{eqnarray*}
and so (\ref{ImprovedRiem}) implies that
$\int_{B(o,r)}|\Riem|^{\frac{n}{2}}$ decays like
$O(r^{\frac{2n}{3}^{-}}).$ Running through the above argument again,
that $|\nabla^s\Riem|=O(r^{-(s+\frac23)^{-}})$.  This yields
actually that $|\nabla\Riem|\in{L}^{\frac35n^{-}}$, and with
\begin{eqnarray*}
\triangle\nabla\Ric &=&\nabla\Riem*\Ric \,+\,\Riem*\nabla\Ric,
\end{eqnarray*}
Now Proposition \ref{LpLemma} implies that $|\nabla\Ric|\in{L}^p$
for all $p$.

Then (\ref{ImprovementNabla2Ric}) implies now that
$|\nabla^2\Ric|\in{L}^{\frac{6}{11}n^{-}}$. But then Proposition
\ref{LpLemma} applied to $\triangle\Riem=\Riem*\Riem+\nabla^2\Ric$
implies that $|\Riem|\in{L}^p$ for all $p$. \qed

\begin{proposition} \label{PropRiemNablaRicLpBounds}
Assuming $\int_{B(o,r)}|\Riem|^{\frac{n}{2}}<\epsilon_0$ we get
\begin{eqnarray*}
\left(\int_{B(o,r/2)}|\Riem|^a\right)^\frac1a
&\le&Cr^{\frac{n}{a}-2}\left(\int_{B(o,r)}|\Riem|^{\frac{n}{2}}\right)
^\frac2n\\
\left(\int_{B(o,r/2)}|\nabla\Ric|^a\right)^\frac1a
&\le&Cr^{\frac{n}{a}-3}\left(\int_{B(o,r)}|\Ric|^{\frac{n}{2}}\right)^
\frac2n
\end{eqnarray*}
for all $a>0$, regardless of the presence of singularities.
\end{proposition}
\underline{\sl Pf}\\
\indent We try to gain estimates for $L^p(|\nabla\Ric|)$.
\begin{eqnarray*}
C\left(\int\phi^{l\gamma}|\nabla\Ric|^{k\gamma}\right)^\frac1\gamma
&\le&\int\phi^{l-2}|\nabla\phi|^2|\nabla\Ric|^k
\,+\,\int\phi^l|\nabla\Ric|^{k-2}|\nabla|\nabla\Ric||^2\\
\int\phi^l|\nabla\Ric|^{k-2}|\nabla|\nabla\Ric||^2
&\le&\int\phi^{l-2}|\nabla\phi|^2|\nabla\Ric|^k
\,-\,\int\phi^l|\nabla\Ric|^{k-2}\left<\nabla\Ric,\triangle\nabla\Ric
\right>.
\end{eqnarray*}
Using a commutator formula on the last term, we get
\begin{eqnarray*}
&&\int\phi^l|\nabla\Ric|^{k-2}|\nabla|\nabla\Ric||^2\\
&&\quad\le\;\int\phi^{l-2}|\nabla\phi|^2|\nabla\Ric|^k
\,+\,\int\phi^l|\nabla\Ric|^k|\Riem|\\
&&\quad\quad+\,l\int\phi^{l-1}|\nabla\phi||\nabla\Ric|^{k-1}|\triangle
\Ric|
\,+\,\int\phi^l|\nabla\Ric|^{k-2}|\triangle\Ric|^2\\
&&\quad\quad+\,(k-2)\int\phi^l|\nabla\Ric|^{k-3}\left<\nabla\Ric,
\nabla|\nabla\Ric|\otimes\triangle\Ric\right>.
\end{eqnarray*}
Therefore
\begin{eqnarray*}
&&\int\phi^l|\nabla\Ric|^{k-2}|\nabla|\nabla\Ric||^2\\
&&\quad\le\;C\int\phi^{l-2}|\nabla\phi|^2|\nabla\Ric|^k
\,+\,C\int\phi^l|\nabla\Ric|^k|\Riem|\\
&&\quad\quad+\,C\int\phi^l|\nabla\Ric|^{k-2}|\triangle\Ric|^2
\end{eqnarray*}
and so
\begin{eqnarray*}
C\left(\int\phi^{l\gamma}|\nabla\Ric|^{k\gamma}\right)^\frac1\gamma
&\le&\int\phi^{l-2}|\nabla\phi|^2|\nabla\Ric|^k
\,+\,\int\phi^l|\nabla\Ric|^{k-2}|\triangle\Ric|^2
\end{eqnarray*}
\begin{eqnarray}
&&C\left(\int\phi^{l\gamma}|\nabla\Ric|^{k\gamma}\right)^\frac1\gamma
\;\le\int\phi^{l-2}|\nabla\phi|^2|\nabla\Ric|^k\nonumber\\
&&\quad\quad+\,\int\phi^l|\nabla\Ric|^{k-2}|\Ric|^2|\Riem|^2
\,+\,\int\phi^l|\nabla\Ric|^{k-2}|\nabla{X}|^2\label
{NablaRicSobolevDecent}
\end{eqnarray}

\noindent To continue we must estimate $\int|\Riem|^a$ locally:
\begin{eqnarray*}
\left(\int\phi^{k\gamma}|\Riem|^{k\gamma}\right)^\frac1\gamma
&\le&\int\phi^{l-2}|\nabla\phi|^2|\Riem|^k
\,+\,\int\phi^l|\Riem|^{k-2}|\nabla|\Riem||^2
\end{eqnarray*}
\begin{eqnarray*}
&&C\int\phi^l|\Riem|^{k-2}|\nabla|\Riem||^2\\
&&\quad\le\;\int\phi^{l-2}|\nabla\phi|^2|\Riem|^k
\,-\,\int\phi^l|\Riem|^{k-2}\left<\Riem,\triangle\Riem\right>\\
&&C\int\phi^l|\Riem|^{k-2}|\nabla|\Riem||^2\\
&&\quad\le\;\int\phi^{l-2}|\nabla\phi|^2|\Riem|^k
\,+\,\int\phi^l|\Riem|^{k+1}
\,-\,\int\phi^l|\Riem|^{k-2}\left<\Riem,\nabla^2\Ric\right>
\end{eqnarray*}
\begin{eqnarray*}
&&C\int\phi^l|\Riem|^{k-2}|\nabla|\Riem||^2\\
&&\quad\le\;\int\phi^{l-2}|\nabla\phi|^2|\Riem|^k
\,+\,\int\phi^l|\Riem|^{k+1}
\,+\,\int\phi^l|\Riem|^{k-2}\left<\Riem,\nabla\phi\otimes\nabla\Ric \right>\\
&&\quad\quad+\,\int\phi^l|\Riem|^{k-2}\left<\Riem,\nabla|\Riem|\otimes
\nabla\Ric\right>
\,+\,\int\phi^l|\Riem|^{k-2}\left<\nabla\Ric,\nabla\Ric\right>
\end{eqnarray*}
\begin{eqnarray*}
&&C\int\phi^l|\Riem|^{k-2}|\nabla|\Riem||^2\\
&&\quad\le\;\int\phi^{l-2}|\nabla\phi|^2|\Riem|^k
\,+\,\int\phi^l|\Riem|^{k+1}
\,+\,\int\phi^l|\Riem|^{k-2}|\nabla\Ric|^2
\end{eqnarray*}
Doing the same integration-by-parts on the last term we get finally
\begin{eqnarray*}
&&C\int\phi^l|\Riem|^{k-2}|\nabla|\Riem||^2\\
&&\quad\le\;\int\phi^{l-2}|\nabla\phi|^2|\Riem|^k
\,+\,\int\phi^l|\Riem|^{k+1}
\,+\,\int\phi^l|\Riem|^{k-2}\left<\Ric,\triangle\Ric\right>
\end{eqnarray*}
Altogether therefore,
\begin{eqnarray*}
\left(\int\phi^{k\gamma}|\Riem|^{k\gamma}\right)^\frac1\gamma
&\le&\int\phi^{l-2}|\nabla\phi|^2|\Riem|^k
\,+\,\int\phi^l|\Riem|^{k+1} \,+\,\int\phi^l|\Riem|^{k-1}|\nabla{X}|
\end{eqnarray*}
Using H\"older's inequality and Lemma \ref{NablaXLemmaLimited} on
the last term, we get that
\begin{eqnarray*}
\left(\int\phi^{k\gamma}|\Riem|^{k\gamma}\right)^\frac1\gamma
&\le&\int\phi^{l-2}|\nabla\phi|^2|\Riem|^k, \label{LPNormaKnowledge}
\end{eqnarray*}
which, using Lemma \ref{Nabla2RicLimited}, holds across
singularities for all $k$. Iterating this gives
\begin{eqnarray*}
\left(\int_{B(o,r)}|\Riem|^a\right)^\frac1a
&\le&r^{\frac{n}{a}-2}\left(\int_{B(o,r/2)}|\Riem|^{\frac{n}{2}}
\right)^\frac2n.
\end{eqnarray*}

Returning to (\ref{NablaRicSobolevDecent}), we get
\begin{eqnarray*}
&&C\left(\int\phi^{l\gamma}|\nabla\Ric|^{k\gamma}\right)^\frac1\gamma\\
&&\quad\le\int\phi^{l-2}|\nabla\phi|^2|\nabla\Ric|^k
\,+\,\left(\int\phi^l|\nabla\Ric|^{k\gamma}\right)^{\frac{k-2}{k}
\frac1\gamma}\left(\int|\nabla{X}|^{\frac{nk}{n+k-2}}\right)^{\frac{2
(k+n-2)}{kn}}\\
&&\quad\quad+\,\left(\int\phi^{l\gamma}|\nabla\Ric|^{k\gamma}\right)^
{\frac{k-2}{k}\frac1\gamma}\left(\int|\Ric|^{k\gamma}\right)^{\frac2k
\frac1\gamma}\left(\int|\Riem|^n\right)^\frac2n
\end{eqnarray*}
knowing what we do about $L^{n}(|\Riem|)$, $L^p(|\Ric|)$, and
$L^p(|\nabla{X}|)$, we iterate to get
\begin{eqnarray*}
C\left(\int_{B(o,r/2)}|\nabla\Ric|^{a}\right)^\frac1a
&\le&r^{\frac{n}{a}-\frac{n}{2}}\left(\int_{B(o,r)}|\nabla\Ric|^2
\right)^\frac12
\,+\,r^{\frac{n}{a}-3}\left(\int_{B(o,r)}|\Ric|^{\frac{n}{2}}\right)^
\frac2n.
\end{eqnarray*}
We can easily estimate $\int|\nabla\Ric|^2$ in terms of $|\Ric|$, so
we can get
\begin{eqnarray*}
\left(\int_{B(o,r/2)}|\nabla\Ric|^{a}\right)^\frac1a
&\le&Cr^{\frac{n}{a}-3}\left(\int_{B(o,r)}|\nabla\Ric|^{\frac{n}{2}}
\right)^\frac2n
\end{eqnarray*}
\qed

\subsection{Removing curvature singularities, $n=4$} \label
{SubsectionUhlenbeck}

In this section we prove that for some $s>0$, $|\Riem|=O(r^{-2+s})$
in dimension 4, where $r$ indicates distance to a singularity.
Although we do not get specific bounds of the sort in Theorem
\ref{LocalCurvatureBounds}, this result is enough to prove the full
removable singularity theorem in Section \ref{WeakCompactness}. Some
parts of the argument are glossed over here; a complete argument can
be found in the thesis of the second author \cite{Web1}.

In dimension 4 the situation is unfortunately less straightforward
than in higher dimensions. Roughly speaking our coupled elliptic
system has the form $\triangle u \ge -fu - g$, where $u\ge0$ is some
curvature quantity. In dimension $4$ the hypothesis that
$f,g\in{L}^{\frac{n}{2}}$ is insufficient for a purely analytical
argument to remove a point singularity.  The counterexample is
$u=-r^{-2}(\log{r})^{-1}$, for which Sibner's lemma also fails.

We look again to the geometry of our manifolds to provide us
additional input. Uhlenbeck's 1982 paper on Yang-Mills connections
introduced what has become a standard technique here, which we
briefly review. After a choice of gauge (local coordinates) the
connection can be written $D=d+A$, with $A$ being an
$\mathfrak{so}(n)$-valued 1-form, and the curvature $F$, by
definition just $D\circ{D}$, can be written $F=DA-\frac12[A,A]$.
Uhlenbeck used the implicit function theorem to show that in the
annulus, if the metric is almost flat and the gauge is chosen so $A$
is small, the gauge can be slightly modified to make $D^{*}A=0$. A
gauge in which this holds is called a Hodge gauge. To get better
control on $L^2(|F|)$, one estimates on the annulus $\Omega$,
\begin{eqnarray*}
\int_\Omega|F|^2 &=&\int_\Omega\left<DA-\frac12[A,A],F\right>\\
&=&\int_\Omega\left<a,D^{*}F\right> \,+\,C\int_\Omega|A|^2|F|
\,+\,Boundary\;Terms.
\end{eqnarray*}
Working in a Hodge gauge has the advantage of making certain
estimates involving $A$ possible; for instance the $\int|A|^2|F|$
term can be estimated. The $D^{*}F$ term is ordinarily
uncontrollable, but whatever advantage one can squeeze out here may
improve the estimate for $\int|F|^2$. For Yang-Mills connections
$D^{*}F=0$ by definition; this is also true in the Einstein case. In
general the second Bianchi identity gives only that $D^{*}F$ is a
combination of $\nabla\Ric$ terms, so in principle better control
over $L^2(|F|)$ can come from better control over
$L^p(|\nabla\Ric|)$. This was essentially the method of \cite{TV1},
where they were able to get improved estimates for $|\Ric|$, and
then for $|\nabla\Ric|$. Assuming that a good $L^p$ estimate for
$|\nabla\Ric|$ is somehow achieved, one gets $L^2(|F|)$ estimates on
the full punctured disk by estimating on successively smaller
annuli, piecing together the boundary terms, and showing that the
residue (the inner boundary term on the shrinking annuli) vanishes.
For details, see \cite{Uhl}, \cite{Tia1}, \cite{TV1}.

In the case most similar to ours, Theorem 6.4 of \cite{TV1}, better
control on $L^2(|\nabla\Ric|)$ is achieved using an {\sl improved
Kato inequality} for the Ricci curvature, which yields an improved
elliptic inequality.  Their inequality relied on the K\"ahler metric
having constant scalar curvature, so their particular estimates are
unavailable to us. In the proof below we essentially take advantage
of the holomorphicity of $X$ to recover some information about the
irreducible $U(n)$ decomposition of derivatives of curvature tensor.
However we only partially recover an improved Kato inequality and
more effort is needed to achieve something useful. Although our
method of proof is standard, we run through it again because the
value of the constant actually turns out to be important.

Assume $V$ is a complex vector space. Let $\mathcal{A}$ be the space
of tensors $A_{i\bar{j}kl}$ of type
$V\otimes\overline{V}\otimes{V}\otimes{V}$ such that $A$ is
trace-free in the first two positions and symmetric in the first and
third positions; that is $\Sigma_s{A}_{s\bar{s}kl}=0$ and
$A_{i\bar{j}kl}=A_{k\bar{j}il}$. Let $\mathcal{B}$ be the space of
tensors $B_{i\bar{j}k}$ of type $V\otimes\overline{V}\otimes{V}$
that are trace-free in the first two positions.
\begin{lemma}
Assume $V$ has complex dimension $m$. Let
$\left<\,,\,\right>:\mathcal{A}\times\overline{\mathcal{B}}\rightarrow
{V}$ denote the trace in the first three positions. Then
$|\left<a,B\right>|^2\le\frac{m-1}{2m}|A|^2|B|^2$ when
$B_{i\bar{s}s}\ne0$.
\end{lemma}
\underline{\sl Pf}\\
\noindent Restricting ourselves to tensors of unit norm, and using
Lagrange's multiplier method, one finds
\begin{eqnarray}
\left<\left<\widetilde{A},B\right>,\left<a,B\right>\right> &=&
a\left<a,\widetilde{A}\right>\label{FirstLagrange}\\
\left<\left<a,\widetilde{B}\right>,\left<a,B\right>\right> &=&
b\left<B,\widetilde{B}\right>,\label{secondLagrange}
\end{eqnarray}
for $\widetilde{A}\in\mathcal{A}$, $\widetilde{B}\in\mathcal{B}$
arbitrary. Clearly $a=b=|\left<a,B\right>|^2$. Letting $\lambda$ be
the vector $\lambda=\frac{1}{|\left<a,B\right>|^2}\left<a,B\right>$,
(\ref{FirstLagrange}) and (\ref{secondLagrange}) can be written
\begin{eqnarray*}
\left<\widetilde{A},B\otimes\lambda\right>
&=&\left<\widetilde{A},A\right>\\
\left<\left<a,\lambda\right>,\widetilde{B}\right>
&=&\left<B,\widetilde{B}\right>.
\end{eqnarray*}
This means that, with $\pi_1$ the projection onto $\mathcal{A}$ and
$\pi_2$ the projection onto $\mathcal{B}$, $B$ satisfies
\begin{eqnarray}
B&=&\pi_2\left<\pi_1(B\otimes\lambda),\lambda\right>.\label
{BKatoEquation}
\end{eqnarray}
For arbitrary $\widetilde{A}\in\mathcal{A}$,
$\widetilde{B}\in\mathcal{B}$,
\begin{eqnarray*}
\pi_1(\widetilde{A})
&=&\frac12\left(\widetilde{A}_{i\bar{j}kl}\,+\,\widetilde{A}_{k\bar{j}
il}\right)
\,-\,\frac{1}{2m}\delta_{i\bar{j}}\left(\widetilde{A}_{s\bar{s}kl}\,+
\,\widetilde{A}_{k\bar{s}sl}\right)\\
\pi_2(\widetilde{B}) &=&\widetilde{B}_{i\bar{j}k}
\,-\,\frac1m\delta_{i\bar{j}}\widetilde{B}_{s\bar{s}k}.
\end{eqnarray*}
Then we compute
\begin{eqnarray*}
\pi_2\left<\pi_1(B\otimes\lambda),\lambda\right>
&=&\frac12\left(B_{i\bar{j}k}+B_{k\bar{j}i}\right)|\lambda|^2
\,-\,\frac{1}{2m}\delta_{i\bar{j}}B_{k\bar{s}s}|\lambda|^2.
\end{eqnarray*}
Tracing both sides of (\ref{BKatoEquation}) in $j$, $k$, gives
\begin{eqnarray}
B_{i\bar{s}s}\left(2-|\lambda|^2+\frac1m|\lambda|^2\right) \;=\;0,
\end{eqnarray}
so either $B$ is trace-free in the second two variables, or
$|\lambda|^2=\frac{2m}{m-1}$. \qed

\begin{proposition}[Improved Kato Inequality] \label{ImprovedKatoLemma}
Let $M$ be an extremal K\"ahler manifold of complex dimension $m$
and of nonconstant scalar curvature. Denote by
$E_{i\bar{j}}=\Ric_{i\bar{j}}-\frac1mh_{i\bar{j}}R$ the trace-free
Ricci tensor. Then
\begin{eqnarray}
2|\nabla|\nabla{E}||^2 &\le& \frac{m-1}{2m}|\nabla^2E|^2
\,+\,|\overline\nabla\nabla{E}|^2,
\end{eqnarray}
where we denote $\nabla^2E=E_{i\bar{j},kl}$ and
$\overline\nabla\nabla{E}=E_{i\bar{j},k\bar{l}}$.
\end{proposition}
\underline{\sl Pf}\\
\noindent Adopting the notation from above, we have
$\nabla^2E\in\mathcal{A}$ and $\nabla{E}\in\mathcal{B}$. Therefore
\begin{eqnarray}
|\left<\nabla^2E,\nabla{E}\right>|^2\le\frac{m-1}{2m}|\nabla^2E|^2|
\nabla{E}|^2
\end{eqnarray}
The result follows from the identity
\begin{eqnarray*}
\nabla|\nabla{E}|^2 &=&\left<\nabla^2E,\nabla{E}\right>
\,+\,\left<\nabla{E},\overline\nabla\nabla{E}\right>.
\end{eqnarray*}
\qed

\begin{lemma}[Improved elliptic inequality]\label {ImprovedEllipticUsingKato}
If $|\alpha-1-\delta|<\frac{\sqrt{32}}{5}$, then
\begin{eqnarray*}
&&\int\phi^2|\nabla{E}|^{\alpha}\triangle|\nabla{E}|^{1-\delta}\\
&&\quad\ge-\frac12(1-\delta)\delta{C}\int|\nabla\phi|^2|\nabla{E}|^{1-
\delta+\alpha}
\,-\,\frac12(1-\delta)\delta{C}\int\phi^2|\Ric||\nabla{E}|^{1-\delta+
\alpha}\\
&&\quad\quad+\frac12(1-\delta)\int\phi^2|\nabla{E}|^{-1-\delta+\alpha}
\left(\left<\triangle\nabla{E},\nabla{E}\right>
\,+\,\left<\nabla{E},\overline\triangle\nabla{E}\right>\right).
\end{eqnarray*}
\end{lemma}
\underline{\sl Pf}\\
\indent Using the improved Kato inequality and, and setting
$\eta=\frac{m-1}{2m}$, we get
\begin{eqnarray}
\triangle|\nabla{E}|^{1-\delta}
&=&\frac12(1-\delta)|\nabla{E}|^{-1-\delta}\left((1+\delta)\eta|
\nabla^2E|^2
\,-\,\delta|\overline\nabla\nabla{E}|^2\right)\nonumber\\
&&+\,\frac12(1-\delta)|\nabla{E}|^{-1-\delta}\left(\left<\triangle
\nabla{E},\nabla{E}\right>\,+\,\left<\nabla{E},\overline\triangle
\nabla{E}\right>\right)\label{LaplacianNablaE}
\end{eqnarray}
We want $(1+\delta)\eta|\nabla^2E|^2
\,-\,\delta|\overline\nabla\nabla{E}|^2\ge0,$ though this does not
seem possible in the pointwise sense. We will have better luck after
integration however. Using integration by parts and a commutator
formula, we get
\begin{eqnarray*}
&&\int\phi^2|\nabla{E}|^\beta|\overline\nabla\nabla{E}|^2\\
&&\quad\le\int\phi^2|\nabla{E}|^\beta|\nabla^2E|^2\\
&&\quad\quad\quad+\,2\int\phi|\nabla\phi||\nabla{E}|^{1+\beta}|
\nabla^2E|
\,+\,2\int\phi|\nabla\phi||\nabla{E}|^{1+\beta}|\overline\nabla\nabla {E}|\\
&&\quad\quad\quad+\,|\beta|\int\phi^2|\nabla{E}|^{\beta}|\nabla^2E||
\nabla|\nabla{E}||
\,+\,|\beta|\int\phi^2|\nabla{E}|^{\beta}|\overline\nabla\nabla{E}||
\nabla|\nabla{E}||\\
&&\quad\quad\quad+\,3\int\phi^2|\nabla{E}|^{\beta+2}|\Ric|
\end{eqnarray*}
Then with
$|\nabla|\nabla{E}||\le\frac{1}{\sqrt{2}}|\nabla^2E|+\frac{1}{\sqrt
{2}}|\overline\nabla\nabla{E}|$ and assuming that
$|\beta|<\frac{\sqrt{32}}{5}$, we get
\begin{eqnarray*}
&&\left(1-\frac{5|\beta|}{4\sqrt{2}}\right)\int\phi^2|\nabla{E}|^
\beta|\overline\nabla\nabla{E}|^2\\
&&\quad\le\left(2+\frac{33|\beta|}{\sqrt{2}}\right)\int\phi^2|\nabla
{E}|^\beta|\nabla^2E|^2\\
&&\quad\quad\quad+\,3\int\phi^2|\nabla{E}|^{\beta+2}|\Ric|
\,+\,\left(1+\frac{8\sqrt{2}}{|\beta|}\right)\int|\nabla\phi|^2|\nabla
{E}|^{2+\beta}.
\end{eqnarray*}
In fact it is only really necessary that $|\beta|<\sqrt{2}$, but
this method does not allow for arbitrary $\beta$.  With
$|\alpha|<\frac{\sqrt{32}}{5}-|1+\delta|$ we therefore get
\begin{eqnarray*}
&&\int\phi^2|\nabla{E}|^{\alpha}\triangle|\nabla{E}|^{1-\delta}\\
&&\quad\ge\frac12(1-\delta)\left(1-\eta(1+\delta)-\delta{C}\right)\int
\phi^2|\nabla{E}|^{-1-\delta+\alpha}|\nabla^2E|^2\\
&&\quad\quad-\frac12(1-\delta)\delta{C}\int|\nabla\phi|^2|\nabla{E}|^
{1-\delta+\alpha}
\,-\,\frac12(1-\delta)\delta{C}\int\phi^2|\Ric||\nabla{E}|^{1-\delta+
\alpha}\\
&&\quad\quad+\frac12(1-\delta)\int\phi^2|\nabla{E}|^{-1-\delta+\alpha}
\left(\left<\triangle\nabla{E},\nabla{E}\right>
\,+\,\left<\nabla{E},\overline\triangle\nabla{E}\right>\right).
\end{eqnarray*}
The first term is positive when $\delta$ is sufficiently small. \qed

The next lemma shows how to use this improved elliptic inequality.
\begin{lemma}
Assume $\supp\phi$ consists of manifold points. There exists an
$\epsilon_0>0$ so that $\int_{\supp\phi}|\Ric|^2<\epsilon_0$ implies
\begin{eqnarray*}
\left(\int\phi^4|\nabla{E}|^2\right)^\frac12
&\le&C\int|\nabla\phi|^2|\nabla{E}| \,+\,C\int\phi^2|\Riem||X|\\
&&+\,C\left(\int|\Riem|^2\right)^{\frac12}\left(\int\phi^2|\nabla\phi|
^2|E|^2\right)^\frac12
\end{eqnarray*}
\end{lemma}
\underline{\sl Pf}\\
\indent Set $u=|\nabla{E}|^{1-\delta}$ and use the Sobolev
inequality to get
\begin{eqnarray*}
C\left(\int\phi^4u^{2\frac{1}{1-\delta}}\right)^\frac12
&\le&\int|\nabla\phi|^2u^{\frac{1}{1-\delta}}
\,-\,\int\phi^2u^{\frac{\delta}{1-\delta}}\triangle{u}
\end{eqnarray*}
Since Lemma \ref{ImprovedEllipticUsingKato} holds for
$\alpha=\delta$ and $\int|\Ric|^2$ is assumed small, we get
\begin{eqnarray*}
\left(\int\phi^4|\nabla{E}|^2\right)^\frac12
&\le&C\int|\nabla\phi|^2|\nabla{E}|
\,+\,C\int\phi^2|\triangle\nabla{E}|
\end{eqnarray*}
for $C=C(C_S)$. We'll use
\begin{eqnarray*}
\triangle\nabla{E} &=&\nabla\Riem*E \,+\,\Riem*\nabla{E}
\,+\,\Riem*X
\end{eqnarray*}
First note that $\Riem*X\in{L}^1$, since both are in $L^2$.  Also,
$\int\phi^2|\Riem||\nabla{E}|$ can be combined into the left side.
Altogether,
\begin{eqnarray*}
C\left(\int\phi^4|\nabla{E}|^2\right)^\frac12
&\le&\int|\nabla\phi|^2|\nabla{E}| \,+\,\int\phi^2|\nabla\Riem||E|
\,+\,\int\phi^2|\Riem||X|
\end{eqnarray*}
The Sobolev inequality directly gives
\begin{eqnarray*}
\left(\int\phi^8|E|^4\right)^\frac12
&\le&C\int\phi^2|\nabla\phi|^2|E|^2 \,+\,C\int\phi^4|\nabla{E}|^2.
\end{eqnarray*}
Since we are working in the smooth case, we may also our previous
result that
\begin{eqnarray*}
\left(\int|\nabla\Riem|^\frac43\right)^\frac34
&\le&C\left(\int|\Riem|^2\right)^\frac12
\end{eqnarray*}
assuming the domain of the second integral is somewhat larger than
the domain of the first. The result immediately follows. \qed

As we now show, this lemma gives us just enough to conclude that
$|\nabla\Ric|\in{L}^{\frac43}$. Due to Theorem
\ref{LocalCurvatureBounds} already have $|\nabla\Ric|\in{L}^{p}$ for
all $p<\frac43$. For a similar argument, see the proof of Theorem
5.8 of \cite{BKN}.
\begin{lemma}\label{UhlenbeckPrelimEst}
Its holds that $|\nabla{E}|\in{L}^{\frac43}$, and in fact given any
$\beta>1$,
\begin{eqnarray*}
\int_{B(o,\rho)}|\nabla\Ric|^{\frac43}
&\le&C\int_{B(o,\beta\rho)-B(o,\rho)}|\nabla\Ric|^\frac43
\,+\,C\int_{B(o,\beta\rho)}|\Riem|^2,
\end{eqnarray*}
$C=C(C_S,\beta)$, despite the possible presence of singularities.
\end{lemma}
\underline{\sl Pf}\\
\indent Choosing any $k\in(1,2]$, the previous lemma and H\"older's
inequality gives
\begin{eqnarray*}
\int\phi^{2k}|\nabla{E}|^k
&\le&C(\Vol\supp\phi)^{1-\frac{k}{2}}\left(\int|\nabla\phi|^{\frac{2k}
{k-1}}\right)^{k-1}\int_{\supp\nabla\phi}|\nabla{E}|^k\\
&&+\,C(\Vol\supp\phi)^{1-\frac{k}{2}}\left(\int|\nabla\phi|^2\right)^
{\frac{k}{2}}\left(\max_{\supp|\nabla\phi|}|E|\right)^{k}\left(\int|
\Riem|^2\right)^{\frac{k}{2}}\\
&&+\,C(\Vol\supp\phi)^{\frac1k-\frac12}\int\phi^2|\Riem||X|.
\label{TestFunctionsCloseInInequality}
\end{eqnarray*}
This holds assuming no singularity lies in $\supp\phi$.  Now we let
the $\phi$ be test functions with support everywhere except for
small balls around the singularities. For simplicity we can assume
there is a single singularity; if there multiple singularities this
method requires that the test functions must close in around all of
them simultaneously. Choose some number $\beta>1$ and let $\phi_i$
be a sequence test function with
$\supp\phi_i\cap{B(o,\beta^{-i-1})}=\varnothing$, with
$\phi_i\equiv1$ in $M-B(o,\beta^{-i})$, and with
$|\nabla\phi_i|\le2\beta^{i+1}$.

Set $A_i=\int_{M-B(o,\beta^{-i})}|\nabla{E}|^k$.  Then
$\int_{\supp|\nabla\phi|}|\nabla{E}|^k\;=\;A_{i+1}-A_i$.  With
$|E|=o(r^{-2})$ near singularities, inequality
(\ref{TestFunctionsCloseInInequality}) takes the form
\begin{eqnarray*}
A_i &\le&C\left(A_{i+1}-A_i\right) \,+\,C\beta^{-i(4-3k)}
\,+\,C\beta^{-i(\frac4k-2)}\\
A_i &\le&\frac{C}{1+C}A_{i+1} \,+\,\frac{C}{1+C}\beta^{-i(4-3k)}
\,+\,\frac{C}{1+C}\beta^{-i(\frac4k-2)}.
\end{eqnarray*}
Iterating, we get
\begin{eqnarray*}
A_i &\le&\left(\frac{C}{1+C}\right)^NA_{i+N}\\
&&+\,\frac{C}{1+C}\beta^{-i(4-3k)}\left(1+\frac{C}{1+C}\beta^{-(4-3k)}
\,+\,\dots\,+\,\left(\frac{C}{1+C}\beta^{-(4-3k)}\right)^{N-1}\right)\\
&&+\,\frac{C}{1+C}\beta^{-i(\frac4k-2)}\left(1+\frac{C}{1+C}\beta^{-
(\frac4k-2)}\,+\,\dots\,+\,\left(\frac{C}{1+C}\beta^{-(\frac4k-2)}
\right)^{N-1}\right)
\end{eqnarray*}
An advantage is possible in the boundary case where $k=\frac43$. In
this case clearly the last two terms are bounded independently of
$N$. In the case $k=\frac43$, we know that $A_{i+N}$ grows slower
than any power of $\beta^N$; therefore as $N\rightarrow\infty$ the
first term vanishes.  Thus
\begin{eqnarray*}
A_i &\le&1+C,
\end{eqnarray*}
which is a bound independent of $i$.  Letting $i\rightarrow\infty$
yields the theorem.

Now with $|\nabla{E}|\in{L}^{\frac43}$, one easily gets that
\begin{eqnarray*}
\left(\int\phi^4|\nabla{E}|^2\right)^\frac12
&\le&C\int|\nabla\phi|^2|\nabla{E}| \,+\,C\int\phi^2|\Riem||X|\\
&&+\,C\left(\int|\Riem|^2\right)^{\frac12}\left(\int\phi^2|\nabla\phi|
^2|E|^2\right)^\frac12
\end{eqnarray*}
holds regardless of singularities.  With
$|\nabla\Ric|^2\le|\nabla{E}|^2+\frac1n|X|^2$ and
$|\nabla{E}|\le{C}(n)|\nabla\Ric|$, we get
\begin{eqnarray*}
\left(\int\phi^4|\nabla\Ric|^2\right)^\frac12
&\le&C\int|\nabla\phi|^2|\nabla\Ric| \,+\,C\int\phi^2|\Riem||X|\\
&&+\,Cr^{-1}\int|\Riem|^2 \,+\,C\left(\int\phi^2|X|^2\right)^\frac12
\end{eqnarray*}
Using H\"older's inequality and that $|X|=o(r^{-3})$ points near
singularities,
\begin{eqnarray*}
\left(\int_{B(o,r/2)}|\nabla\Ric|^{\frac43}\right)^\frac34
&\le&C\left(\int_{B(o,r)}|\nabla\Ric|^\frac43\right)^\frac34
\,+\,C\int_{B(o,r)}|\Riem|^2.
\end{eqnarray*}
Using that $\int|\Riem|^2$ is presumed small, we get the lemma. \qed

\begin{lemma}[Uhlenbeck's method]\label{UhlenbeckEstimate}
Assume $|F|=o(r^{-2})$ near singularities and that
$|\nabla\Ric|\in{L}^{\frac43}$. If $o$ is a singuarity, we can
choose $\rho$ small enough and and $\beta$ large enough so that
\begin{eqnarray*}
&&\int_{B(o,\rho)}|F|^2 \;\le\frac12\int_{B(o,\rho\beta)}|F|^2
\,+\,C\int_{B(o,\rho\beta)}|\nabla\Ric|^\frac43.
\end{eqnarray*}
where $C$ is some universal constant.
\end{lemma}
\underline{\sl Pf}\\
\indent We will not present Uhlenbeck's argument in its entirety
here, but our use of it will be unique enough that we must repeat
some of the proof. Uhlenbeck first proves that a gauge can be found
on the annulus so that ${D}^{*}A=0$, where $A$ is the connection
1-form.  One of the main advantages of computing in this special
gauge is that integral norms of $|A|$ are bounded in terms of those
of $|F|$ (see \cite{Tia1}, \cite{Uhl}).  In fact given a domain
$\Omega$, we can get
\begin{eqnarray}
\int_\Omega|A|^2 \;\le\;C\int_\Omega|F|^2
\label{UhlenbeckAInTermsOfF}
\end{eqnarray}
where $C=C(\Omega)$.  Our computation is similar to those in
\cite{Uhl}, \cite{Tia1}, and \cite{TV1}, but we use test function
methods rather than try to control the boundary terms. We will use
the second Bianchi identity $ D^{*}F=D\Ric$.  We get
\begin{eqnarray}
&&\int\phi^2|F|^2
\;=\;\int\phi^2\left<DA-\frac12[A,A],F\right>\label {UhlenbeckMainEstimate}\\
&&\quad\quad=\;\int\phi^2\left<DA,F\right>
\,-\,\frac12\int\phi^2\left<[A,A],F\right>\nonumber\\
&&\quad\quad=\;-2\int\phi\left<\nabla\phi\otimes{A},F\right>
\,-\,\int\phi^2\left<a,D^{*}F\right>
\,-\,\frac12\int\phi^2\left<[A,A],F\right>\nonumber\\
&&\quad\quad=\;-2\int\phi\left<\nabla\phi\otimes{A},F\right>
\,-\,\int\phi^2\left<a,D\Ric\right>
\,-\,\frac12\int\phi^2\left<[A,A],F\right>\nonumber
\end{eqnarray}
In a Hodge gauge it is possible to estimate $\int|A|^4$ in terms of
$\int|F|^2$. The Sobolev inequality gives
\begin{eqnarray}
\left(\int|A|^4\right)^\frac12 &\le&C_S\int|DA|^2
\,+\,(\Vol\supp\phi)^{-\frac12}\int|A|^2\nonumber\\
&&\le\;C_S\int|F|^2
\,+\,\left(C_S\sup_{\supp\phi}|A|^2+(\Vol\supp\phi)^{-\frac12}\right)
\int|A|^2\nonumber\\
&&\le\;C\int|F|^2
\end{eqnarray}
where we have used (\ref{UhlenbeckAInTermsOfF}). Here $C$ depends on
the Sobolev constant and on $\sup_{\supp\phi}|F|$. We now get from
(\ref{UhlenbeckMainEstimate}),
\begin{eqnarray}
\int\phi^2|F|^2 &\le&c\int|\nabla\phi|^2|A|^2
\,+\,c\int\phi^2|D\Ric|^\frac43 \label{FirstFEstimate}
\end{eqnarray}
We want to estimate $\int|A|^2$ back in terms of $\int|F|^2$, but we
need to control the coefficient. We get becomes
\begin{eqnarray*}
\int\phi^2|F|^2
&\le&c\left(\int|\nabla\phi|^4\right)^\frac12\left(\int_{\supp\phi}|A|
^4\right)^\frac12
\,+\,c\int\phi^2|D\Ric|^\frac43\\
&\le&c\left(\int|\nabla\phi|^4\right)^\frac12\int|F|^2\,+\,c\int
\phi^2|D\Ric|^\frac43
\end{eqnarray*}
Assuming $\phi$ is defined in the annulus $B(o,1)-B(o,\beta^{-1})$,
we can make $\int|\nabla\phi|^4$ very small by making $\beta$ large;
in fact we can make $\int|\nabla\phi|^4\backsim(\log\beta)^{-3}$.
This done, we get
\begin{eqnarray*}
\int\phi^2|F|^2 &\le& \epsilon\int_{\supp\phi}|F|^2
\,+\,c\int\phi^2|D\Ric|^\frac43
\end{eqnarray*}
We can choose $\beta$ large enough so that $\epsilon<\frac14$. A
significant subtlety is that as the annulus goes to zero
$B(o,1)-B(o,\beta^{-1})$ degenerates to a punctured disk, the
estimate (\ref{FirstFEstimate}) does not degenerate. It is possible
to prove this with a modification of the argument on pg. 129 of
\cite{Tia1}; see \cite{Web1} for the details.

We now piecing together successively smaller annuli, in order to
close in around the singularity. Let $\phi_i$ be a test function
with $\phi_i\equiv1$ in $B(o,\beta^{-i-1})-B(o,\beta^{-i-2})$,
$\phi_i\equiv0$ in $B(o,\beta^{-i-3})$ and outside
$B(o,\beta^{-i})$, and also
$\int|\nabla\phi_i|^4\le{C}(\log\beta)^{-3}$.  Then our inequality
reads
\begin{eqnarray*}
&&\int_{B(o,\beta^{-i-1})\cap{B}(o,\beta^{-i-2})}|F|^2\\
&&\quad\le\epsilon\int_{B(o,\beta^{-i-2})-B(o,\beta^{-i-3})}|F|^2\\
&&\quad\quad+\,\epsilon\int_{B(o,\beta^{-i})-B(o,\beta^{-i-1})}|F|^2\\
&&\quad\quad+\,c\int_{B(o,\beta^{-i})-B(o,\beta^{-i-3})}|D\Ric|^\frac43
\end{eqnarray*}
Now summing both sides from $i=N$ to $\infty$ gives
\begin{eqnarray*}
&&\int_{B(o,\beta^{-N-1})}|F|^2
\;\le\frac12\int_{B(o,\beta^{-N})}|F|^2
\,+\,3c\int_{B(o,\beta^{-N})}|D\Ric|^\frac43.
\end{eqnarray*}
\qed

\begin{proposition} \label{Dim4RiemImprovement}
Assuming Lemmas \ref{UhlenbeckPrelimEst} and
\ref{UhlenbeckEstimate}, we have that $|F|\in{L}^p$ for some $p>2$.
\end{proposition}
\underline{\sl Pf}\\
\indent Propositions \ref{UhlenbeckPrelimEst} and
\ref{UhlenbeckEstimate} now give
\begin{eqnarray*}
&&\int_{B(o,\rho)}|F|^2 \;\le\;\int_{B(o,\rho\beta)-B(o,\rho)}|F|^2
\,+\,3c\int_{B(o,\rho\beta)}|D\Ric|^\frac43.
\end{eqnarray*}
\begin{eqnarray*}
\int_{B(o,\rho)}|\nabla\Ric|^\frac43
&\le&C\int_{B(o,\rho\beta)-B(o,\rho)}|\nabla\Ric|^{\frac43}
\,+\,C\int_{B(o,\rho\beta)}|F|^2.
\end{eqnarray*}
Setting
\begin{eqnarray*}
A_i\;=\;\int_{B(o,\rho\beta^{-i})}|F|^2\quad\quad\quad\quad
B_i\;=\;\int_{B(o,\rho\beta^{-i})}|\nabla\Ric|^{\frac43}.
\end{eqnarray*}
we have
\begin{eqnarray*}
A_i&=&C(A_{i-1}-A_i) \,+\,CB_{i-1}\\
B_i&=&C(B_{i-1}-B_i) \,+\,CA_{i-1}
\end{eqnarray*}
It is possible to iterate these to get
\begin{eqnarray*}
A_i&=&\left(\frac{C}{1+C}\right)^i\left(A_{0}\,+\,B_{0}\right).
\end{eqnarray*}
Thus choosing $s>1$ so that $\beta^{-s}=\frac{C}{1+C}$ we get
\begin{eqnarray*}
\int_{B(o,\rho\beta^{-i})}|F|^2 &\le&C'\beta^{-si}.
\end{eqnarray*}
This proves the existence of an $s>0$ so that
$\int_{B(o,r)}|F|^2=O(r^s)$.  Using elliptic regularity (Theorem
\ref{RiemBounds}) we get that $|\Riem|=O(r^{-2+s})$ near
singularities.  Therefore $|\Riem|\in{L}^p$ for any
$p<\frac{4}{2-s}$. \qed

\begin{theorem} \label{LocalCurvatureBoundsDim4}
Assume $g$ is an extremal K\"ahler metric on a Riemannian
manifold-with-singularities of dimension 4.  When $a>2$, there
exists $\epsilon_0 = \epsilon_0(C_S,a)$ and $C=C(C_S,a)$ so that
$$
\int_{B(o,r)}|\Riem|^{2} \;\le\; \epsilon_0
$$
implies
\begin{eqnarray}
&&\left(\int_{B(o,r/2)}|\Riem|^a\right)^\frac1a
\;\le\;Cr^{-2+\frac{4}{a}}\left(\int_{B(o,r)}|\Riem|^2\right)^\frac12.
\end{eqnarray}
\end{theorem}
\underline{\sl Pf}\\
\indent Now that we know $|\Riem|\in{L}^p$ for some $p>2$, we can
use Sibner's Lemma and Proposition \ref{LpLemma}, and repeat the
proof of Proposition \ref{PropRiemNablaRicLpBounds}. Actually, in
the case $n=4$ Proposition \ref{LpLemma} is not quite strong enough
as stated. Referring to the notation from the statement of
\ref{LpLemma}, it is required that $q>2$. This can be changed to
allow equality however. First one makes the following change to the
statement of Lemma \ref{L2Lemma}: Given $k>\frac12\frac{n}{n-2}$, if
$u^k\in{L}^2$ then given any $l\le{k}$ it holds that
$\nabla{u^l}\in{L}^2_{loc}$. \qed

We would also like to point out that our improved Kato inequality is
sufficient for proving an improved curvature decay rate at infinity.
That is
\begin{eqnarray*}
|\Riem| &=&O(r^{-2-s})
\end{eqnarray*}
as $r\rightarrow\infty$, for some $s>0$. This can be proven using
Uhlenbeck's method, essentially just by taking annular regions
extending out to infinity rather than in toward a singularity. But
considering our rather bulky use of the improved Kato inequality, it
is unlikely our method will allow computation of the optimal decay
rate.

\section{Weak Compactness} \label{WeakCompactness}

In this section we assume our manifolds satisfy a local volume
growth upper bound, which is a significant assumption without
global, pointwise lower bounds on the Ricci curvature.  Following
\cite{TV2}, we can use the convergence result proved here to turn
around and actually {\sl prove} the volume growth assumption, which
is done by scaling the manifolds so the local growth condition does
hold, and then applying the results of this section. Essentially the
possibility of large local volume ratios is counterbalanced by the
freedom, in the following argument, to let diameters be as large as
desired or even infinite.

We shall adopt the following definition of asymptotically locally
Euclidean manifolds: a complete manifold will be called ALE if there
exists a compact set $K\subset M$ so that each component of $M-K$ is
diffeomorphic to $(\mathbb{R}^n-B)/\Gamma$ for some ball $B\in
\mathbb{R}^n$ and some subgroup $\Gamma\subset SO(n)$ (depending on
the end), and so that under this identification, the metric
components satisfy
\begin{eqnarray*}
&&g_{ij}\;=\;\delta_{ij}\,+\,o(1)\\
&&\partial^{k}(g_{ij})\;=\;o(r^{-k}),
\end{eqnarray*}
where $\partial^k$ indicates any partial derivative of order $k$. In
\cite{TV1} for instance, such a manifold is called ALE of order $0$.

In this section, we assume
$\{(M_\alpha,g_\alpha,x_\alpha)\}_{\alpha\in A}$ is a family of
compact, pointed $n$-dimensional Riemannian manifolds that satisfy
\begin{itemize}
\item[{\it i})] Upper bounds on energy: $\int_{M_\alpha}|\Riem|^{\frac {n}{2}}
\le \Lambda$
\item[{\it ii})] Lower bounds on volume: $\Vol_{g_\alpha} M_\alpha  \ge \nu$
\item[{\it iii})] Weak regularity: $\int_{B_r}|\Riem|^{\frac{n}{2}} <
\epsilon_0 \;\Rightarrow\;\sup_{B_{r/2}}|\nabla^p \Riem| \le Cr^{-
p-2} \left(\int_{B_r}|\Riem|^{\frac{n}{2}}\right)^\frac2n$
\item[{\it iv})] Bounded Sobolev constants $C_M<C_S$.
\item[{\it v})] Upper bound on local volume growth: $\Vol_{g_\alpha} B (p,r)
\le\overline{v}r^n$ for $0 \le r \le 1$
\end{itemize}

\begin{proposition}\label{OpenLimits}
Let $\{(M_\alpha,g_\alpha,x_\alpha)\}_{\alpha\in A}$ be a family of
pointed, compact Riemannian manifolds that satisfy the above
conditions. Then a subsequence $\{(M_i,g_i,x_i)\}_{i=i}^\infty$
converges in the pointed Gromov-Hausdorff topology to a complete
pointed Riemannian manifold-with-singularities $(M_\infty, g_\infty,
x_\infty)$ with at most $\Lambda/\epsilon_0$ singularities.  If
$M_\infty$ is noncompact, it is $ALE$.
\end{proposition}
\underline{\sl Pf}

Similar arguments appear frequently in the literature, so we briefly
describe the main steps. Choose a small radius $r>0$. Let $K \subset
M$ be the (compact) set of points $p\in M$ where
$\int_{B(p,r)}|\Riem|^{\frac{n}{2}} \ge \epsilon_0$. Cover $K$ by
balls $B(p_i,2r)$ such that the $B(p_i,r)$ are disjoint; there can
be no more than $\Lambda/\epsilon_0$ balls in such a covering. Set
$\Omega_{i,r,R} = \left(M_i - \bigcup_j B(p_j,2r)\right) \,\cap\,
B(x_i,R)$. Notice that when $r$ is small enough, the local volume
growth bounds give $\Vol\Omega_{i,r,R} \ge \Vol_{g_i} B(x_i,R) -
\epsilon$.

On $\Omega_{i,r,R}$ we have $|\nabla^k \Riem| \le C r^{-k-2}$. The
lower bound on volume growth together with the curvature estimate
imply the Cheeger lemma (\cite{Che}) which gives injectivity radius
bounds. Therefore we can take a pointed limit along a subsequence of
the sets $\Omega_{i,r,R}$ to get a smooth limiting
manifold-with-boundary $\Omega_{\infty,r,R}$.

This convergence is smooth in the topology, and $C^{k+1,\alpha}$ in
the metric by our $L^{\infty}$ bounds on the $k^{th}$ derivative of
curvature. We get diffeomorphisms
$\Phi_{i,r,R}:\Omega_{\infty,r,R}\rightarrow\Omega_{i,r,R}$ for
large $i$ such that the pullback metrics $\Phi_{i,r,R}^{*}g_i$
converge smoothly to $g_\infty$. Adjusting $r$ will change the limit
manifold, but the limit manifolds naturally embed in one another.
Put
\begin{eqnarray*}
\Omega_{\infty,R} &=& \bigcup_{0<r}\Omega_{\infty,r,R} \\
\Omega_\infty &=& \bigcup_{R < \infty}\Omega_{\infty,R}.
\end{eqnarray*}
The local upper bound on volume growth insures that $\Omega_\infty$
can be completed by adding discrete points, which constitute the
singular set $S$, which has cardinality at most
$\Lambda/\epsilon_0$. The result is a complete
manifold-with-singularities $M_\infty = \Omega_\infty \cup S$. It is
possible that $S$ is empty, or that some points of $S$ might be
smooth points of $M$.

In theorem 4.1 of \cite{TV1}, Tian-Viaclovsky show essentially that
a complete manifold-with-singularities with Euclidean volume growth
and $|\Riem|\cong o(r^{-2})$ at infinity is in fact ALE. Their
method of proof is geometric and will hold in any dimension for
manifolds-with-singularities, though it is stated for
$4$-dimensional smooth manifolds (see theorem 4.1 of \cite{TV1} and
the comment immediately afterwards). In our setting, the volume
growth lower bound is implied by the Sobolev constant bound and
quadratic curvature decay is ensured by condition ({\it iii}). An
assumption on $b_1(M)$ is not necessary due to the improvements in
\cite{TV3}. Thus our limit manifold, if noncompact, will be ALE.
\hfill$\Box$

Next we examine the curvature singularities that arise in the limit,
and, following the existing literature, sketch a proof that they are
indeed $C^\infty$ (possibly nonreduced) orbifold points.  As is
common in the literature, we say an orbifold possesses some
structure if the structure exists at smooth points and, after
lifting to the smooth orbifold cover of any point, it can be
completed. For instance, an orbifold is called extremal K\"ahler if
the lift of the metric to the orbifold cover of any point can be
completed.

We also consider the {\sl order} of the multifold points that arise
in the limit.  We define order as follows: if $o$ is a multifold
point with tangent cone $T$ at $o$, the order of $o$ is just the
cardinality of the set of components of $T-\{o\}$. We will often use
the terms orbifold and multifold interchangably. When it is
important to distinguish, we shall call an orbifold {\sl reduced}
when each singular point has order 1.

\begin{proposition} \label{CurvatureSingularityStructure}
Assume $M$ is a Riemannian manifold-with-singularities and that $M$
carries an extremal K\"ahler metric at every smooth point. Then the
singularities are $C^\infty$ Riemannian multifold points, and $M$ is
an extremal K\"ahler multifold. Further, the cardinality of any
orbifold group $\Gamma$ has a bound $|\Gamma| \le C(C_S)$, and the
order of the multifold points are bounded by
$C=C(C_S,\overline{v})$.
\end{proposition}
\underline{\sl Pf}\\
\indent This is a local proof; we need only consider neighborhoods
of singularities. Most of the work here is identical to that found
elsewhere in the literature. Let $o$ be a singularity. First choose
a locally connected component $N$ of $M-\{o\}$; by this it is meant
that $(N\cap{B}(o,r))-\{o\}$ is connected for all $r>0$.  We know
that $|\Riem|=o(r^{-2})$ on $\overline{N}$, where $r$ is the
distance to $o$, and so the proof Lemma (5.13) of \cite{BKN} yields
that $\overline{N}$ has a unique tangent cone at $o$ that is
diffeomorphic to $\mathbb{R}^{n}/\Gamma$ where $\Gamma$ is some
isometric action on $\mathbb{R}$ whose only fixed point is $\{o\}$.
Since the Sobolev inequality holds on $\overline{N}$ and hence
$\overline{N}$ has local volume growth lower bounds, there is a
bound on the cardinality $|\Gamma|$ of the orbifold group that
depends only on $n$ and $C_S$.  Lastly, at any singularity point
$\{o\}$, any small ball $B(o,r)$ must have at most a uniformly
bounded number of components.  This is because each component has
local volume growth lower bounds, so many components together would
give a very large local volume growth; this is impossible by
assumption.

Now we examine the regularity of the metric on the orbifold cover of
any component of a multifold point. Let $B=B(o,\epsilon)$ be a small
ball around $o$ diffeomorphic to $T$. Choose one component of
$B-\{o\}$ and consider its orbifold cover (a neighborhood of the
origin in $\mathbb{R}^n$). Lifting the metric to this neighborhood,
we must analyze the regularity of the metric at a deleted point of
$\mathbb{R}^n$.

With bounded curvature and dimension $n>2$, elementary arguments
show the metric is $C^0$. A less elementary argument suffices to
construct $C^{1,1}$ coordinates; for instance the construction of
\cite{BKN} beginning on pg 342 shows this to be possible.  We are
able to cite this result in the higher dimensional case by Theorem
\ref{LocalCurvatureBounds}, and in dimension 4 by Theorem
\ref{Dim4RiemImprovement}. With $C^{1,1}$ coordinates, it is
possible to construct harmonic coordinates, as in \cite{DK}.

In harmonic coordinates, we have the coupled system
\begin{eqnarray*}
&&\triangle(g_{ij}) \;=\; {\Ric}_{ij} \,+\, Q(g,\partial g)\\
&&\triangle\Ric \;=\; \Riem*\Ric \,+\,\nabla X\\
&&\triangle X \;=\; \Ric*X.
\end{eqnarray*}
A bootstrapping argument is possible using the $L^p$ theory.  Since
$\Ric\in{L}^p$ for some $p>2$ and $g$ and $\partial g$ are bounded,
the first equation gives $g_{ij}\in W^{2,p}$ for all $p$, so
$\Ric\in W^{0,p}$.  Then the third equation gives $(\nabla
X)_{ij}\in W^{1,p}$ and so the second equation gives $\Ric_{ij}\in
W^{2,p}$. Then the first equation gives $g_{ij}\in W^{4,p}$ and
therefore by Sobolev imbedding $g_{ij}\in C^{3,\alpha}$, so
$\Riem\in C^{1,\alpha}$. Now we turn to the Schauder theory.  In
harmonic coordinates the coefficients of $\triangle$ are
$C^{3,\alpha}$, so the last equation gives $X_i\in C^{3,\alpha}$, so
$(\nabla X)_{ij} \in C^{2,\alpha}$, so the middle equation gives
$\Ric_{ij}\in C^{3,\alpha}$. Then with $\partial g\in C^{2,\alpha}$,
the first equation gives $g_{ij}\in C^{4,\alpha}$, an improvement in
regularity. Bootstrapping like this gives $g_{ij}\in C^{k,\alpha}$
for all $k$, so $g\in C^\infty$.  All of this is standard elliptic
theory; see for instance chapter 5 of \cite{Evn} for Sobolev
embedding, and chapters 6 and 9 of \cite{GT} for the Schauder theory
and $L^p$ theory.  This completes the proof that our curvature
singularities are $C^{\infty}$ Riemannian multifold points.

Finally we check the complex structure on the orbifold covers. Since
the tensor $J$ is harmonic (indeed, covariant constant) and of
bounded norm, its lift will extend smoothly over the deleted point.
The completed complex structure is clearly integrable, since the
Nijenhuis tensor is smooth and is assumed to vanish everywhere
except at the origin, and so it vanishes everywhere. Also, since
$R_{,\bar{i}\bar{j}}=0$ outside the singularity and $R$ is
$C^\infty$, after completion $R_{,\bar{i}\bar{j}}=0$ everywhere, so
the metric on the orbifold cover is extremal K\"ahler. \qed

\begin{proposition}[Limits are reduced orbifolds] \label
{LimitsAreReducesOrbifolds} Suppose
$(M_\alpha,g_\alpha,x_\alpha)_{\alpha\in A}$ is a family of
$n$-dimensional extremal K\"ahler manifolds that satisfy the
conditions (i)-(v) of this section, and which also have a local
volume growth upper bound. Then a subsequence converges to a reduced
extremal K\"ahler orbifold. If $\Gamma$ is an orbifold group then
$\Gamma\subset{U}(n)$, and there is a bound on its order,
$|\Gamma|\le C(C_S,n)$. There is a bound the number of orbifold
points, given by $C=C(n,\Lambda,C_S)$.
\end{proposition}
\underline{\sl Pf}\\
\indent Proposition (\ref{CurvatureSingularityStructure}) shows that
any manifold-with-singularities constructed in the proof of
proposition (\ref{OpenLimits}) is indeed a Riemannian multifold.  We
need only pass to a further subsequence to get a converging almost
complex structure.  The limiting complex structure is clearly
integrable, since the Nijenhuis tensor will continue to vanish at
all smooth points of the limit. $C^1$ convergence at smooth points
implies also $d\omega = 0$ (where $\omega$ is the K\"ahler form), so
the limiting multifold is K\"ahler at smooth points, and $C^4$
convergence guarantees that $R_{,\bar{i}\bar{j}}=0$, so the
multifold metric is extremal at smooth points.  The orbifold group
is a subgroup of $U(n)$ because its action on the cover preserves
$J$.

Finally we establish that the limit is actually a {\sl reduced}
orbifold, meaning that $B(o,r)-\{o\}$ has only one component
regardless of $o$ or $r$. To this end we do a blowup analysis at a
forming singularity in order to capture a two (or more) ended
singularity model.  Assume $\{o\}$ is a singularity that locally
separates $M_\infty$. Let $p_i\in{M}_i$ be a sequence of points with
$p_i\rightarrow\{o\}$, and let $B(p_i,r_i)$ be balls with the
following property: $\partial{B}(p_i,r_i)$ has one component, but
whenever $r_i<\rho<i\,r_i$ then $\partial{B}(p_i,\rho)$ has more
than one component (one must generally pass to a subsequence here).

Now rescale the manifolds $M_i$ by setting $\bar{g}_i=r_i^{-2}g_i$,
and take a limit.  By Proposition \ref{OpenLimits} we know the limit
is an ALE manifold-with-singularities, which we know are $C^\infty$
orbifold points by Proposition \ref{CurvatureSingularityStructure}.
Since $B(o,0.99)\subset{M}_\infty$ does not separate $M_\infty$ but
any ball $B(o,r)$ of radius $r>1$ has more than one boundary
component, we know the limit has more than one end.  If the limit
has a locally separating singularity, repeat the process until we
arrive at a limiting object whose singularities do not locally
separate.

We can use Theorem 4.1 of \cite{LT1} to conclude that $M$ has at
most one non-parabolic end, and therefore at least one end is
parabolic. However the result of Holopainen-Koskela or of Li-Tam
(Theorem 1.4 of \cite{HK}, Theorem 1.9 of \cite{LT2}) imply that
none of our ends are parabolic. This contradiction establishes the
proof. \qed

There has been a great deal of work relating the function-theoretic
aspects of manifolds to their Riemannian or K\"ahlerian geometry.
We'd like to mention the nice survey article \cite{Li} by Peter Li.

Finally we are able to complete the proof of Theorem
\ref{FinalCompactnessTheorem0} or \ref{FinalCompactnessTheorem} by
removing condition ({\it v}) from the list at the beginning of this
section.  Following the proof of \cite{TV2}, we prove Theorem 5.5
for the case of extremal K\"ahler metrics.  First we cite a volume
comparison lemma; see for instance Proposition 20 of \cite{Bor}.

\begin{lemma}[Orbifold volume comparison] \label{MultiVolComp}
Assume $M^n$ is a smooth Riemannian orbifold.  Let $B=B(p,r)\subset
M$ be any ball. If $\Ric \ge -(n-1)H$ in $B - S$, then $\Vol B(p,r)
\le {\Vol}_{-H} B(r)$.
\end{lemma}

Define the maximal volume ratio $MV_t(M)$ of $M^n$ at scale $t$ to
be
$$
MV_t(M) = \sup_{x\in M,\,0<r<t} r^{-n} \Vol B(x,r).
$$
$MV_\infty(M)$ is of course an upper bound on the volume ratio of
balls in $M$.  We will also denote by $\Vol_c B(t)$ the volume of
the ball of radius $t$ in the space form of constant sectional
curvature $c$.

\begin{theorem}[Upper bound on volume growth] \label{VolGrowthTheorem}
Let $(M_{\lambda},g_{\lambda})_{\lambda\in A}$ be a family of
compact, extremal K\"ahler manifolds. Assume
$\Vol_{g_\lambda}(M_\lambda)\ge\nu$,
$\Diam_{g_\lambda}(M_\lambda)\le\delta$,
$\|\Riem\|_{L^{\frac{n}{2}}} \le \Lambda$, and Sobolev constants
$C_{M_\lambda}$ bounded above by $C_S < \infty$.  Then there exists
a bound on $MV_\infty(M_\lambda)$ depending on $C_b$, $\nu$,
$\delta$, $\Lambda$, and $C_S$.
\end{theorem}
\underline{\sl Pf} \\
\noindent Assume no such bound exists, so there is a sequence of
such Riemannian manifolds $M_i = \{M, g_i\}$ with $MV_\infty
(M_i)\rightarrow\infty$.

First, $\int_{B(x,2r)} |\Riem|^2 \le \epsilon_0$ implies $|\Riem|
\le C\epsilon_0r^{-2}$ in $B(x,r)$, so assuming (without loss of
generality) that $r\le\delta$, Bishop volume comparison gives $\Vol
B(x,r)\le r^n\delta^{-n}\Vol_{-C\epsilon_0}B(\delta)\triangleq
Ar^n$.

Choose points $x_i\in M_i$ and radii $r_i$ so that $\Vol B(x_i,r_i)
= 2A (r_i)^4$, and $r_i$ has the following minimality property:
whenever $p\in M_i$ and $r\le r_i$, we have $\Vol B(p,r) \le 2A
r^4$. In other words, $MV_{r_i}(M_i) = 2A$. Note also that
$\int_{B(x_i,2r_i)}|\Riem|^2 \ge \epsilon_0$.

Set $x_i^{(1)} = x_i$, $r_i^{(1)} = r_i$, and $A^{(1)} = 2A =
2\Vol_{-C\epsilon_0} B(1)$. For an induction argument assume that
$k$ sequences of balls $\{B_i^{(1)}\}_{i=1}^{\infty},\,\dots,\,
\{B_i^{(k)}\}_{i=1}^{\infty}$ have been chosen, where $B_{i}^{(j)}
\triangleq B(x_i^{(j)},r_i^{(j)})$ and $B_i^{(j)} \subset M_i$, and
assume the balls satisfy the following assumptions:
\begin{itemize}
\item the balls in the $j^{th}$ sequence, $B_i^{(j)}$, have volume  ratio
$(r_i^{(j)})^{-n} \cdot \Vol B_i^{(j)} \triangleq A^{(j)}$  fixed
independent of $i$
\item $A^{(j+1)} \ge 2 A^{(j)}$
\item each ball $B_i^{(j)}$ has the largest volume ratio among all  balls in
$M_i$ of equal or smaller radius.
\item for large $i$, $\int_{\bigcup_{j=1}^{k}B_i^{(j)}} |\Riem|^{\frac {n}{2}}
\ge k\epsilon_0,$
\end{itemize}
For an induction argument we will show it is possible to extract a
$(k+1)^{th}$ sequence with the same assumptions.

Choose one of the sequences $\{B_i^{(l)}\}_{i=1}^\infty$; we garner
geometric information around the points $x_i^{(l)}$ by blowing up
with $x_i^{(l)}$ as the basepoint. Scale each manifold $M_i$ so that
the $i^{th}$ ball has radius $1$, by setting $\tilde{g}_i =
(r_i^{(l)})^{-2}g_i$. With the new metrics, we have an upper bound
on the volume ratio for balls of radius $\le 1$, so after passing to
a subsequence we get convergence to a limit multifold $(M_\infty,
g_\infty)$. We know the limit multifold is ALE and therefore has a
global upper bound on volume growth, meaning $MV_\infty(M_\infty)
\triangleq L^{(l)} < \infty$. Obviously $L^{(l)} \ge A^{(l)}$, so
also $L^{(l)} \ge 2^l A$. We will denote the scaled radii
$\tilde{r}^{(j)}_i = r^{(j)}_i / r^{(l)}_i$ and the scaled balls
$\tilde{B}^{(j)}_i = \tilde{B}(x^{(j)}_i,\tilde{r}^{(j)}_i)$. Of
course $\tilde{B}^{(l)} = \tilde{B}(x^{(l)}_i,1)$.

Return now to the unscaled manifolds. Choose a $(k+1)^{th}$ sequence
of balls, $B_i^{(k+1)} = B(x^{(k+1)}_i,r^{(k+1)}_i)$, with volume
ratio $A^{(k+1)} = 2 \cdot 6^4 \cdot L^{(k)} (\ge 2A^{(k)})$, and so
that $B^{(k+1)}_i$ has the largest volume ratio among all balls of
equal or smaller radius.  We will prove that for large $i$,
$\int_{\bigcup_{j=1}^{k+1} B(x_i^{(j)},2r^{(j)}_i)} |\Riem|^2 \ge
(k+1)\epsilon_0$. This completes the induction argument and yields a
contradiction with the $L^{\frac{n}{2}}$ curvature bound on the
$M_i$.

We know that $r_i^{(k+1)} > r_i^{(l)}$, $l \le k$.  But we don't
know whether $r_i^{(k+1)}/r_i^{(l)}$ is bounded.

\underline{\sl Case I: If $r_i^{(k+1)}/r_i^{(l)}$ is bounded,
$B(x_i^{(k+1)},2r_i^{(k+1)})$ is eventually disjoint} \newline
\underline{\sl from $B(x_i^{(l)},2r_i^{(l)})$}

Assuming $r_i^{(k+1)}/r_i^{(l)} = \tilde{r}_i^{(k+1)} \le N$,we can
prove that for large $i$, $\tilde{B}(x_i^{(l)},2)$ and
$\tilde{B}(x_i^{(k+1)},2\tilde{r}_i^{(k+1)})$ are disjoint. When $i$
is large, regions in $(M,\tilde{g}_i)$ around $x_i^{(l)}$ are very
close to the limiting orbifold, so we have volume ratio bounds for
balls near $x_i^{(l)}$: in particular, with $i$ large we have
$\Vol_{\tilde{g}_i} \tilde{B}(x_i^{(l)},\tilde{r}) < 2L^{(l)}
\tilde{r}^4$ for arbitrary $\tilde{r}\le6N$.

Balls of bounded radius a bounded distance away from $x_i^{(l)}$
must have curvature ratio nearly bounded by $L^{(l)}$, the global
volume ratio of the limit orbifold. Specifically, if
$\tilde{B}(x^{(l)}_i,2)\cap
\tilde{B}(x_i^{(k+1)},2\tilde{r}_i^{(k+1)}) \ne \varnothing$, then
$\tilde{B}(x_i^{(k+1)},2\tilde{r}_i^{(k+1)}) \subset
\tilde{B}(x_i^{(l)},6\tilde{r}_i^{(k+1)})$ so $\Vol
\tilde{B}(x_i^{(k+1)},\tilde{r}_i^{(k+1)}) \le \Vol
\tilde{B}(x_i^{(l)},6\tilde{r}_i^{(k+1)}) < 2 L^{(l)}
6^4(\tilde{r}_i^{(k+1)})^4$, a contradiction. Unscaling now, we have
that $B(x_i^{(k+1)}, 2r_i^{(k+1)})$ is disjoint from $B(x_i^{(l)},
2r_i^{(l)})$.

If this argument works for all $l\le k$, we have
$\int_{\bigcup_{l=1}^{k+1} B(x_i^{(l)}, r_i^{(l)})} |\Riem|^2 \ge
(k+1)\epsilon_0$ as desired. If not, for any $l$ with $r^{(k+1)}_i /
r^{(l)}_i$ unbounded, we move to the second case.

\underline{\sl Case II: If $r_i^{(k+1)}/r_i^{(l)}$ is not bounded,
$B(x_i^{(k+1)},r_i^{(k+1)})$ eventually has a} \newline
\underline{\sl region of high curvature that is disjoint from any of
the other balls}

Passing to a subsequence, we can assume
$\tilde{r}_i^{(k+1)}\rightarrow\infty$. The idea is that if
$\tilde{r}_i^{(k+1)}$ becomes unboundedly large, we might not have
disjointness of $B(x_i^{(k+1)},2\tilde{r}_i^{(k+1)})$ from the
smaller balls, but the smaller balls are nearly multifold points and
therefore can only multiply the overall volume ratio by a controlled
amount. Since the volume ratio is very large in $B^{(k+1)}_i$,
volume comparison forces some region inside $B^{(k+1)}_i$ to have
large $|\Riem|$ and therefore large $L^{\frac{n}{2}}$-norm of
curvature disjoint from the other balls.

Assume $\int_{B(x^{(k+1)}_i,2r^{(k+1)}_i) - \bigcup_{j=1}^l
B(x_i^{(j)},2r_i^{(j)})} |\Riem|^2 \le \epsilon_0$ for all large
$i$, for if not we are done. We know the unscaled radii
$r^{(k+1)}_i$ are bounded as $i\rightarrow\infty$, so we can do a
blowup limit by scaling $\tilde{g}_i = \left(r^{(k+1)}_i\right)^{-2}
g_i$ and reach a limit orbifold
$(\tilde{M}_\infty,\tilde{g}_\infty)$ with basepoint
$x^{(k+1)}_\infty$. We can assume $\tilde{B}(x^{(k+1)}_\infty,2)$
has $k$ or fewer orbifold points, which correspond to the limits of
the centers of the balls $\tilde{B}_i^{(j)}$. For if there are than
$k+1$ multifold points, there are at least $k+1$ curvature
concentration points in $B(x^{(k+1)}_i,2r^{(k+1)}_i)$, and we are
done.

Letting $S$ be the singular set in $\tilde{M}_\infty$, Fatou's lemma
ensures we have $\int_{\tilde{B}(x^{(k+1)}_\infty,2)-S}|\Riem|^2 \le
\epsilon_0$. The Moser iteration argument works despite the presence
of orbifold points, and we have that $|\Riem| \le C\epsilon_0$ on
$\tilde{B}$, so orbifold volume comparison now guarantees that
$$
\Vol \tilde{B} \;\le\; {\Vol}_{-C\epsilon_0} B(1),
$$
which violates the fact that we chose the $B^{(k+1)}_i$ to have
volume ratio $\ge 2^{k+1}A = 2^{k+1}\Vol_{-C\epsilon_0} B(1)$. In
the unscaled manifold we must have
\begin{eqnarray*}
\int_{B(x^{(k+1)}_i,2r^{(k+1)}_i) - \bigcup_{j=1}^{k}
B(x^{(j)}_i,2r^{(j)}_i)} |\Riem|^2 &\ge&\epsilon_0,
\end{eqnarray*}
so therefore $\int_{\bigcup_{j=1}^{k+1} B(x^{(j)}_i,2r^{(j)}_i)}
|\Riem|^2 \ge (k+1)\epsilon_0$ \qed

\begin{theorem}[Orbifold compactness] \label{FinalCompactnessTheorem}
Any family $\{M_\alpha,J_\alpha,g_\alpha\}_{\alpha\in A}$ of
extremal manifolds satisfying conditions ({\it i}) - ({\it iv}) of
the introduction contains a subsequence $\{M_i,J_i,g_i\}$ that
converges in the Gromov-Hausdorff topology to a reduced compact
extremal K\"ahler orbifold. Further, there is a bound
$C_1=C_1(\Lambda,C_S,n)$ on the number of singularities, and a bound
$C_2=C_2(C_S,n)$ on the order of any orbifold group.
\end{theorem}
\underline{\sl Pf} \\
\indent In light of Theorem \ref{VolGrowthTheorem}, Proposition
\ref{LimitsAreReducesOrbifolds} now goes through without a separate
assumption on local volume growth. \qed

In light of our results so far, an almost trivial consequence is the
following gap theorem. Such a theorem is useful, for instance, in
constructing bubble-trees. We state it here for convenience.
\begin{corollary}[Gap theorem]
There exists an $\epsilon_0=\epsilon_0(n,C_S)$ with the following
property. Assuming $(M,g,J)$ is an extremal K\"ahler orbifold
(possibly nonreduced) and that
\begin{eqnarray*}
\int_M|\Riem|^{\frac{n}{2}} \le\epsilon_0,
\end{eqnarray*}
then $(M,g)$ is flat.
\end{corollary}
\underline{\sl Pf}\\
\indent If the $\epsilon$-regularity theorem can be shown to hold,
namely that
\begin{eqnarray*}
\int_{B(o,r)}|\Riem|^{\frac{n}{2}}\le\epsilon_0 \quad\implies\quad
|\Riem(o)|\le{C}r^{-2}\left(\int_{B(o,r)}|\Riem|^{\frac{n}{2}}\right)^
\frac2n,
\end{eqnarray*}
then the result follows. If singularities are present, it is
possible that the Moser iteration technique will fail, due to its
reliance on integration by parts. However Theorems
\ref{LocalCurvatureBounds} and \ref{LocalCurvatureBoundsDim4},
combined with Sibner's lemma, Proposition \ref{L2Lemma}, will ensure
that residues will not crop up. Thus the Moser iteration goes
through, and we get our result. \qed

\section{Appendix: Local integral bounds for curvature at smooth points}

\subsection{Statement of the technical estimates} The following
propositions hold when $\supp{\phi}$ consists of smooth points, and
the real dimension is 4 or higher.
\begin{lemma}\label{NablaTNabla^2TEstimate}
Assume $0\le\phi$ and $\int|\Riem|^{\frac{n}{2}}$ has been chosen
small compared to $C_S$, $k$, and $l$.  With $C=C(k,l,C_S)$ we have
the two estimates
\begin{eqnarray}
\int\phi^l|\nabla{T}|^{k-2}|\nabla^2T|^2
&\le&C\int\phi^{l-2}|\nabla\phi|^2|\nabla{T}|^k\label{NablaSquared1}\\
&&+\,C\int\phi^l|\nabla{T}|^{k-2}|\triangle{T}|^2\nonumber\\
&&+\,C\int\phi^l|\nabla{T}|^{k-2}|T|^2|\Riem|^2\nonumber
\end{eqnarray}
\begin{eqnarray}
\int\phi^l|\nabla{T}|^{k-2}|\nabla^2T|^2
&\le&C\int\phi^{l-2}|\nabla\phi|^2|\nabla{T}|^k\label{NablaSquared2}\\
&&+\,C\int\phi^l|\nabla{T}|^{k-2}|\triangle{T}|^2\nonumber\\
&&+\,C\int\phi^l|\nabla{T}|^{k-1}|T||\nabla\Ric|\nonumber
\end{eqnarray}
\end{lemma}

\begin{lemma}\label{NablaTkGammaEstimate}
Assume $0\le\phi$ and $\int|\Riem|^{\frac{n}{2}}$ has been chosen
small compared to $C_S$, $k$, and $l$.  With $C=C(k,l,C_S)$ we have
the two estimates
\begin{eqnarray}
\left(\int\phi^{l\gamma}|\nabla{T}|^{k\gamma}\right)^\frac1\gamma
&\le&C\int\phi^{l-2}|\nabla\phi|^2|\nabla{T}|^k\label{kGamma1}\\
&&+\,C\int\phi^l|\nabla{T}|^{k-2}|\triangle{T}|^2\nonumber\\
&&+\,C\int\phi^l|\nabla{T}|^{k-2}|T|^2|\Riem|^2\nonumber
\end{eqnarray}
\begin{eqnarray}
\left(\int\phi^{l\gamma}|\nabla{T}|^{k\gamma}\right)^\frac1\gamma
&\le&C\int\phi^{l-2}|\nabla\phi|^2|\nabla{T}|^k\label{kGamma2}\\
&&+\,C\int\phi^l|\nabla{T}|^{k-2}|\triangle{T}|^2\nonumber\\
&&+\,C\int\phi^l|\nabla{T}|^{k-1}|T||\nabla\Ric|\nonumber
\end{eqnarray}
\end{lemma}

\begin{lemma}\label{NablaTkEstimate}
Assume $0\le\phi$, $|\nabla\phi|\le\frac2r$, and
$\int|\Riem|^{\frac{n}{2}}$ has been chosen small compared to $C_S$,
$k$, and $l$.  With $C=C(k,l,C_S)$ we have the two estimates
\begin{eqnarray}
\int\phi^l|\nabla{T}|^k
&\le&r^{-2}C\int\phi^{l-2}|\nabla{T}|^{k-2}|T|^2\label{kEst1}\\
&&+\,r^2C\int\phi^{l+2}|\nabla{T}|^{k-2}|\triangle{T}|^2\nonumber\\
&&+\,r^2C\int\phi^{l+2}|\nabla{T}|^{k-2}|T|^2|\Riem|^2\nonumber
\end{eqnarray}
\begin{eqnarray}
\int\phi^l|\nabla{T}|^k
&\le&r^{-2}C\int\phi^{l-2}|\nabla{T}|^{k-2}|T|^2\label{kEst2}\\
&&+\,r^2C\int\phi^{l+2}|\nabla{T}|^{k-2}|\triangle{T}|^2\nonumber\\
&&+\,r^2C\int\phi^{l+2}|\nabla{T}|^{k-1}|T||\nabla\Ric|\nonumber
\end{eqnarray}
\end{lemma}

\subsection{Proof of the technical estimates} We achieve the
estimates in a number of stages.  Our spaghetti-like argument
involves obtaining partial estimates for one quantity in order to
estimate a second, and using the second to get a better estimate for
the first, etc.  The steps involved are standard, so exposition is
kept to a minimum.

\noindent\underline{\sl Initial estimate for
$\int\phi^l|\nabla{T}|^{k-2}|\nabla^2T|^2$}
\begin{eqnarray*}
\int\phi^l|\nabla{T}|^{k-2}|\nabla^2T|^2
&=&-l\int\phi^{l-1}|\nabla{T}|^{k-2}\left<\nabla^2T,\nabla\phi\otimes
\nabla{T}\right>\\
&&-\,(k-2)\int\phi^l|\nabla{T}|^{k-2}|\nabla|\nabla{T}||^2\\
&&-\,\int\phi^l|\nabla{T}|^{k-2}\left<\triangle\nabla{T},\nabla{T}
\right>
\end{eqnarray*}
\begin{eqnarray*}
\int\phi^l|\nabla{T}|^{k-2}|\nabla^2T|^2
&\le&2l^2\int\phi^{l-2}|\nabla\phi|^2|\nabla{T}|^k
\,-\,2\int\phi^l|\nabla{T}|^{k-2}\left<\triangle\nabla{T},\nabla{T}
\right>
\end{eqnarray*}

\noindent\underline{\sl Estimate for
$-\int\phi^l|\nabla{T}|^{k-2}\left<\triangle\nabla{T},\nabla{T}\right>$}
Commutator formula:
\begin{eqnarray*}
\triangle\nabla{T} &=&\nabla\triangle{T} \,+\,\nabla(\Riem*T)
\,+\,\Riem*\nabla{T}.
\end{eqnarray*}
\begin{eqnarray*}
-\int\phi^l|\nabla{T}|^{k-2}\left<\triangle\nabla{T},\nabla{T}\right>
&=&-\int\phi^l|\nabla{T}|^{k-2}\left<\nabla\triangle{T},\nabla{T} \right>\\
&&-\,\int\phi^l|\nabla{T}|^{k-2}\left<\nabla(\Riem*T),\nabla{T}\right>\\
&&-\,\int\phi^l|\nabla{T}|^{k-2}\left<\Riem*\nabla{T},\nabla{T}\right>.
\end{eqnarray*}
We estimate the three terms individually. First term:
\begin{eqnarray*}
&&-\int\phi^l|\nabla{T}|^{k-2}\left<\nabla\triangle{T},\nabla{T} \right>\\
&&\quad=\;l\int\phi^{l-1}|\nabla{T}|^{k-2}\left<\nabla\phi\otimes
\triangle{T},\nabla{T}\right>\\
&&\quad\quad+\;(k-2)\int\phi^l|\nabla{T}|^{k-3}\left<\nabla|\nabla{T}|
\otimes\triangle{T},\nabla{T}\right>\\
&&\quad\quad+\;\int\phi^l|\nabla{T}|^{k-2}|\triangle{T}|^2
\end{eqnarray*}
\begin{eqnarray*}
&&-\int\phi^l|\nabla{T}|^{k-2}\left<\nabla\triangle{T},\nabla{T} \right>\\
&&\quad\le\;\frac{l^2}{2}\int\phi^{l-2}|\nabla\phi|^2|\nabla{T}|^k\\
&&\quad\quad+\;\mu\int\phi^l|\nabla{T}|^{k-2}|\nabla^2T|^2\\
&&\quad\quad+\;\left(\frac32+\frac{2(k-2)^2}{\mu}\right)\int\phi^l|
\nabla{T}|^{k-2}|\triangle{T}|^2
\end{eqnarray*}
Second term:
\begin{eqnarray*}
&&-\int\phi^l|\nabla{T}|^{k-2}\left<\nabla(\Riem*T),\nabla{T}\right>\\
&&\quad=\,l\int\phi^{l-1}|\nabla{T}|^{k-2}\left<\nabla\phi\otimes
(\Riem*T),\nabla{T}\right>\\
&&\quad\quad+\,(k-2)\int\phi^l|\nabla{T}|^{k-3}\left<\nabla|\nabla{T}|
\otimes(\Riem*T),\nabla{T}\right>\\
&&\quad\quad+\,\int\phi^l|\nabla{T}|^{k-2}\left<\Riem*T,\triangle{T}
\right>
\end{eqnarray*}
\begin{eqnarray*}
&&-\int\phi^l|\nabla{T}|^{k-2}\left<\nabla(\Riem*T),\nabla{T}\right>\\
&&\quad\le\;\frac{l^2}{2}\int\phi^{l-2}|\nabla\phi|^2|\nabla{T}|^k\\
&&\quad\quad+\;\mu\int\phi^l|\nabla{T}|^{k-2}|\nabla^2T|^2\\
&&\quad\quad+\,\left(1+\frac{2(k-2)^2}{\mu}\right)\int\phi^l|\nabla
{T}|^{k-2}|T|^2|\Riem|^2\\
&&\quad\quad+\,\frac12\int\phi^l|\nabla{T}|^{k-2}|\triangle{T}|^2
\end{eqnarray*}
Third term:
\begin{eqnarray*}
&&-\int\phi^l|\nabla{T}|^{k-2}\left<\Riem*\nabla{T},\nabla{T}\right>
\;\le\;\int\phi^l|\nabla{T}|^k|\Riem|
\end{eqnarray*}
Therefore
\begin{eqnarray*}
-\int\phi^l|\nabla{T}|^{k-2}\left<\triangle\nabla{T},\nabla{T}\right>
&\le&l^2\int\phi^{l-2}|\nabla\phi|^2|\nabla{T}|^k\\
&&+\,2\mu\int\phi^l|\nabla{T}|^{k-2}|\nabla^2T|^2\\
&&+\,\left(1+\frac{2(k-2)^2}{\mu}\right)\int\phi^l|\nabla{T}|^{k-2}|T|
^2|\Riem|^2\\
&&+\,\left(2+\frac{2(k-2)^2}{\mu}\right)\int\phi^l|\nabla{T}|^{k-2}|
\triangle{T}|^2\\
&&+\,\int\phi^l|\nabla{T}|^k|\Riem|.
\end{eqnarray*}

\noindent\underline{\sl Alternative estimate for
$-\int\phi^l|\nabla{T}|^{k-2}\left<\triangle\nabla{T},\nabla{T}\right>$}

Commutator formula:
\begin{eqnarray*}
\triangle\nabla{T} &=&\nabla\triangle{T} \,+\,\nabla\Ric*T
\,+\,\Riem*\nabla{T}.
\end{eqnarray*}
\begin{eqnarray*}
-\int\phi^l|\nabla{T}|^{k-2}\left<\triangle\nabla{T},\nabla{T}\right>
&=&-\int\phi^l|\nabla{T}|^{k-2}\left<\nabla\triangle{T},\nabla{T} \right>\\
&&-\,\int\phi^l|\nabla{T}|^{k-2}\left<\nabla\Ric*T,\nabla{T}\right>\\
&&-\,\int\phi^l|\nabla{T}|^{k-2}\left<\Riem*\nabla{T},\nabla{T}\right>.
\end{eqnarray*}
We deal with the terms individually. First term remains the same:
\begin{eqnarray*}
&&-\int\phi^l|\nabla{T}|^{k-2}\left<\nabla\triangle{T},\nabla{T} \right>\\
&&\quad\le\;\frac{l^2}{2}\int\phi^{l-2}|\nabla\phi|^2|\nabla{T}|^k\\
&&\quad\quad+\;\mu\int\phi^l|\nabla{T}|^{k-2}|\nabla^2T|^2\\
&&\quad\quad+\;\left(\frac32+\frac{2(k-2)^2}{\mu}\right)\int\phi^l|
\nabla{T}|^{k-2}|\triangle{T}|^2
\end{eqnarray*}
Second term:
\begin{eqnarray*}
-\int\phi^l|\nabla{T}|^{k-2}\left<\nabla\Ric*T,\nabla{T}\right>
&\le&\int\phi^l|\nabla{T}|^{k-1}|T||\nabla\Ric|
\end{eqnarray*}
The third term remains the same:
\begin{eqnarray*}
&&-\int\phi^l|\nabla{T}|^{k-2}\left<\Riem*\nabla{T},\nabla{T}\right>
\;\le\;\int\phi^l|\nabla{T}|^k|\Riem|
\end{eqnarray*}
Therefore
\begin{eqnarray*}
-\int\phi^l|\nabla{T}|^{k-2}\left<\triangle\nabla{T},\nabla{T}\right>
&\le&\frac{l^2}{2}\int\phi^{l-2}|\nabla\phi|^2|\nabla{T}|^k\\
&&+\,\mu\int\phi^l|\nabla{T}|^{k-2}|\nabla^2T|^2\\
&&+\,\left(\frac32+\frac{2(k-2)^2}{\mu}\right)\int\phi^l|\nabla{T}|^
{k-2}|\triangle{T}|^2\\
&&+\,\int\phi^l|\nabla{T}|^k|\Riem|\\
&&+\,\int\phi^l|\nabla{T}|^{k-1}|T||\nabla\Ric|
\end{eqnarray*}

\noindent\underline{\sl Two estimates for
$\int\phi^l|\nabla{T}|^{k-2}|\nabla^2T|^2$} First:
\begin{eqnarray*}
\int\phi^l|\nabla{T}|^{k-2}|\nabla^2T|^2
&\le&4l^2\int\phi^{l-2}|\nabla\phi|^2|\nabla{T}|^k\\
&&+\,\left(4+16(k-2)^2\right)\int\phi^l|\nabla{T}|^{k-2}|\triangle{T}| ^2\\
&&+\,\left(2+16(k-2)^2\right)\int\phi^l|\nabla{T}|^{k-2}|T|^2|\Riem|^2\\
&&+\,2\int\phi^l|\nabla{T}|^k|\Riem|
\end{eqnarray*}
Second:
\begin{eqnarray*}
\int\phi^l|\nabla{T}|^{k-2}|\nabla^2T|^2
&\le&3l^2\int\phi^{l-2}|\nabla\phi|^2|\nabla{T}|^k\\
&&+\,\left(3+8(k-2)^2\right)\int\phi^l|\nabla{T}|^{k-2}|\triangle{T}| ^2\\
&&+\,2\int\phi^l|\nabla{T}|^{k-1}|T||\nabla\Ric|\\
&&+\,2\int\phi^l|\nabla{T}|^k|\Riem|
\end{eqnarray*}

\noindent\underline{\sl Two Estimates for
$\left(\int\phi^{l\gamma}|\nabla{T}|^{k\gamma}\right)^\frac1\gamma$}\\

\noindent Sobolev Inequality:
\begin{eqnarray*}
\frac{1}{2C_S^2}\left(\int\phi^{l\gamma}|\nabla{T}|^{k\gamma}\right)^
\frac1\gamma &\le&l^2\int\phi^{l-2}|\nabla\phi|^2|\nabla{T}|^k
\,+\,k^2\int\phi^l|\nabla{T}|^{k-2}|\nabla^2T|^2
\end{eqnarray*}
Assume
\begin{eqnarray*}
\left(\int|\Riem|^{\frac{n}{2}}\right)^\frac2n\le\frac{1}{8k^2C_S^2}.
\end{eqnarray*}
First Estimate:
\begin{eqnarray*}
\left(\int\phi^{l\gamma}|\nabla{T}|^{k\gamma}\right)^\frac1\gamma
&\le&16k^2l^2C_S^2\int\phi^{l-2}|\nabla\phi|^2|\nabla{T}|^k\\
&&+\,32k^4C_S^2\int\phi^l|\nabla{T}|^{k-2}|\triangle{T}|^2\\
&&+\,32k^4C_S^2\int\phi^l|\nabla{T}|^{k-2}|T|^2|\Riem|^2\\
&&+\,4k^2C_S^2\int\phi^l|\nabla{T}|^k|\Riem|
\end{eqnarray*}
\begin{eqnarray*}
\left(\int\phi^{l\gamma}|\nabla{T}|^{k\gamma}\right)^\frac1\gamma
&\le&32k^2l^2C_S^2\int\phi^{l-2}|\nabla\phi|^2|\nabla{T}|^k\\
&&+\,64k^4C_S^2\int\phi^l|\nabla{T}|^{k-2}|\triangle{T}|^2\\
&&+\,64k^4C_S^2\int\phi^l|\nabla{T}|^{k-2}|T|^2|\Riem|^2
\end{eqnarray*}
Second Estimate:
\begin{eqnarray*}
\left(\int\phi^{l\gamma}|\nabla{T}|^{k\gamma}\right)^\frac1\gamma
&\le&8k^2l^2C_S^2\int\phi^{l-2}|\nabla\phi|^2|\nabla{T}|^k\\
&&+\,16k^4C_S^2\int\phi^l|\nabla{T}|^{k-2}|\triangle{T}|^2\\
&&+\,4k^2C_S^2\int\phi^l|\nabla{T}|^{k-1}|T||\nabla\Ric|\\
&&+\,4k^2C_S^2\int\phi^l|\nabla{T}|^k|\Riem|
\end{eqnarray*}
\begin{eqnarray*}
\left(\int\phi^{l\gamma}|\nabla{T}|^{k\gamma}\right)^\frac1\gamma
&\le&16k^2l^2C_S^2\int\phi^{l-2}|\nabla\phi|^2|\nabla{T}|^k\\
&&+\,32k^4C_S^2\int\phi^l|\nabla{T}|^{k-2}|\triangle{T}|^2\\
&&+\,8k^2C_S^2\int\phi^l|\nabla{T}|^{k-1}|T||\nabla\Ric|
\end{eqnarray*}

\noindent\underline{\sl Final estimate for $\int\phi^l|\nabla{T}|^
{k-2}|\nabla^2T|^2$}\\

\noindent First:
\begin{eqnarray*}
\int\phi^l|\nabla{T}|^{k-2}|\nabla^2T|^2
&\le&C(k,l,C_S)\int\phi^{l-2}|\nabla\phi|^2|\nabla{T}|^k\\
&&+\,C(k,l,C_S)\int\phi^l|\nabla{T}|^{k-2}|\triangle{T}|^2\\
&&+\,C(k,l,C_S)\int\phi^l|\nabla{T}|^{k-2}|T|^2|\Riem|^2\\
\end{eqnarray*}
Second:
\begin{eqnarray*}
\int\phi^l|\nabla{T}|^{k-2}|\nabla^2T|^2
&\le&C(k,l,C_S)\int\phi^{l-2}|\nabla\phi|^2|\nabla{T}|^k\\
&&+\,C(k,l,C_S)\int\phi^l|\nabla{T}|^{k-2}|\triangle{T}|^2\\
&&+\,C(k,l,C_S)\int\phi^l|\nabla{T}|^{k-1}|T||\nabla\Ric|\\
\end{eqnarray*}

\noindent\underline{\sl Initial estimate for
$\int\phi^l|\nabla{T}|^k$}\\
\begin{eqnarray*}
\int\phi^l|\nabla{T}|^k
&=&\int\phi^l|\nabla{T}|^{k-2}\left<\nabla{T},\nabla{T}\right>\\
&=&-l\int\phi^{l-1}|\nabla{T}|^{k-2}\left<\nabla{T},\nabla\phi\otimes
{T}\right>\\
&&-(k-2)\int\phi^l|\nabla{T}|^{k-3}\left<\nabla{T},\nabla|\nabla{T}|
\otimes{T}\right>\\
&&-\int\phi^l|\nabla{T}|^{k-2}\left<\triangle{T},T\right>
\end{eqnarray*}
\begin{eqnarray*}
\frac12\int\phi^l|\nabla{T}|^k
&\le&\frac{l^2}{2}\int\phi^{l-2}|\nabla\phi|^2|\nabla{T}|^{k-2}|T|^2\\
&&+r^2\frac{\mu}{2}\int\phi^{l+2}|\nabla{T}|^{k-2}|\nabla^2T|^2\\
&&+r^{-2}\frac{(k-2)^2}{2\mu}\int\phi^{l-2}|\nabla{T}|^{k-2}|T|^2\\
&&+r^{-2}\frac12\int\phi^{l-2}|\nabla{T}|^{k-2}|T|^2\\
&&+r^2\frac12\int\phi^{l+2}|\nabla{T}|^{k-2}|\triangle{T}|^2
\end{eqnarray*}
\begin{eqnarray*}
\int\phi^l|\nabla{T}|^k
&\le&\int\phi^{l-2}\left(l^2|\nabla\phi|^2+\mu^{-1}(k-2)^2r^{-2}+r^
{-2}\right)|\nabla{T}|^{k-2}|T|^2\\
&&+r^2\mu\int\phi^{l+2}|\nabla{T}|^{k-2}|\nabla^2T|^2\\
&&+r^2\int\phi^{l+2}|\nabla{T}|^{k-2}|\triangle{T}|^2
\end{eqnarray*}

\noindent \underline{\sl Two estimates for $\int\phi^l|\nabla{T}|^k$}\\

\noindent In this section we assume $\phi\le1$ and
$|\nabla\phi|\le\frac2r$.

\noindent First estimate:
\begin{eqnarray*}
\int\phi^l|\nabla{T}|^k
&\le&\int\phi^{l-2}\left(l^2|\nabla\phi|^2+\mu^{-1}(k-2)^2r^{-2}+r^
{-2}\right)|\nabla{T}|^{k-2}|T|^2\\
&&+\,r^2\mu4(l+2)^2\int\phi^l|\nabla\phi|^2|\nabla{T}|^k\\
&&+\,r^2\mu\left(4+16(k-2)^2\right)\int\phi^{l+2}|\nabla{T}|^{k-2}|
\triangle{T}|^2\\
&&+\,r^2\mu\left(2+16(k-2)^2\right)\int\phi^{l+2}|\nabla{T}|^{k-2}|T|
^2|\Riem|^2\\
&&+\,r^2\mu2\int\phi^{l+2}|\nabla{T}|^k|\Riem|\\
&&+r^2\int\phi^{l+2}|\nabla{T}|^{k-2}|\triangle{T}|^2
\end{eqnarray*}
\begin{eqnarray*}
\int\phi^l|\nabla{T}|^k
&\le&r^{-2}C(k,l)\int\phi^{l-2}|\nabla{T}|^{k-2}|T|^2\\
&&+\,r^2C(k,l)\int\phi^{l+2}|\nabla{T}|^{k-2}|\triangle{T}|^2\\
&&+\,r^2C(k,l)\int\phi^{l+2}|\nabla{T}|^{k-2}|T|^2|\Riem|^2\\
&&+\,r^2C(k,l)\left(\int\phi^{(l+2)\gamma}|\nabla{T}|^{k\gamma}\right)
^\frac1\gamma\left(\int|\Riem|^{\frac{n}{2}}\right)^\frac2n
\end{eqnarray*}
Continue:
\begin{eqnarray*}
\int\phi^l|\nabla{T}|^k
&\le&r^{-2}C(k,l)\int\phi^{l-2}|\nabla{T}|^{k-2}|T|^2\\
&&+\,r^2C(k,l)\int\phi^{l+2}|\nabla{T}|^{k-2}|\triangle{T}|^2\\
&&+\,r^2C(k,l)\int\phi^{l+2}|\nabla{T}|^{k-2}|T|^2|\Riem|^2\\
&&+\,C(k,l,C_S)\left(\int\phi^l|\nabla{T}|^k\right)\left(\int|\Riem|^
{\frac{n}{2}}\right)^\frac2n\\
&&+\,r^2C(k,l,C_S)\left(\int\phi^{l+2}|\nabla{T}|^{k-2}|\triangle{T}|
^2\right)\left(\int|\Riem|^{\frac{n}{2}}\right)^\frac2n\\
&&+\,r^2C(k,l,C_S)\left(\int\phi^{l+2}|\nabla{T}|^{k-2}|T|^2|\Riem|^2
\right)\left(\int|\Riem|^{\frac{n}{2}}\right)^\frac2n
\end{eqnarray*}
\begin{eqnarray*}
\int\phi^l|\nabla{T}|^k
&\le&r^{-2}C(k,l,C_S)\int\phi^{l-2}|\nabla{T}|^{k-2}|T|^2\\
&&+\,r^2C(k,l,C_S)\int\phi^{l+2}|\nabla{T}|^{k-2}|\triangle{T}|^2\\
&&+\,r^2C(k,l,C_S)\int\phi^{l+2}|\nabla{T}|^{k-2}|T|^2|\Riem|^2\\
\end{eqnarray*}

\noindent Second estimate:
\begin{eqnarray*}
\int\phi^l|\nabla{T}|^k
&\le&\int\phi^{l-2}\left(l^2|\nabla\phi|^2+\mu^{-1}(k-2)^2r^{-2}+r^
{-2}\right)|\nabla{T}|^{k-2}|T|^2\\
&&+\,r^2\mu3(l+2)^2\int\phi^l|\nabla\phi|^2|\nabla{T}|^k\\
&&+\,r^2\mu\left(3+8(k-2)^2\right)\int\phi^{l+2}|\nabla{T}|^{k-2}|
\triangle{T}|^2\\
&&+\,r^2\mu2\int\phi^{l+2}|\nabla{T}|^{k-1}|T||\nabla\Ric|\\
&&+\,r^2\mu2\int\phi^{l+2}|\nabla{T}|^k|\Riem|\\
&&+r^2\int\phi^{l+2}|\nabla{T}|^{k-2}|\triangle{T}|^2
\end{eqnarray*}
\begin{eqnarray*}
\int\phi^l|\nabla{T}|^k
&\le&r^{-2}C(k,l)\int\phi^{l-2}|\nabla{T}|^{k-2}|T|^2\\
&&+\,r^2C(k,l)\int\phi^{l+2}|\nabla{T}|^{k-2}|\triangle{T}|^2\\
&&+\,r^2C(k,l)\int\phi^{l+2}|\nabla{T}|^{k-1}|T||\nabla\Ric|\\
&&+\,r^2C(k,l)\left(\int\phi^{(l+2)\gamma}|\nabla{T}|^{k\gamma}\right)
^\frac1\gamma\left(\int|\Riem|^{\frac{n}{2}}\right)^\frac2n
\end{eqnarray*}
Continue:
\begin{eqnarray*}
\int\phi^l|\nabla{T}|^k
&\le&r^{-2}C(k,l)\int\phi^{l-2}|\nabla{T}|^{k-2}|T|^2\\
&&+\,r^2C(k,l)\int\phi^{l+2}|\nabla{T}|^{k-2}|\triangle{T}|^2\\
&&+\,r^2C(k,l)\int\phi^{l+2}|\nabla{T}|^{k-1}|T||\nabla\Ric|\\
&&+\,r^2C(k,l,C_S)\left(\int\phi^l|\nabla\phi|^2|\nabla{T}|^k\right)
\left(\int|\Riem|^{\frac{n}{2}}\right)^\frac2n\\
&&+\,r^2C(k,l,C_S)\left(\int\phi^{l+2}|\nabla{T}|^{k-2}|\triangle{T}|
^2\right)\left(\int|\Riem|^{\frac{n}{2}}\right)^\frac2n\\
&&+\,r^2C(k,l,C_S)\left(\int\phi^{l+2}|\nabla{T}|^{k-1}|T||\nabla\Ric|
\right)\left(\int|\Riem|^{\frac{n}{2}}\right)^\frac2n
\end{eqnarray*}
\begin{eqnarray*}
\int\phi^l|\nabla{T}|^k
&\le&r^{-2}C(k,l,C_S)\int\phi^{l-2}|\nabla{T}|^{k-2}|T|^2\\
&&+\,r^2C(k,l,C_S)\int\phi^{l+2}|\nabla{T}|^{k-2}|\triangle{T}|^2\\
&&+\,r^2C(k,l,C_S)\int\phi^{l+2}|\nabla{T}|^{k-1}|T||\nabla\Ric|\\
\end{eqnarray*}

\subsection{The induction argument in the smooth case}

In this section we assume the following:
\begin{hypothesis} \label{InductionAssumption}
Assume $a\ge\frac{n}{2}$, $n\ge4$.  There exist
$\epsilon_0=\epsilon_0(C_S,p,a,n)$ and $C=C(C_S,p,a,n)$ so that if
$\int_{B(o,r)}|\Riem|^{\frac{n}{2}}\le\epsilon_0$, then
\begin{eqnarray}
\left(\int_{B(o,r/2)}|\nabla^{p-1}\Riem|^a\right)^\frac1a
&\le&Cr^{-1-p+\frac{n}{a}}\left(\int_{B(o,r)}|\Riem|^{\frac{n}{2}}
\right)^\frac2n\\
\left(\int_{B(o,r/2)}|\nabla^p\Ric|^a\right)^\frac1a
&\le&Cr^{-2-p+\frac{n}{a}}\left(\int_{B(o,r)}|\Ric|^{\frac{n}{2}}
\right)^\frac2n\\
\left(\int_{B(o,r/2)}|\nabla^pX|^a\right)^\frac1a
&\le&Cr^{-3-p+\frac{n}{a}}\left(\int_{B(o,r)}|R|^{\frac{n}{2}}\right)^
\frac2n,
\end{eqnarray}
\end{hypothesis}
and prove:
\begin{theorem} \label{NablaPRiemTheorem}
Assume Hypothesis \ref{InductionAssumption}.  There exist
$\epsilon_0=\epsilon_0(C_S,p,a,n)$ and $C=C(C_S,p,a,n)$ so that if
$\int_{B(o,r)}|\Riem|^{\frac{n}{2}}\le\epsilon_0$, then
\begin{eqnarray}
\left(\int_{B(o,r/2)}|\nabla^p\Riem|^a\right)^\frac1a
&\le&Cr^{-2-p+\frac{n}{a}}\left(\int_{B(o,r)}|\Riem|^{\frac{n}{2}}
\right)^\frac2n\\
\left(\int_{B(o,r/2)}|\nabla^{p+1}\Ric|^a\right)^\frac1a
&\le&Cr^{-3-p+\frac{n}{a}}\left(\int_{B(o,r)}|\Ric|^{\frac{n}{2}}
\right)^\frac2n\\
\left(\int_{B(o,r/2)}|\nabla^{p+1}X|^a\right)^\frac1a
&\le&Cr^{-4-p+\frac{n}{a}}\left(\int_{B(o,r)}|R|^{\frac{n}{2}}\right)^
\frac2n.
\end{eqnarray}
\end{theorem}
Now Propositions \ref{GoodXEstimate}, \ref{GoodRicEstimate},
\ref{GoodRiemEstimate}, and \ref{GoodNablaRicAndXEstimate} together
show that Hypothesis \ref{InductionAssumption} holds in the case
$p=0$, so the conclusion is true for all $p\in\mathbb{N}$.

\noindent\underline{\sl Pf of Theorem \ref{NablaPRiemTheorem}}

\noindent First we use the commutator formula
(\ref{GeneralCommutator}) to get three estimates:
\begin{eqnarray}
&&|\triangle\nabla^{p-1}\Riem|^2
\;\le\;C\sum_{i=0}^{p-1}|\nabla^i\Riem|^2|\nabla^{p-1-i}\Riem|^2
\,+\,C|\nabla^{p+1}\Ric|^2\label{TriangleReductionRiem}\\
&&|\triangle\nabla^p\Ric|^2
\;\le\;C\sum_{i=0}^{p}|\nabla^i\Riem|^2|\nabla^{p-i}\Ric|^2
\,+\,C|\nabla^{p+1}X|^2\label{TriangleReductionRic}\\
&&|\triangle\nabla^pX|^2
\;\le\;C\sum_{i=0}^{p-1}|\nabla^i\Riem|^2|\nabla^{p-i}X|^2
\,+\,C\sum_{i=0}^p|\nabla^i\Ric|^2|\nabla^{p-i}X|^2.\label
{TriangleReductionX}
\end{eqnarray}
Notice that $\triangle\nabla^pX$ involves at most the $(p-1)^{th}$
derivative of $\Riem$.

\noindent\underline{\sl Step I: Estimating $|\nabla^{p+1}X|$ norms}

We deal with $\nabla^{p+1}X$ first. We can use Technical Lemma
\ref{NablaTkEstimate} to get
\begin{eqnarray*}
\int\phi^k|\nabla^{p+1}X|^k &\le& Cr^{-k}\int\phi^k|\nabla^pX|^k\\
&&\quad+\,Cr^2\left(\int\phi^{k\gamma}|\nabla^pX|^{k\gamma}\right)^
\frac1\gamma\left(\int_{\supp(\phi)}|\nabla\Ric|^{\frac{n}{2}}\right)^
{k\frac2n}\\
&&\quad+\,Cr^2\int\phi^k|\nabla^{p+1}X|^{k-2}|\triangle\nabla^pX|^2\\
&\le& Cr^{-k}\int\phi^k|\nabla^pX|^k\\
&&\quad+\,Cr^2\left(\int\phi^{k\gamma}|\nabla^pX|^{k\gamma}\right)^
\frac1\gamma\left(\int_{\supp(\phi)}|\nabla\Ric|^{\frac{n}{2}}\right)^
{k\frac2n}\\
&&\quad+\,Cr^2\sum_{i=0}^{p-1}\int\phi^k|\nabla^{p+1}X|^{k-2}|\nabla^i
\Riem|^2|\nabla^{p-i}X|^2\\
&&\quad+\,Cr^2\sum_{i=0}^p\int\phi^k|\nabla^{p+1}X|^{k-2}|\nabla^i
\Ric|^2|\nabla^{p-i}X|^2.
\end{eqnarray*}
The induction hypothesis yields estimates for each of the
$\nabla^i\Riem$, $\nabla^i\Ric$, and $\nabla^{p-i}X$ integral terms.
Then using H\"older's inequality and collecting terms will give us
\begin{eqnarray*}
\int\phi^k|\nabla^{p+1}X|^k &\le& Cr^{-k}\int\phi^k|\nabla^pX|^k.
\end{eqnarray*}
Using the induction hypothesis again yields
\begin{eqnarray}
\left(\int_{B(o,r/2)}|\nabla^{p+1}X|^k\right)^\frac1k &\le&
Cr^{-4-p+\frac{n}{k}}\left(\int_{B(o,r)}|R|^{\frac{n}{2}}\right)^
\frac2n.\label{FinalXResult}
\end{eqnarray}

\noindent\underline{\sl Step II: Estimating $|\nabla^{p+1}\Ric|$
norms} \\

\noindent Now it is necessary to bound $\int|\nabla^{p+1}\Ric|^k$.
Using \ref{NablaTkEstimate} again, we get
\begin{eqnarray*}
\int\phi^k|\nabla^{p+1}\Ric|^k &\le& Cr^{-2}\int\phi^{k-2}|\nabla^{p
+1}\Ric|^{k-2}|\nabla^p\Ric|^2\\
&&\quad+\,Cr^2\int\phi^{k+2}|\nabla^{p+1}\Ric|^{k-2}|\nabla^p\Ric||
\nabla\Ric|\\
&&\quad+\,Cr^2\int\phi^{k+2}|\nabla^{p+1}\Ric|^{k-2}|\triangle\nabla^p
\Ric|^2,
\end{eqnarray*}
which, with the induction hypothesis and H\"older's inequality,
becomes
\begin{eqnarray*}
\int\phi^k|\nabla^{p+1}\Ric|^k &\le& Cr^{-k}\int|\nabla^p\Ric|^k\\
&&+\,Cr^2\int\phi^{k+2}|\nabla^{p+1}\Ric|^{k-2}|\triangle\nabla^p\Ric|
^2.
\end{eqnarray*}
Use the formula for $\triangle\nabla^p\Ric$.  Noting that all of the
terms appearing in $\triangle\nabla^p\Ric$ are estimable except the
$|\nabla^p\Riem|$ term, we can use H\"older's inequality to actually
get
\begin{eqnarray}
\int\phi^k|\nabla^{p+1}\Ric|^k &\le&
Cr^{-k}\int_{\supp\phi}|\nabla^p\Ric|^k\label{NablaP1RicEstimate}\\
&&+Cr^2\int\phi^{k+2}|\nabla^{p+1}\Ric|^{k-2}|\nabla^p\Riem|^2|\Ric|.
\nonumber
\end{eqnarray}
We work with the final term:
\begin{eqnarray*}
&&Cr^2\int\phi^{k+2}|\nabla^{p+1}\Ric|^{k-2}|\nabla^p\Riem|^2|\Ric|\\
&&\quad\le\;Cr^2\left(\int\phi^{(k+2)\gamma}|\nabla^{p+1}\Ric|^{k
\gamma}\right)^{\frac{k-2}{k}\frac{n-2}{n}}\left(\int\phi^k|\nabla^p
\Riem|^k\right)^\frac2k\left(\int\phi^{\frac{2n}{k-2}}|\Ric|^{\frac{k}
{k-2}\frac{n}{2}}\right)^{\frac{k-2}{k}\frac2n}
\end{eqnarray*}
Now we must work with the
$\left(\int\phi^{(k+2)\gamma}|\nabla^{p+1}\Ric|^{k\gamma}\right)^
{\frac{k-2}{k}\frac{n-2}{n}}$ factor.  Using Technical Lemma
\ref{NablaTkGammaEstimate} we get
\begin{eqnarray*}
\left(\int\phi^{(k+2)\gamma}|\nabla^{p+1}\Ric|^{k\gamma}\right)^{\frac
{n-2}{n}}
&\le&C\int\phi^k|\nabla\phi|^2|\nabla^{p+1}\Ric|^k\\
&&+\,C\int\phi^k|\nabla^{p+1}\Ric|^{k-2}|\triangle\nabla^p\Ric|^2\\
&&+\,C\int\phi^k|\nabla^{p+1}\Ric|^{k-1}|\nabla^p\Ric||\nabla\Ric|.
\end{eqnarray*}
The integral norms of all quantities except $|\nabla^{p+1}\Ric|$ are
estimable, and we get
\begin{eqnarray*}
\left(\int\phi^{(k+2)\gamma}|\nabla^{p+1}\Ric|^{k\gamma}\right)^{\frac
{n-2}{n}} &\le&Cr^{-2}\int\phi^k|\nabla\phi|^2|\nabla^{p+1}\Ric|^k.
\end{eqnarray*}
Putting this back into (\ref{NablaP1RicEstimate}) and using the
induction hypothesis gives
\begin{eqnarray}
\int\phi^k|\nabla^{p+1}\Ric|^k &\le&
Cr^{-k}\int_{\supp\phi}|\nabla^p\Ric|^k\label {NablaP1RicBetterEstimate}\\
&&+\,Cr^{-k}\left(\int_{\supp\phi}|\Ric|^{\frac{n}{2}}\right)^
{\frac2n}\left(\int\phi^{k}|\nabla^p\Riem|^k\right).\nonumber
\end{eqnarray}

\noindent\underline{\sl Step III: Estimating $|\nabla^p\Riem|$
norms} \\
\noindent Using Technical Lemma \ref{NablaTkEstimate} we get
\begin{eqnarray*}
\int\phi^k|\nabla^p\Riem|^k
&\le&Cr^{-2}\int\phi^{k-2}|\nabla^p\Riem|^{k-2}|\nabla^{p-1}\Riem|^2\\
&&+Cr^2\int\phi^{k+2}|\nabla^p\Riem|^{k-2}|\nabla^{p-1}\Riem|^2|\Riem| ^2\\
&&+Cr^2\int\phi^{k+2}|\nabla^p\Riem|^{k-2}|\triangle\nabla^{p-1}\Riem|
^2.
\end{eqnarray*}
Applying H\"older's inequality and the induction hypothesis, we get
\begin{eqnarray*}
\int\phi^k|\nabla^p\Riem|^k &\le&Cr^{-k}\int|\nabla^{p-1}\Riem|^k
\,+\,Cr^k\int\phi^{k+2}|\triangle\nabla^{p-1}\Riem|^k.
\end{eqnarray*}
In fact the integral norms of all quantities appearing in
$\triangle\nabla^{p-1}\Riem$ are estimable by the induction
hypothesis.  We are left with only
\begin{eqnarray*}
\int\phi^k|\nabla^p\Riem|^k &\le&Cr^{-k}\int|\nabla^{p-1}\Riem|^k
\,+\,Cr^k\int\phi^{k+2}|\nabla^{p+1}\Ric|^k,
\end{eqnarray*}
which we can estimate with (\ref{NablaP1RicBetterEstimate}) to get
\begin{eqnarray*}
\int\phi^k|\nabla^p\Riem|^k &\le&Cr^{-k}\int|\nabla^{p-1}\Riem|^k.
\end{eqnarray*}
Finally the induction hypothesis gives
\begin{eqnarray*} \int_{B(o,r/2)}|\nabla^p\Riem|^k
&\le&Cr^{-pk}\int_{B(o,r)}|\Riem|^k.
\end{eqnarray*}
Finally also equation (\ref{NablaP1RicEstimate}) and the result of
Step II give the final estimate for $\int|\nabla^{p+1}\Ric|^k$. \qed

\end{document}